\documentclass{article}
\title{Locally symmetric spaces and K-theory of number fields}
\author{Thilo Kuessner}
\date{}
\usepackage[arrow, matrix, curve]{xy}
\usepackage{amsthm}
\usepackage{hyperref}
\input{amssym.tex}
\newtheorem{rem}{Remark}
\newtheorem{thm}{Theorem}
\newtheorem{lem}{Lemma}
\newtheorem{cor} {Corollary}

\newtheorem{pro} {Proposition}
\newtheorem{df} {Definition}

\newtheorem{obs}{Observation}

	\newcommand{\g}{\underline{g}}
	\newcommand{\kk}{\underline{k}}
	\newcommand{\p}{\underline{p}}
	\newcommand{\h}{\underline{h}}
	\newcommand{\tq}{\underline{t}}
\newcommand{\ff}{{\Bbb F}}

	\newenvironment{pf}{{\it Proof:}\quad}{\hfill$QED$}

\begin{document}
\maketitle
	\begin{abstract}
For a closed locally symmetric space $M=\Gamma\backslash G/K$ and
a representation $\rho:G\rightarrow GL\left(N,{\Bbb C}\right)$ we consider the push-forward of the fundamental class in $H_*\left(BGL\left(\overline{\Bbb Q}\right)\right)$ and a related invariant in $K_*\left(\overline{\Bbb Q}\right)\otimes{\Bbb Q}$. We discuss the nontriviality of this invariant and we generalize the construction to cusped locally symmetric spaces of ${\Bbb R}$-rank one.
%We describe a fundamental class construction which associates a group-homological or K-theoretic invariant to certain finite-volume locally symmetric spaces $\Gamma\backslash G/K$ and representations of $G$, and we discuss its nontriviality.
%For a locally symmetric space of noncompact type $M=\Gamma\backslash G/K$, either closed or cusped of rank 1, 
%We describe an invariant of flat bundles over locally symmetric spaces with 
%values in
%the K-theory of number fields and discuss the nontriviality and ${\Bbb Q}$-independence of its values.
		\end{abstract}

		\section{Introduction}

While elements in topological K-theory $K^{-*}\left(X\right)$ are,
by definition, represented by (virtual) vector bundles over the space $X$,
it is less evident what the topological meaning of elements in algebraic K-theory $K_*\left(
A\right)$, for a commutative ring $A$, may be. An approach, which can be found 
e.g.\ in the appendix of \cite{ka}, is to consider elements 
in $K_d\left(A\right)$ associated to a flat $GL\left(A\right)$-bundle over a
$d$-dimensional homology sphere
$M$. Namely, let 
%This is done via a fundamental class construction: let $\left[M\right]\in H_n\left(M;{\Bbb Z}\right)$ be the fundamental class, and $
$$\rho:\pi_1M\rightarrow GL\left(A\right)$$ 
be the monodromy representation
of the flat 
bundle,
then its plusification
$$\left(B\rho\right)^+:M^+\rightarrow BGL^+\left(A\right)$$ 
can, in view of $M^+\simeq{\Bbb S}^d$, be
considered as an element in algebraic K-theory $$ K_d\left(A\right):=\pi_dBGL^+\left(A\right).$$
It was proved by Hausmann and Vogel (see \cite{hv})
that for a finitely generated, commutative, unital ring $A$ and $d\ge 5$ or $d=3$, all elements in $K_d\left(A\right)$ arise from such a construction.

If the manifold
$M$ is not a homology sphere, but still possesses a fundamental class $\left[M\right]\in H_d\left(M;{\Bbb Q}\right)$, one
can still consider $$\left(B\rho\right)_*\left[M\right]\in H_d\left(BGL\left(A\right);{\Bbb Q}\right)$$ 
and can use a suitably defined projection (see Section 2.4) to the primitive part of the homology to obtain
$$\gamma\left(M\right)\in PH_d\left(BGL\left(A\right);{\Bbb Q}\right)\cong K_d\left(A\right)\otimes
{\Bbb Q}.$$

An interesting
special case, which has been studied by Dupont-Sah and others, is $K_3\left({\Bbb C}\right)$. By a theorem of Suslin (\cite{su}),
$K_3^{ind}\left({\Bbb C}\right)$ is, up to torsion, 
isomorphic to the Bloch group $B\left({\Bbb C}\right)$. On the other hand, each ideally triangulated hyperbolic 3-manifold yields, in a very natural way, 
an element in $B\left({\Bbb C}\right)$, the Bloch invariant. By \cite{ny}, this element does not depend on the chosen ideal triangulation.

A generalization to higher-dimensional hyperbolic manifolds was provided by Goncharov in \cite{gon}. 
To an odd-dimensional hyperbolic manifold $M^{2n-1}$
and the flat bundle coming from
a half-spinor representation he associates an element $\gamma\left(M\right)\in K_{2n-1}\left(\overline{\Bbb Q}\right)\otimes {\Bbb Q}$, 
and proves its nontriviality by showing that evaluation
of the Borel class yields (a fixed multiple of) the 
volume. 
	
It thus arises as a natural question, whether other locally symmetric spaces and different flat bundles give nontrivial elements in the K-theory of number fields (and eventually how much of algebraic K-theory in
odd degrees can be represented by locally symmetric spaces and representations
of
their fundamental groups). 
%Goncharov (\cite{gon}) proved that,
%for all (2m-1)-dimensional hyperbolic manifolds, the representation of $\pi_1$ coming
%		from the spin representation leads to a nontrivial element in 
%		$K_{2m-1}\left(\overline{\Bbb Q}\right)\otimes{\Bbb R}

In Section 2,
we generalize the argument from \cite{gon} to the extent that, for a compact locally
symmetric space $M^{2n-1}=\Gamma\backslash G/K$ of noncompact type
and a representation $\rho:G\rightarrow GL\left(N,{\Bbb C}\right)$, nontriviality
of the associated element $\gamma\left(M\right)\in K_{2n-1}\left(
\overline{\Bbb Q}\right)
\otimes {\Bbb Q}$ is (independently of $\Gamma$)
equivalent to nontriviality of the Borel class $\rho^*b_{2n-1}$.

It does, in general, not work to associate elements in algebraic K-theory to flat bundles
over manifolds with boundary.
Nonetheless we succeed in Section 4 to associate an element $\gamma\left(M\right)\in K_{2n-1}\left(
\overline{\Bbb Q}\right)
\otimes {\Bbb Q}$ to flat bundles over locally 
rank one symmetric spaces of finite volume. 
(\cite{gon} did this for hyperbolic manifolds and half-spinor representations, but implicitly assuming that $\partial M$ 
be connected.) 

Nontriviality of classes in $K_{2n-1}\left(
\overline{\Bbb Q}\right)
\otimes {\Bbb Q}\cong PH_{2n-1}\left(GL\left(N,\overline{\Bbb Q}\right);{\Bbb Q}\right)$ will be checked
by pairing with the Borel classes
%, that is the cohomology classes
$b_{2n-1}\in H_c^{2n-1}\left(GL\left(N,{\Bbb C}\right);{\Bbb R}\right)$.
% which 
%correspond to the generators of $H^*\left(U\left(N\right);{\Bbb Z}\right)\cong\Lambda_{\Bbb Z}\left(b_1,b_3,b_5,\ldots\right)$ under the van
%Est isomorphism $H^*\left(U\left(N\right);{\Bbb R}\right)\cong 
%H_c^{*}\left(GL\left(N,{\Bbb C}\right);{\Bbb R}\right)$, see Section 2.4 for a detailled description. 
The results of Section 2 (for closed manifolds) and Section 4 (for cusped manifolds) are subsumed as follows.\\
\\
{\bf Theorem.} {\em For each symmetric space $G/K$ of noncompact type and odd dimension $d=2n-1$, and to
each representation $\rho:G\rightarrow GL\left(N,{\Bbb C}\right)$ with $\rho^*b_{2n-1}\not=0$, there exists a constant $c_\rho\not=0$, such that the following holds.

If $M=\Gamma\backslash G/K$ is a finite-volume, orientable, locally symmetric space and either $M$ is compact or $rk\left(G/K\right)=1$,
%and $\rho\left(\pi_1\partial_iM\right)$ is unipotent for each boundary component $\partial_i M$, 
then
there is an element $$\gamma\left(M\right)\in K_{2n-1}\left( \overline{\Bbb Q}\right)\otimes{\Bbb Q}$$
such that 
%the Borel regulator $r_{2n-1}:K_{2n-1}\left(\overline{\Bbb Q}\right)\otimes{\Bbb Q}
%\rightarrow{\Bbb R}$ fulfills} 
application of the Borel class $b_{2n-1}$ yields}
$$<b_{2n-1},\gamma\left(M\right)>=c_\rho vol\left(M\right).$$
\noindent
%If $G/K$ is a symmetric space of noncompact type and a representation $\rho:G\rightarrow GL\left(N,{\Bbb C}\right)$ satisfies 
In particular, if $\rho^*b_{2n-1}\not=0$, then locally symmetric spaces $\Gamma\backslash G/K$ of ${\Bbb Q}$-independent volume
give ${\Bbb Q}
$-independent elements in $K_{2n-1}\left(
\overline{\Bbb Q}\right)
\otimes
{\Bbb Q}$. 

(In many cases           
one actually associates an element in $K_{2n-1}\left({\Bbb F}\right)\otimes
{\Bbb Q}$, for some number field ${\Bbb F}$, see \hyperref[Thm2]{Theorem \ref*{Thm2}} in Section 2.6.)\\

In Section 3, we work out the list of fundamental representations
$\rho:G\rightarrow GL\left(N,{\Bbb C}\right)$ for which $\rho^*b_{2n-1}\not=0$ holds true. It is easy to prove that $\rho^*b_{2n-1}\not=0
	$ is always true if $2n-1\equiv 3\ mod\ 4$ (and $\rho$ is not the trivial representation). We work out, for which fundamental representations
	$\rho^*b_{2n-1}\not=0$ holds if $2n-1\equiv 1\ mod\ 4$.\\
%(In \cite{gon} it was stated that the half-spinor representations would seem to be the
%only fundamental representations of $Spin\left(2n-1,1\right)$ that yield
%nontrivial invariants of odd-dimensional hyperbolic manifolds. This is
%however not the case. Indeed, if $2n-1=dim\left(M\right)\equiv 3\ mod\ 4$, then each nontrivial 
%representation of $Spin\left(2n-1,1\right)$ yields nontrivial invariants.)
%The proof uses only standard Lie algebra and representation
%theory. The result reads as follows.\\
\\
{\bf Theorem.} {\em 
The following is a complete list of irreducible symmetric spaces $G/K$ of noncompact type and fundamental representations $\rho:G\rightarrow GL\left(N,{\Bbb C}\right)$ with $\rho^*b_{2n-1}\not=0$ for $2n-1:=dim\left(G/K\right)$.\\
\begin{tabular}[ht]{|l|c|}
  \hline
Symmetric Space & Representation\\
  \hline\hline
 $SL_l\left({\Bbb R}\right)/SO_l, l\equiv 0,3,4,7\ mod\ 8$& any fundamental representation\\
 $SL_l\left({\Bbb C}\right)/SU_l, l\equiv 0\ mod\ 2$& any fundamental representation\\
 $SL_{2l}\left({\Bbb H}\right)/Sp_l, l\equiv 0\ mod\ 2$&  any fundamental representation\\
 $Spin_{p,q}/\left(Spin_p\times Spin_q\right), p,q\equiv 1\ mod\ 2, p\not\equiv q\ mod\ 4$&
any fundamental representation\\
 $Spin_{p,q}/\left(Spin_p\times Spin_q\right), p,q\equiv 1\ mod\ 2, p\equiv q\ mod\ 4$&
positive and negative half-spinor representation\\
% $Spin_l\left({\Bbb C}\right)/Spin_l, l\equiv 3\ mod\ 4$&
%the spinor representation and its conjugate\\
 $SO_l\left({\Bbb C}\right)/SO_l, l\equiv 3\ mod\ 4$& any fundamental representation\\
 $Sp_l\left({\Bbb C}\right)/Sp_l, l\equiv 1\ mod\ 4$& any fundamental representation\\
 $E_7\left({\Bbb C}\right)/E_7$& any fundamental representation\\

\hline
\end{tabular}\\
}\\

In this list, the only examples with $rank\left(G/K\right)=1$ are the hyperbolic spaces $H^d=Spin\left(d,1\right)/Spin\left(d\right)$ with $d$ odd. Thus in the noncompact case we only get invariants for hyperbolic manifolds. (In forthcoming work with Inkang Kim we will generalize the construction to ${\Bbb Q}$-rank 1 lattices in symmetric spaces of higher rank.)

For hyperbolic manifolds and half-spinor representations, the construction of $\gamma\left(M\right)$ is due to Goncharov. (Though the proof in 
\cite{gon} implicitly assumes that $\partial M$ be connected.) For hyperbolic 3-manifolds, another construction is due to Cisneros-Molina and Jones in \cite{cj}. (It was related in \cite{cj} to the construction of Neumann-Yang in \cite{ny}.) The latter has the advantage that the number of boundary components does not impose technical 
problems, contrary to the group-homological approach in \cite{gon}.

Our construction for closed locally symmetric spaces in Section 2 is a straightforward generalization of \cite{gon}. 
In the case of cusped locally symmetric spaces (with possibly more than one cusp) it would have
seemed more natural to stick to the approach of \cite{cj}, and in fact this approach generalizes to locally symmetric spaces in a completely straightforward way. However, we did 
not succeed to evaluate the Borel class (in order to discuss nontriviality of the obtained invariants) in this approach. On the other hand, Goncharov's approach, even in the case of only one cusp, uses very special properties of the spinor representation, which can not be generalized to other representations.

Therefore, our argument is sort of a mixture of both approaches. On the one hand it is closer in spirit to the arguments of \cite{gon} (but with the cuspidal completion in Section 4.2 memorizing the geometry of {\em distinct} cusps), 
on the other hand the argument in Section 4.3 uses arguments from \cite{cj} to circumvent the very special group-homological arguments that were applied in \cite{gon} in the special setting of the half-spinor representations.

Of course, it should be interesting to relate the different constructions more directly.

		\section{The closed case}

		The results of this section are fairly straightforward generalizations of the results in \cite{gon} from hyperbolic manifolds to locally symmetric spaces of noncompact type. 
(Similar constructions have also appeared in work of other authors, mainly for hyperbolic
$3$-manifolds.)
We will define a notion of representations with
nontrivial Borel class and will, mimicking the arguments in \cite{gon}, 
		show that representations with nontrivial Borel class
give rise to nontrivial elements in algebraic K-theory of number fields. The problem of constructing representations 
with nontrivial Borel class
will be tackled in Section 3.

\subsection{Preparations}
{\bf Classifying space.} For a group G, its classifying space $BG$ (with respect to the discrete topology on $G$) is the simplicial set $BG$ defined as follows:

- the k-simplices
of $BG$ are the k-tuples $\left(g_1,\ldots,g_k\right)$ with $g_1,\ldots,g_k\in G$,

- the operator $\partial:S_k\left(BG\right)\rightarrow S_{k-1}\left(BG\right)$ is defined by

$\partial\left(g_1,\ldots,g_k\right)=\left(g_2,\ldots,g_k\right)+\sum_{i=1}^{k-1}\left(-1\right)^i\left(g_1,\ldots,g_ig_{i+1},\ldots,g_k\right)+\left(-1\right)^k\left(g_1,\ldots,g_{k-1}\right)$,

- the degeneracy maps are defined by $s_j\left(g_1,\ldots,g_k\right)=\left(g_1,\ldots,g_j,1,g_{j+1},\ldots,g_k\right)$.

The simplicial chain complex of $BG$ will be denoted $C_*^{simp}\left(BG\right)$. The group homology with 
coefficients in a ring $R$ is $H_*\left(G;R\right)=H_*^{simp}\left(BG;R\right):=H_*\left(C_*^{simp}\left(BG\right)\otimes_{\Bbb Z}R,\partial\otimes 1\right)$. 

A homomorphism $\rho:\Gamma\rightarrow G$ induces 
%a simplicial map $B\rho:B\Gamma\rightarrow BG$ 
%and thus 
$\left(B\rho\right)_*:H_*\left(\Gamma;R\right)\rightarrow H_*\left(G;R\right)$.\\
\\
{\bf Straight simplices.} Let $M$ be a Riemannian manifold of nonpositive sectional curvature, $C_*\left(M\right)$ the chain complex of singular simplices.

Let $\pi:\widetilde{M}\rightarrow M$ be the universal covering.
Fix a point $x_0\in M$ and a lift $\tilde{x}_0\in\widetilde{M}$. Then $\Gamma:=\pi_1\left(M,x_0\right)$ acts isometrically by deck transformations 
on $\widetilde{M}$.

In a simply
                connected space of nonpositive sectional curvature each {\em ordered}
                $\left(k+1\right)$-tuple of vertices $\left(y_0,\ldots,y_k\right)$ determines a unique straight $k$-simplex
$str\left(y_0,\ldots,y_k\right)$. In particular, for $g_0,g_1,\ldots,g_k\in \Gamma=\pi_1\left(M,x_0\right)$ there is a unique straight simplex 
$$str\left(g_0\tilde{x}_0,g_1\tilde{x}_0,\ldots,g_k\tilde{x}_0\right)$$ in $\widetilde{M}$. A simplex $\sigma\in C_*\left(M\right)$ is said to be straight if some (hence any) lift $\tilde{\sigma}\in C_*\left(\widetilde{M}\right)$ with $\pi\tilde{\sigma}=\sigma$ is a straight simplex. 
%(All lifts of $\sigma$ belong to the $\Gamma$-orbit of some lift $\tilde{\sigma}$. Since $\Gamma$ acts by isometries, straightness of some lift $\tilde{\sigma}$ implies straightness of each other lift.)

Let $C_*^{str,x_0}\left(M\right)$
be the chain complex of straight simplices with
all vertices in $x_0$. There is a canonical chain map $$\Psi:C_*^{simp}\left(B\Gamma\right)\rightarrow C_*^{str,x_0}\left(M\right)$$ given by 
$$\Psi\left(g_1,\ldots,g_k\right):=\pi\left(str\left(\tilde{x}_0,g_1\tilde{x}_0,g_1g_2\tilde{x}_0,\ldots,g_1\ldots g_k\tilde{x}_0\right)\right).$$
%The homomorphism $\Psi$ is a chain map because
%$$\Psi\partial\left(g_1,\ldots,g_k\right)=\Psi\left(g_2,\ldots,g_k\right)+\sum_{i=1}^{k-1}\left(-1\right)^i\Psi\left(g_1,\ldots,g_ig_{i+1},\ldots,g_k\right)+\left(-1\right)^k \Psi\left(g_1,\ldots,g_{k-1}\right)$$
%$$=\pi\left(str\left(\tilde{x}_0,g_2\tilde{x}_0,\ldots,g_2\ldots g_k\tilde{x}_0\right)\right)+\sum_{i=1}^{k-1}\left(-1\right)^i\pi\left(str\left(\tilde{x}_0,\ldots,g_1\ldots g_{i-1}\tilde{x}_0,g_1\ldots g_ig_{i+1}\tilde{x}_0,\ldots,g_1\ldots g_k\tilde{x}_0\right)\right)+$$
%$$\left(-1\right)^k\pi\left(str\left(\tilde{x}_0,g_1\tilde{x}_0,\ldots,g_1\ldots g_{k-1}\tilde{x}_0\right)\right)
%=\pi\left(str\left(g_1\tilde{x}_0,g_1g_2\tilde{x}_0,\ldots,g_1g_2\ldots g_k\tilde{x}_0\right)\right)+$$
%$$\sum_{i=1}^{k-1}\left(-1\right)^i\pi\left(str\left(\tilde{x}_0,\ldots,g_1\ldots g_{i-1}\tilde{x}_0,g_1\ldots g_ig_{i+1}\tilde{x}_0,\ldots,g_1\ldots g_k\tilde{x}_0\right)\right)+\left(-1\right)^k\pi\left(str\left(\tilde{x}_0,g_1\tilde{x}_0,\ldots,g_1\ldots g_{k-1}\tilde{x}_0\right)\right)$$
%$$=\pi\left(\partial str\left(\tilde{x}_0,g_1\tilde{x}_0,\ldots,g_1\ldots g_k\tilde{x}_0\right)\right)=
%\partial \Psi\left(g_1,\ldots,g_k\right),$$ where we have used that 
%$\pi\left(str\left(\tilde{x}_0,g_2\tilde{x}_0,\ldots,g_2\ldots g_k\tilde{x}_0\right)\right)=\pi\left(str\left(g_1\tilde{x}_0,g_1g_2\tilde{x}_0,\ldots,g_1g_2\ldots g_k\tilde{x}_0\right)\right)$ for each deck transformation $g_1\in \Gamma$.

Let $w_0,\ldots,w_k$ be the vertices of the standard simplex $\Delta^k$. For $j=0,\ldots,k$ let $\gamma_j\subset\Delta^k$ be the sub-1-simplex with
$\partial\gamma_j=w_{j}
-w_{j-1}$ for $j=1,\ldots,k$. Then there is a homomorphism $$\Phi: C_*^{str,x_0}\left(M\right)\rightarrow C_*^{simp}\left(B\Gamma\right)$$
defined by $\Phi\left(\sigma\right)=\left(g_1,\ldots,g_k\right)$, where 
%$\sigma\in
%C_k^{str,x_0}\left(M\right)$ is a continuous map $\sigma:\Delta^k\rightarrow M$ with $\sigma\left(w_j\right)=x_0$ for $j=0,\ldots,k$,
%and 
$g_j\in\Gamma=\pi_1\left(M,x_0\right)$ is the homotopy class (rel.\ vertices) of $\sigma\mid_{\gamma_j}$ for $j=1,\ldots,k$.

Clearly $\Phi\left(\pi\left(str\left(\tilde{x}_0,g_1\tilde{x}_0,g_1g_2\tilde{x}_0,\ldots,g_1\ldots g_k\tilde{x}_0\right)\right)\right)=\left(g_1,\ldots,g_k\right)$, thus $\Phi\Psi=id$. On the other hand, 
%each $\sigma\in
%C_k^{str,x_0}\left(M\right)$ is a continuous map 
a straight simplex $\sigma:\Delta^k\rightarrow M$ with all vertices in $x_0$ is uniquely determined by the homotopy classes (rel.\ vertices) of
$g_j=\left[\sigma\mid_{\gamma_j}\right]$ for $j=1,\ldots,k$, because its lift to $\widetilde{M}$ must be in the $\Gamma$-orbit of 
$str\left(\tilde{x}_0,g_1\tilde{x}_0,g_1g_2\tilde{x}_0,\ldots,g_1\ldots g_k\tilde{x}_0\right)$. Thus $\Psi\Phi=id$. This shows that $\Psi$ and $\Phi$ are chain isomorphisms, inverse to each other.\\
%$\partial\gamma_j=w_{j}
%-w_{j-1}$ for $j=1,\ldots,k$. Indeed, if $g_j:=\left[\sigma\mid_{\gamma_j}\right]\in
%\Gamma=\pi_1\left(M,x_0\right)$, then $\sigma=\pi\left(str\left(\tilde{x}_0,g_1\tilde{x}_0,\ldots,g_k\tilde{x}_0\right)\right)$.
%Thus
%$$C_*^{str,x_0}\left(M\right)\cong C_*^{simp}\left(B\Gamma\right),$$
%where the bijection maps
%$\sigma$ to $\left(g_1,\ldots,g_k\right)$. It follows from the definition of the boundary operator on $B\Gamma$ that this bijection is a chain map, thus an isomorphism of chain complexes.
\\
{\bf Eilenberg-MacLane map.} 
%Let $M$ be a Riemannian
%manifold of nonpositive sectional curvature. 
%As usual we denote by $C_*\left(M\right):=C_*^{sing}\left(M\right)$ the singular chain complex. Moreover we denote 
Let $C_*^{x_0}\left(M \right)\subset C_*\left(M\right)$ 
be the subcomplex generated by singular simplices with all vertices in $x_0$.
The inclusions $$
%C_*^{simp}\left(B\Gamma\right)\simeq 
C_*^{str,x_0}\left(M\right)\subset C_*^{x_0}\left(M\right)\subset C_*\left(M\right)$$
are chain homotopy equivalences\footnotemark\footnotetext[1]{Pictorially
the chain homotopy inverse $$str:C_*\left(M\right)\rightarrow C_*^{str,x_0}\left(M\right)$$ of the inclusion
$C_*^{str,x_0}\left(M\right)\subset  C_*\left(M\right)$
first homotopes all vertices of a given cycle into $x_0$ and then straightens the so-obtained cycle (by induction on dimension of subsimplices, depending on the 
given order of vertices)
as in in \cite[Lemma C.4.3]{bp}.}.
For the first inclusion this is proved (for arbitrary aspherical manifolds, but with an isomorphic image
of $C_*^{simp}\left(B\Gamma\right)$ instead of the in this generality not defined $C_*^{str,x_0}\left(M\right)$) in \cite[Theorem 1a]{em}.
For the second inclusion it is proved in Paragraph 31 of \cite{eil}. 

The composition of the chain isomorphism $\Psi:C_*^{simp}\left(B\Gamma\right)\rightarrow C_*^{str,x_0}\left(M\right)$ with the inclusion 
$C_*^{str,x_0}\left(M\right)\rightarrow C_*\left(M\right)$ is thus a chain homotopy equivalence
$$C_*^{simp}\left(B\Gamma\right)\rightarrow C_*
\left(M\right),$$
%which we denote the Eilenberg-MacLane map.
the induced isomorphism
$$EM: H_*^{simp}\left(B\Gamma;{\Bbb Z}\right)\rightarrow H_*\left(M;{\Bbb Z}\right)$$
will be called the Eilenberg-MacLane map.
The chain homotopy inverse is given by the composition of $str$ with the chain isomorphism $\Phi$, thus $$EM^{-1}=\Phi_*\circ str_*.$$
The geometric realization $\mid B\Gamma\mid$ of $B\Gamma$ in the sense of \cite{mil} is aspherical by
% Then $\mid B\Gamma\mid$ is a $K\left(\Gamma,1\right)$, that is $\pi_1\mid B\Gamma\mid\cong \Gamma$ and $\pi_n\mid B\Gamma\mid=0$ for $n\ge 2$, see 
\cite[p.128]{may}.
Given a manifold $M$ and an isomorphism $I:\pi_1M\cong\Gamma$, there is an up to homotopy unique continuous mapping $h^M:M\rightarrow
\mid B\Gamma\mid$ which induces $I$ on $\pi_1$, see \cite[p.177]{may}. The map $h^M$ (rather its homotopy class) is called the classifying map for $\pi_1M$.
%, the set of homotopy classes of maps from maps from $M$ to $\mid B\Gamma\mid$ is given by $\left[M,\mid B\Gamma\right]\cong Hom\left(\pi_1M,\Gamma\right)$, where the bijection maps each homotopy class to the induced homomorphism of fundamental groups (with appropriate base points)
If $M$ is aspherical and has the homotopy type of a CW-complex 
then $h^M$
is a homotopy equivalence,
and
% $h^M:M\rightarrow \mid B\Gamma\mid$ such that
$h^M_*:H_*\left(M;{\Bbb Z}\right)\rightarrow H_*\left(\mid B\Gamma\mid;{\Bbb Z}\right)$ is the composition of $EM^{-1}$ with the isomorphism 
$i_*:H_*^{simp}\left(B\Gamma;{\Bbb Z}\right)\rightarrow
H_*\left(\mid B\Gamma\mid;{\Bbb Z}\right)$ that is induced by the inclusion $i$
of the simplicial into the singular chain complex.
%(If $M$ is a smooth manifold, then one has a triangulation $M=\tau_1
%\cup\ldots\cup\tau_r$ and one can explicitly realize $h^M:M\rightarrow 
%\mid B\Gamma\mid$ by mapping the simplex $\tau_i$ 
%to the simplex $\Phi\left(str\left(\tau_i\right)\right)$ in $B\Gamma$.)

		\subsection{Construction of elements in algebraic K-theory}

Throughout this paper, a {\em ring} $A$ will mean a {\em commutative ring with unit}. In all applications $A$ will be a subring (with unit) of the ring of complex numbers: $A\subset {\bf C}$.

One defines $GL\left(A\right)$ as the increasing union $GL\left(A\right)=\cup_{N\in{\bf N}}GL\left(N,A\right),$
where $GL\left(N,A\right)$ is considered as a subgroup of $GL\left(N+1,A\right)$ via the canonical embedding as $N\times N$-block matrices with complementary $1\times 1$-block having entry $1$.
We consider the simplicial set
$BGL\left(A\right)$ as defined in Section 2.1, and $\mid BGL\left(A\right)\mid$ its geometrical realisation. 
%(In the language of algebraic topologists, 
%$\mid BGL\left(A\right)\mid$ is the classifying space for $GL\left(A\right)^\delta$,
%that is for the group $GL\left(A\right)$
%with the discrete topology. Thus $\pi_1 \mid BGL\left(A\right)\mid=GL\left(A\right)$.)

%Assume that $M$ is a closed, orientable, connected $n$-manifold with $\Gamma:=\pi_1M$.
%		Assume that we are given a commutative ring $A$ with unit and a representation
A representation		$\rho:\Gamma
\rightarrow GL\left(A\right)$ induces 
%a simplicial map $B\rho:B\Gamma\rightarrow BGL\left(A\right)$ and thus 
a continuous map
$$\mid B\rho\mid: \mid B\Gamma\mid \rightarrow \mid BGL\left(A\right)\mid.$$ 
%If $\Gamma=\pi_1M$, then composition with the classifying map $h^M:M\rightarrow \mid B\Gamma\mid$ yields
%a continuous map $$\mid B\rho\mid \circ h^M:M\rightarrow 
\begin{df} Let $M$ be a topological space with $\Gamma:=\pi_1\left(M,x_0\right),x_0\in M$, let $A$ be a commutative ring with unit and let $\rho:\Gamma
\rightarrow GL\left(A\right)$ be a homomorphism. Then
%, for $d\in{\Bbb N}$, 
we define 
$$\left(H\rho\right)_*:H_d\left(M;{\Bbb Q}\right)\rightarrow H_d^{simp}\left(BGL\left(A\right);{\Bbb Q}\right)$$
as the composition of $\mid B\rho\mid$ with the classifying map $h^M:M\rightarrow \mid B\Gamma\mid$:
$$
\begin{xy}
\xymatrix{H_d\left(M;{\Bbb Q}\right)\ar[r]
^{h^M_*}
& H_d\left(\mid B\Gamma\mid;{\Bbb Q}\right)\ar[r]
^{\mid B\rho\mid_*}
& H_d\left(\mid BGL\left(A\right)\mid;{\Bbb Q}\right)&
\cong H_d^{simp}\left(BGL\left(A\right);{\Bbb Q}\right)}\end{xy}.$$
%where $h^M:M\rightarrow \mid B\Gamma\mid$ is the classifying map.
%where the last arrow is the inverse of the isomorphism induced by the inclusion $C_*^{simp}\left(BGL\left(A\right)\right)\rightarrow C_*\left(\mid BGL\left(A\right)\mid;{\Bbb Q}\right)$.\end{df}
\end{df}
(We will use without mention that inclusion $i:C_*^{simp}\left(BGL\left(A\right)\right)\rightarrow C_*\left(\mid BGL\left(A \right)\mid\right)$ induces an isomorphism $i_*:H_*^{simp}\left(BGL\left(A\right)\right)\rightarrow H_*\left(\mid BGL\left(A \right)
\mid\right)$, see \cite[Lemma 5]{mil}.)\\
%The continuous map $\mid B\rho\mid \circ h^M$ induces a 
%	homomorphism $$\left(\mid B\rho\mid \circ h^M\right)_d:H_d\left(M;{\Bbb Q}\right)\rightarrow H_d\left(\mid BGL\left(A\right)\mid;{\Bbb Q}\right)\cong H_d^{simp}\left(BGL\left(A\right);{\Bbb Q}\right).$$ 
%Since 
%If $M$ is a closed, oriented, connected $d$-manifold, then we have a fundamental class 

If $M$ is
%not necessarily
%a homology sphere, but
a closed, orientable, connected d-manifold, and the ring $A$ satisfies mild assumptions (see Section 2.5), e.g.\ for $A={\overline{\Bbb Q}}$,
% or a (not totally real) number field or its ring of integers,
then we will now explain how to construct an element
%in $K_d\left(A\right)$ but rather in
in $K_d\left( A\right)\otimes{\Bbb Q}$.

Let $ \left[M\right]\in H_d\left(M;{\Bbb Q}\right)$ be the fundamental class.
%, which is the image of a generator 
%of $H_d\left(M;{\Bbb Z}\right)$ under the change-of-rings homomorphism associated to the
%inclusion ${\Bbb Z}\rightarrow{\Bbb Q}$.
%We may consider the image of the fundamental class $\left[M\right]\in H_d\left(M;{\Bbb Q}\right)$ 
Consider
		$$\left(H\rho\right)_*\left[M\right]\in 
H_d\left(\mid BGL\left(A\right)\mid;{\Bbb Q}\right)\cong H_d\left(\mid BGL\left(A\right)\mid^+;{\Bbb Q}\right).$$
By the Milnor-Moore Theorem, the Hurewicz homomorphism $$K_d\left(A\right):=\pi_d \left(\mid BGL\left(A\right)\mid^+\right)\rightarrow  H_d
		\mid \left(BGL\left(A\right)\mid^+;{\Bbb Z}\right)$$ 
gives, after tensoring
with ${\Bbb Q}$, an injective
homomorphism $$K_d\left(A\right)
\otimes{\Bbb Q}=\pi_d \left(\mid BGL\left(A\right)\mid^+\right)\otimes{\Bbb Q}
\rightarrow  H_d                             \left(\mid BGL\left(A\right)\mid^+;{\Bbb Q}\right).$$
%\cong H_d^{simp} \left(BGL\left(A\right);{\Bbb Q}\right).$$

Its image, again by the Milnor-Moore theorem, 
is the subgroup 
$PH_d                             \left(\mid BGL\left(A\right)\mid^+;{\Bbb Q}\right)$
of primitive elements, which 
%we denote by 
%$PH_d                             \left(\mid BGL\left(A\right)\mid^+;{\Bbb Q}\right)$.
we will henceforth identify with $K_d\left(A\right)
\otimes{\Bbb Q}$.

%One of the defining properties of Quillen's plus construction (see \cite{ro}) is that
By Quillen 
%, inclusion induces an isomorphism $Q_*:H_*\left(\mid BGL\left(A\right)\mid;{\Bbb Q}\right)\rightarrow H_*\left(\mid BGL\left(A\right)\mid^+;{\Bbb Q}\right)$, 
(compare \cite[Section 9.1]{bu}), inclusion induces an isomorphism $$Q_*:PH_*\left(\mid BGL\left(A\right)\mid;{\Bbb Q}\right)\rightarrow PH_*\left(\mid BGL\left(A\right)\mid^+;{\Bbb Q}\right)\cong K_*\left(A\right)\otimes{\Bbb Q}.$$

(If $d$ is even and $A$ is a ring of integers in any number field, then $PH_d^{simp}\left(BGL\left(A\right);{\Bbb Q}\right)
=0$, cf.\ \cite[Theorem 9.9]{bu}. Therefore one is only interested in the case $d=2n-1$.)
% For some reasons to be seen later we will not tensor with ${\Bbb Q}$ but with ${\Bbb R}\left(-m\right):=\left(2\pi i\right)^{-m}{\Bbb R}$. Therefore we
%consider the injective homomorphism
%$$I_{2m-1}:K_{2m-1}\left(A\right)
%\otimes{\Bbb R}\left(-m\right)\rightarrow H_{2m-1}\left(BGL\left(A\right),{\Bbb R}\left(-m\right)\right).$$

%We note that there is a canonical 
Whenever we have a fixed projection $$pr_*: H_*^{simp}\left(BGL\left(A\right);{\Bbb Q}\right)\rightarrow
PH_*^{simp}\left(BGL\left(A\right);{\Bbb Q}\right)\cong PH_*\left(\mid BGL\left(A\right)\mid;{\Bbb Q}\right),$$
%	Namely, the kernel of the Borel regulator $r_n:H_n\left(BGL\left(A\right);{\Bbb Q}\right)\rightarrow {\Bbb R}$ (defined by pairing with the Borel
%class $b_n$, see section 2.3.) is complementary to $P_n\left(BGL\left(A\right);{\Bbb Q}\right)$,
%thus we may define $pr_n$ as the projection along the kernel of the Borel
%regulator.
we can define an element $\gamma\left(M\right)\in
K_d\left(A\right)\otimes{\Bbb Q}$ as
$$\gamma\left(M\right):=Q_* pr_*\left(H\rho\right)_*\left[M\right].$$
%where $Pi_d:PH_*^{simp}\left(BGL\left(A\right);{\Bbb Q}\right)\rightarrow PH_*\left(\mid BGL\left(A\right)\mid;{\Bbb Q}\right)$ is the restriction of 
%$i_d:H_*^{simp}\left(BGL\left(A\right);{\Bbb Q}\right)\rightarrow H_*\left(\mid BGL\left(A\right)\mid;{\Bbb Q}\right)$
%to the subgroup of primitive elements.

In Section 2.5 we are going to show that
e.g.\ for $A={\overline{\Bbb Q}}$ (and also for many other rings) the projection $pr_*$ can be chosen such that the 
evaluations of the Borel class on
$h$ and $pr_*\left(h\right)$ agree for all $h\in H_*^{simp}\left(BGL\left({\overline{\Bbb Q}}\right);{\Bbb Q}\right)$. In particular, to check nontriviality of $\gamma\left(M\right)$ it will then suffice to apply the Borel 
class to $\left(H\rho\right)_*\left[M\right]$.\\
%Of course, the so-constructed element $P\left(B\rho\right)_*\left[M\right]$ depends on $P$. However, the value of the Borel regulator of $P\left(B\rho\right)_*\left[M\right]$ does not depend on $P$, because the Borel regulator vanishes on
%	decomposable elements (i.e.\ is well-defined on K-theory). Since, by Borel's Theorem, the K-theory of number rings is determined by values of Borel regulators, we may thus neglect the dependence on the projection $P$.

If $M$ is a (compact, orientable) manifold with {\em nonempty} boundary, then
there is no general construction of an element in algebraic K-theory. However, we will show in Section 4 that for 
finite-volume locally rank one symmetric spaces one can generalize the above construction and again construct an 
invariant $\gamma\left(M\right)\in K_d\left(\overline{\Bbb Q}\right)\otimes{\Bbb Q}$.

		\subsection{The volume class in $H_c^d\left(Isom\left(\widetilde{M}\right)\right)$}

	{\bf Volume class.}
		For a Lie group $G$, let $C_c\left(G^{*+1},{\Bbb R}\right)$ be the {\em continuous} $G$-invariant mappings from $G^{*+1}$ to ${\Bbb R}$, 
$\delta$ the usual coboundary operator and 
		$H^*_c\left(G;{\Bbb R}\right)$
		the cohomology of $\left(C_c\left(G^{*+1},{\Bbb R}\right)^G,\delta\right)$.
There is a comparison map $comp:H_c^*\left(G;{\Bbb R}\right)\rightarrow
		H^*\left(G;{\Bbb R}\right)$
defined by the cochain map
$$comp\left(f\right)\left(g_1,\ldots,g_k\right):=f\left(1,g_1,g_1g_2,\ldots,g_1g_2\ldots g_k\right).$$
% for $f\in C_c^k\left(G\right)$. In particular, elements of $H_c^*\left(G;{\Bbb R}\right)$ can be evaluated on $H_*^{simp}\left(BG;{\Bbb R}\right)$.

\begin{rem} For $f\in H_c^*\left(G;{\Bbb R}\right)$ and $c\in H_*\left(G;{\Bbb R}\right)$, we will denote
$$<f,c>=comp\left(f\right)\left(c\right).$$\end{rem}

  Let $\widetilde{M}=G/K$ be a symmetric space of noncompact type. The Riemannian metric and in particular the volume form are given via the
Killing form and are thus canonical.
%, with isometry group $G$. 
It is well-known (\cite[Chapter V, Theorem 3.1]{hel}) that $\widetilde{M}$ has nonpositive sectional curvature. One can assume that $G$ is semisimple and acts by orientation-preserving isometries on $\widetilde{M}$.\\

Fix an arbitrary
point $\tilde{x}\in
                \widetilde{M}=G/K$.
		The volume class $$v_d\in H_c^d\left(G;{\Bbb R}\right)$$ is defined as follows. 
%In
%a simply
%                connected space of nonpositive sectional curvature each {\em ordered}
%                $k+1$-tuple of vertices $\left(p_0,\ldots,p_k\right)$ determines a unique straight $k$-simplex $str\left(p_0,\ldots,p_k\right)$.
We define a simplicial cochain
		$c\nu_d\in C^d_{simp}\left(BG\right)$ by $$c\nu_d\left(g_1,\ldots,g_d\right) =
		algvol\left(str\left(\tilde{x},g_1\tilde{x},\ldots,g_1\ldots g_d\tilde{x}\right)\right):=\int_{str\left(\tilde{x},g_1\tilde{x},\ldots,g_d\tilde{x}\right)}dvol,$$
%as defined in \cite[p.107]{bp}, p.107, 
that is the signed volume $algvol$ (see \cite[p.107]{bp})
                of the straight simplex with vertices $\tilde{x},g_1\tilde{x},\ldots,g_1\ldots g_d\tilde{x} $.
%Since
%		the volume form is invariant under isometries (and isometries map straight simplices to
%		straight simplices), this is a $G$-invariant mapping. 
From Stokes' Theorem and $algvol\left(str\left(g_1\tilde{x}, g_1g_2\tilde{x},\ldots,g_1\ldots g_d\tilde{x}\right)\right)=algvol\left(str\left(\tilde{x},g_2\tilde{x},\ldots,g_2\ldots,g_d\tilde{x}\right)\right)$ one can conclude
%implies that $c\nu_d$ 
%is a cocycle: 
$$\delta c\nu_d\left(g_1,\ldots,g_{d+1}\right)=
%c\nu_d\left(\left(g_2,\ldots,g_{d+1}\right)
%+\sum_{i=1}^{d}\left(-1\right)^i\left(g_1,\ldots,g_ig_{i+1},\ldots,g_{d+1}\right)+\left(-1\right)^{d+1}\left(g_1,\ldots,g_{d}\right)\right)=$$
%$$
%algvol\left(str\left(g_1\tilde{x}, g_1g_2\tilde{x},\ldots,g_1\ldots g_{d+1}\tilde{x}\right)\right)+\sum_{i=1}^{d}\left(-1\right)^i
%algvol\left(str\left(\tilde{x}, \ldots,\widehat{g_1\ldots g_{i}\tilde{x}},\tilde{x},\ldots,g_1\ldots g_{d+1}\tilde{x}\right)\right)+$$
%$$+\left(-1\right)^{d+1} algvol\left(str\left(\tilde{x}, \ldots,g_1\ldots g_{d}\tilde{x}\right)\right)
%= \sum_{i=0}^{d+1}\left(-1\right)^i algvol\left(str\left(\tilde{x},\ldots\widehat{g_1\ldots g_i\tilde{x}},\ldots,g_1\ldots g_{d+1}\tilde{x}\right)\right)$$
% $$=
\int_{\partial str\left(\tilde{x},g_1\tilde{x},\ldots,g_1\ldots g_{d+1}\tilde{x}\right)}dvol=
\int_{str\left(\tilde{x},g_1\tilde{x},\ldots,g_{d+1}\tilde{x}\right)}d\left(dvol\right)=0.$$ 
Thus $c\nu_d$ is a simplicial cocycle on $BG$. 
%(Its cohomology class does not depend on $\tilde{x}$, because any $\tilde{x}\in G/K$ can be mapped to any other $\tilde{x}_0\in G/K$ by some $g\in G$ and the action of $G$ on $G/K$ is homotopic to the identity and therefore preserves cohomology classes.) 

Consider the cocycle $\nu_d\in C_c^d\left(G;{\Bbb R}\right)$ given by the (clearly continuous) mapping
$$\nu_d\left(g_0,\ldots,g_d\right)=c\nu_d\left(g_0^{-1}g_1,\ldots,g_{d-1}^{-1}g_d\right)=\int_{str\left(\tilde{x},g_0^{-1}g_1\tilde{x},
\ldots,g_{d-1}^{-1} g_d\tilde{x}\right)}dvol=\int_{str\left(g_0\tilde{x},g_1\tilde{x},\ldots,g_d\tilde{x}\right)}dvol.$$
It defines a cohomology 
		class $v_d:=\left[\nu_d\right]\in H^d_c\left(G;{\Bbb R}\right)$ such that $comp\left(v_d\right)\in H^d\left(BG;{\Bbb R}\right)$ 
is represented by $c\nu_d$. The {\em volume class} $v_d$ does not depend on the chosen $\tilde{x}\in G/K$.
%, it will be called the {\em volume class}.

		\begin{thm}\label{Thm1} Let $M=\Gamma\backslash
		G/K$ be a closed, oriented, connected, d-dimensional locally
		symmetric space of noncompact type,
let $j:\Gamma\rightarrow G$ be
the inclusion
of $\Gamma=\pi_1M$ and $Bj_*:H_*^{simp}\left(B\Gamma;{\Bbb Z}\right)\rightarrow H_*^{simp}\left(BG;{\Bbb Z}\right)$ the induced homomorphism.
%$j_*:H_d\left(M;{\Bbb Z}\right)\cong H_d\left(\Gamma;{\Bbb Z}\right)\rightarrow
%	H_d\left(G;{\Bbb Z}\right)$ be the induced homomorphism, and
Let $\left[M\right]\in H_d\left(M;{\Bbb Z}\right)$ 
be	the fundamental class of $M$. Then $$vol\left(M\right)=<v_d,Bj_*EM^{-1}
	\left[M\right]>.$$\end{thm}
	\begin{pf} 
Let $\sum_{i=1}^r a_i\sigma_i$ represent $\left[M\right]$. 
%(One may choose e.g.\ a triangulation $\sigma_1+\ldots+\sigma_r$.) Then $vol\left(M\right)=\sum_{i=1}^r a_i 
%algvol\left(\sigma_i\right)$. 
Fix $\tilde{x}_0\in\widetilde{M}$ and $x_0:=\pi\left(x_0\right)\in M$.
%By the discussion in the last part of Section 2.1,
Then also $\sum_{i=1}^r a_i\tau_i:=\sum_{i=1}^r a_i str\left(\sigma_i\right)\in C_*^{str,x_0}\left(M\right)$ represents $\left[M\right]$, and $vol\left(M\right)=\sum_{i=1}^r a_i algvol\left(\tau_i\right)$ from Stokes' Theorem. 
%(Possibly after straightening some simplices overlap, so we may
%not get a triangulation. However, it will be sufficient to have a fundamental 
%cycle consisting of straight simplices.)
%From Section 2.1, the isomorphism
%$\Phi:C_*^{str,x_0}\left(M\right)\rightarrow C_*^{simp}\left(B\Gamma\right)$ maps 
%$\tau_i$ to $\left(\gamma_1^i,\ldots,\gamma_d^i\right)\in C_d^{simp}\left(B\Gamma\right)$,
%where 
Let $\gamma_j^i\in\Gamma$ be the homotopy class (rel.\ vertices) of the closed edge from
$\tau_i\left(w_{j-1}\right)$ to $\tau_i\left(w_j\right)$, then
$\tau_j=\pi\left(str\left(\tilde{x}_0,\gamma_1^i\tilde{x}_0,\ldots,\gamma_1^i\ldots\gamma_d^i\tilde{x}_0\right)\right).$
Thus from $EM^{-1}=\Phi_*str_*$ we have that
$Bj_*EM^{-1}\left[M\right]\in H_d\left(G;{\Bbb Z}\right)$ is represented by $$\sum_{i=1}^r a_i \left(1,\gamma_1^i,\ldots,\gamma_d^i\right)\in C_d^{simp}\left(BG\right).$$
But $$c\nu_d\left(\gamma_1^i,\ldots,\gamma_d^i\right)=
%\int_{\left(1,\gamma^i_1,\ldots,\gamma^i_d\right)}p^*dvol=\int_{p\left(\gamma^i_1,\ldots,\gamma^i_d\right)}dvol
\int_{str\left(\tilde{x}_0,\gamma_1^i\tilde{x}_0,\ldots,\gamma_1^i\ldots\gamma_d^i\tilde{x}_0\right)}dvol=\int_{\tau_i}dvol=algvol\left(\tau_i\right),$$
%is the volume of the straight simplex with vertices 
%$\tilde{x}_0,\gamma_1^i\tilde{x}_0,\ldots,\gamma_n^i\tilde{x}_0$, 
%i.e.\ of the lift of $\tau_i$ to $\widetilde{M}$ with first vertex 
%$\tilde{x}_0$. Hence $<v_n,\left(1,\gamma_1^i,\ldots,\gamma_n^i\right)>
%=vol\left(\tau_i\right)$, 
which implies
$$<v_d,Bj_*EM^{-1}\left[M\right]>
%=comp\left(v_d\right)\left(Bj_*EM^{-1}\left[M\right]\right)
=c\nu_d\left(\sum_{i=1}^r a_i\left(\gamma_1^i,\ldots,\gamma_d^i\right)\right)
=\sum_{i=1}^r a_i algvol\left(\tau_i\right)=vol\left(M\right).$$
		\end{pf}

		\subsection{Borel classes}
\subsubsection{Dual symmetric space and dual representations}
		Let $\widetilde{M}=G/K$ be a symmetric space of noncompact type.
		Then $G$ is a semisimple, connected Lie group and $K$ is
		a maximal compact subgroup, see \cite[Chapter VI.1]{hel}.

		Let $\g$ be the Lie algebra of $G$ and $\g=\kk\oplus\p$ its Cartan decomposition.
%, and
%$\kk\subset\g$ be the Lie algebra of $K$. 
%		There is the Cartan decomposition $\g=\kk\oplus\p$ with $\left[\kk,\kk\right]\subset\kk,\left[\kk,\p\right]\subset \p,\left[\p,\p\right]\subset\kk$. 
The Killing
		form $B\left(X,Y\right)=Tr\left(ad\left(X\right)\circ ad\left(Y\right)
\right)$ is negatively definite on $\kk$, positively definite on $\p$.

		The dual symmetric space to $G/K$
is $G_u/K$, where $G_u$ is the simply connected Lie group with                                                                       Lie algebra $\g_u=\kk\oplus i\p\subset \g\otimes{\Bbb C}$,
cf.\ \cite[Chapter V.2.]{hel}. The Killing form on $\g_u$ is negatively definite, thus $G_u/K$ is a compact symmetric space.   

		The Lie algebra cohomology                                                     $H^*\left(\g\right)$ is the cohomology of the complex                          $\left(\Lambda^*\g,d\right)$ with $d\phi\left(X_0,\ldots,X_n\right)=
\sum_{i<j}\left(-1\right)^{i+j}\phi\left(\left[X_i,X_j\right],X_0,\ldots,\hat{X_i},\ldots,\hat{X_j},\ldots,X_n\right)$. 

The relative Lie algebra cohomology $H^*\left(\g,\kk\right)$ is the cohomology of the             
                 subcomplex $\left(C^*\left(\g,\kk\right),d\right)\subset \left(\Lambda^*\g,d\right)$ with 
$C^*\left(\g,\kk\right)=\left\{\phi\in \Lambda^*\g: 
i\left(X\right)\phi=0, ad\left(X\right)\phi=\phi\ \ \forall X\in\kk\right\}$, where $i\left(X\right)$ means insertion as first variable, cf.\ \cite[Section 5.5]{bu}. 

If $G/K$ is a symmetric 
space of noncompact type, and $G_u/K$ its compact 
dual, then there is an obvious isomorphism $H^*\left(\g,\kk\right)\rightarrow
H^*\left(\g_u,\kk\right)$, 
%given by multiplication with $i^{dim\left(\p\right)}$.  
dual to the obvious
${\Bbb R}$-linear map $\kk\oplus i\p\rightarrow \kk\oplus \p$.

Moreover, $H^*\left(\g,\kk\right)$ is the cohomology of the complex            of $G$-invariant differential forms 
on $G/K$, cf.\ \cite[Example 5.39]{bu}. Since $G_u$ is
		compact and connected, there is an isomorphism 
		$ H^*\left(G_u/K;{\Bbb R}\right)\rightarrow  H^*\left(\g_u,\kk\right)$, defined by averaging over $G_u$. 
%(Each closed form representing a deRham cohomology class on $G_u/K$ can be averaged over the compact group 
%$G_u$ to obtain a cohomologous $G_u$-invariant form.)

%		For example, $$H^*\left(spin\left(d,1\right),spin\left(d\right)\right)\cong H^*\left(Spin\left(d+1\right)/Spin\left(d\right);{\Bbb R}\right)=H^*\left({\Bbb S}^d;{\Bbb R}\right).$$
%{\bf Dualizing representations.} 
%Let $\rho:G\rightarrow GL\left(N,{\Bbb C}\right)$ be a representation.
%Upon conjugation we can assume that
                %$\rho$
                %sends $K$ to $U\left(N\right)$.

                \begin{df}\label{dualhom}  Let $\widetilde{M}=G/K$ be a symmetric space of 
noncompact type.
                Let $\rho:\left(G,K\right)\rightarrow \left(GL\left(N,{\Bbb C}\right),
                U\left(N\right)\right)$ be a smooth representation.
                We denote $$D_e\rho:\left(\g,\kk\right)\rightarrow \left(gl\left(
N,{\Bbb C}\right),u\left(N\right)\right)$$ the associated Lie-algebra homomorphism,
and, with $\g=\kk\oplus\p, \g_u:=\kk\oplus i\p$,
                $$D_e\rho_u:\left(\g_u,\kk\right)\rightarrow \left(u\left(N\right
)\oplus u\left(N\right),u\left(N\right)\right)$$
                the induced homomorphism on $\kk\oplus i\p$. The corresponding Lie
                group homomorphism $\rho_u:\left(G_u,K\right)\rightarrow \left(U
\left(N\right)\times U\left(N\right),U\left(N\right)\right)$
                will be called the dual homomorphism to $\rho$. Denote $\overline{\rho}_u:G_u/K\rightarrow U\left(N\right)$ the
induced smooth map.\end{df}
\noindent
                Here $\g_u$,
                $\kk$ and $i\p$ are to be understood as subsets of the 
complexification $\g\otimes{\Bbb C}$.
%, and $G_u$ is the simply connected Lie group with Lie
%algebra $\g_u$.
                In particular, the complexification of $gl_N{\Bbb C}$ is
isomorphic to $gl_N{\Bbb C}\oplus gl_N{\Bbb C}$, and $i\p\simeq u\left(
N\right)$ in this case. We emphasize that {\em $\rho_u$ sends $K$ to the
first factor of $U\left(N\right)\times U\left(N\right)$},
and not to the diagonal subgroup as has been claimed in \cite[page 586]{gon}.

		\subsubsection{ Van Est Theorem}
The van Est Theorem \cite[Theorem 6.9]{bu} states
%, for a connected Lie group $G$ and a maximal compact subgroup $K$,
that there is a natural isomorphism
$$\nu_G:H_c^*\left(G;{\Bbb R}\right)\rightarrow H^*\left(\g,\kk\right).$$

If $\rho:\left(G,K\right)\rightarrow \left(GL\left(N,{\Bbb C}\right),U\left(N\right)\right)$ is a representation
%, sending $K$ to $U\left(N\right)$, 
then we obtain the following commutative diagram, where all vertical arrows are isomorphisms

$$\begin{xy}
\xymatrix{        H_c^*\left(GL\left(N,{\Bbb C}\right);{\Bbb R}\right)
        \ar[r]^{\rho^*}&H_c^*\left(G;{\Bbb R}\right)\\
        H^*\left(gl\left(N,{\Bbb C}\right),u\left(N\right)\right)
\ar[r]^{D_e\rho^*} \ar[u]^{\nu_{GL\left(N,{\Bbb C}\right)}^{-1}}_{\cong}
        &H^*\left(\g,\kk\right) \ar[u]^{\nu_G^{-1}}_{\cong}\\
        H^*\left(u\left(N\right)\oplus u\left(N\right),
u\left(N\right)\right)
\ar[r]^{\mbox{\ \ \ \ \ \ }D_e\rho_u^*} \ar[u]_{\cong}
        &H^*\left(\g_u,\kk\right)\ar[u]_{\cong}\\
        H^*\left(U\left(N\right);{\Bbb R}\right) \ar[r]^{\overline{\rho}_u^*} \ar[u]_{\cong}&
        H^*\left(G_u/K;{\Bbb R}\right)\ar[u]_{\cong}}
        \end{xy}$$

%For a connected, semisimple Lie group $G$ 
%with maximal compact subgroup $K$ we 
%We will dennote by $D_G:H^*\left(G_u/K;{\Bbb R}\right)\rightarrow H_c^*\left(G;{\Bbb R}\right)$ the isomorphism given by the right-hand column of the above diagram.
%If $dim\left(G/K\right)=d$, then $G_u/K$ is a $d$-dimensional, compact, 
%orientable manifold and we have $H^d_c\left(G;{\Bbb R}\right)\cong
%H^d\left(G_u/K;{\Bbb R}\right)\cong
%{\Bbb R}$. The volume 
%class $\left[v_d\right]$
\begin{cor}\label{volumeform} Let $G$ be a connected, semisimple Lie group, $K$ a maximal compact subgroup, $d=dim\left(G/K\right)$, $v_d\in H_c^d\left(G; {\Bbb R}\right)$ the volume class, 
%corresponds (under the van Est isomorphism) to the 
$\left[dvol\right]\in H^d\left(G_u/K;{\Bbb R}\right)$ the de Rham 
cohomology class of the volume form on $G_u/K$ and $$D_G:H^*\left(G_u/K;{\Bbb R}\right)\rightarrow H^*_c\left(G;{\Bbb R}\right)$$ the isomorphism given by 
the right-hand column of
the above diagram. Then $$D_G\left(\left[dvol\right]\right)=v_d.$$ \end{cor}
	\begin{pf} By \cite[Proposition 1.5]{dup} $\nu_G$ maps $v_d$ to the 
	class of the volume form in $H_d\left(G/K;{\Bbb R}\right)\cong H^d\left(\g,\kk\right)$. The Riemannian metrics on $G/K$ and $G_u/K$ are defined by the negative of the
Killing form. The ${\bf R}$-linear map $\kk\oplus i\p\rightarrow \kk\oplus\p$ clearly
preserves the Killing form, thus the isomorphism $H^d\left(\g,\kk\right)\simeq H^d
\left(\g_u,\kk\right)$ maps the volume form of $G/K$ to that of $G_u/K$.\end{pf}\\

		\noindent
\subsubsection{Chern classes and Borel classes}

Let $H$ be a {\em compact} connected
Lie group. Let $I_S^n\left(H\right)$ resp.\ $I_A^n\left(H\right)$
be the $ad$-invariant symmetric resp.\ antisymmetric multilinear $n$-forms on its Lie algebra $\h$.
By \cite[Proposition 5.2]{bu} we have the isomorphism $\Phi_A:I_A^n\left(H\right)\rightarrow H^{n}\left(H;{\Bbb R}\right)$. Moreover,
we remind (\cite[Theorem 5.23]{bu}) that there is the Chern-Weil isomorphism 
$\Phi_S:I_S^n\left(H\right)\rightarrow H^{2n}\left(BH;{\Bbb R}\right)$, where {\em in this section (contrary to the remainder of the paper) $BH$
means the classifying space for $H$ with its Lie group topology}. 

For $H=U\left(N\right)$ we 
%It is well-known (\cite[Theorem 4.12]{bur}) that $H^*\left(BU\left(N\right);{\Bbb Z}\right)={\Bbb Z}\left[C_1,\ldots,C_n\right]$, where $C_i\in H^{2i}\left(BU\left(N\right);{\Bbb Z}\right)$
%are the universal integer Chern classes for $i=1,\ldots,N$. 
%The universal twisted Chern classes are defined (see \cite[Definition 4.10]{bur}) as follows: for $p\in\left\{1,\ldots,n\right\}$ let
%$$Z\left(p\right):=\left(2\pi i\right)^p{\Bbb Z}\subset{\Bbb C}$$
%and define
%$$c_p=\left(2\pi i\right)^p C_p\in H^{2p}\left(BU\left(N\right);{\Bbb Z}\left(p\right)\right).$$
consider the 
%$ad$-invariant 
symmetric polynomial $Tr_n\in I_S^n\left(U\left(N\right)\right)$ defined
by $$
Tr_n\left(A_1,\ldots,A_n\right)=\frac{1}{\left(2\pi i\right)^{n}}\frac{1}{n!}\sum_{\sigma\in S_n}Tr\left(A_{\sigma\left(1\right)}\ldots A_{\sigma\left(n\right)}\right)\in I_S^n\left(U\left(N\right)\right).$$ 
%then we have by \cite[Proposition 5.27]{bu} that 
$$ch_{n}:=\Phi_S\left(Tr_n\right)\in H^{2n}\left(BU\left(N\right);{\Bbb Q}\right),$$
is the $2n$-th component of the {\em universal Chern character}. (We consider the rational valued Chern character whose $2n$-th component is 
obtained by multiplication with $\frac{1}{\left(2\pi i\right)^n}$ from that of
the twisted Chern character considered in \cite[Proposition 5.27]{bu}.)

%There is a fibration $H\rightarrow EG_u\rightarrow
%BG_u$ and an associated 
There is a 'transgression map' $\tau$ which maps a subspace of 
$H^{2n-1}\left(H;{\Bbb Z}\right)
$ (whose elements are the so-called transgressive elements) to 
%the quotient 
$H^{2n}\left(BH;{\Bbb Z}\right)$, cf.\ \cite[Example 4.16]{bu}.
%/ker\left(s\right)$, where $s$ is the so-called suspension homomorphism.
%, cf.\ \cite[page 410]{bo2}. 
%If $G=U\left(N\right)$, then, 
%according to Borel
%\cite[p.412]{bo2} one has $$H^*\left(U\left(N\right);{\Bbb Z}\right)\cong \Lambda_{\Bbb Z}\left(b_1,b_3,b_5,\ldots,b_{2N-1}\right),
%$$
%for transgressive elements $b_{2n-1}\in
%H^{2n-1}\left(U\left(N\right);{\Bbb Z}\right)$ which
%satisfy $$\tau\left(b_{2n-1}\right)=C_{n}.$$
%The classes $b_{2n-1}$ are called Borel classes.
%By definition, Borel classes exist only in odd degrees.
%(We will use the same notation for $b_{2n-1}\in H^{2n-1}\left(U\left(N\right);{\Bbb Z}\right)$ and for its image in
%$H^{2n-1}\left(U\left(N\right);{\Bbb R}\right)$.)\\
By \cite{ca} there is a homomorphism 
$$R:I_S^{n}\left(H\right)\rightarrow I_A^{2n-1}\left(H\right),$$
such that the image of $\Phi_A\circ R$ in
%corresponds (after the isomorphism $I_A^{2n-1}\left(G\right)\cong
$H^{2n-1}\left(H;{\Bbb R}\right)$) consists precisely of the transgressive elements and such that $\tau\circ\Phi_A\circ R= \Phi_S$.
%, where $\pi: 
%H^{2n}\left(BH;{\Bbb Z}\right)\rightarrow
%H^{2n}\left(BH;{\Bbb Z}\right)/ker\left(s\right)$ is the projection. 

For $H=U\left(N\right)$,
\cite[Example 5.37]{bu} gives an explicit representative for the {\em Borel classes}
$$b_{2n-1}:=\Phi_A\left(R\left(Tr_n\right)\right)\in H^{2n-1}\left(U\left(N\right);{\Bbb R}\right)\cong H^{2n-1}\left(u\left(N\right)\right)$$ by the Lie algebra cocycle whose value on $X_1,\ldots,X_{2n-1}\in u\left(N\right)$ is
%$$b_{2n-1}\left(X_1,\ldots,X_n\right):=R\left(c_n\right)\left(X_1,\ldots,X_{2n-1}\right)=$$
$$\frac{1}{\left(2\pi i\right)^n}\frac{\left(-1\right)^{n-1}\left(n-1\right)!}{\left(2n-1\right)!}\sum_{\sigma\in S_{2n-1}}\left(-1\right)^\sigma Tr\left(X_{\sigma\left(1\right)}\left[X_{\sigma\left(2\right)},
X_{\sigma\left(3\right)}\right]\ldots\left[X_{\sigma\left(2n-2\right)},
X_{\sigma\left(2n-1\right)}\right]\right).$$
%with $X_1,\ldots,X_{2n-1}\in u\left(N\right)$.

%According to Borel
%(\cite[p.412]{bo2}) one has $$H^*\left(U\left(N\right);{\Bbb Z}\right)\cong 
%\Lambda_{\Bbb Z}\left(b_1,b_3,b_5,\ldots,b_{2N-1}\right).
%$$
%Again we have multiplied the formula in \cite{bu} by $\frac{1}{\left(2\pi i\right)^n}$ to work with a 
%real-valued class. Thus $b_{2n-1}$ equals $\frac{1}{\left(2\pi i\right)^n}\Phi_{2n-1}$ in the notation of \cite[Section 9.7]{bu}. 
%The {\em Borel element} $Bo_n\in C^*\left(gl\left(N,{\Bbb C}\right), u\left(N\right); {\Bbb R}\left(n-1\right)\right)$ is
%defined in \cite[Section 9.7]{bu} by $Bo_n\left(\wedge_{j=1}^{2n-1}x_j\right)=\Phi_{2n-1}\left(\wedge_{j=1}^{2n-1}
%\left(\overline{x}_j^t+x_j\right)\right))$ and pairing with the Borel element defines the Borel regulator. For this paper it would actually be sufficient
%to work with $b_{2n-1}$, but we will give explicit computations of the Borel regulator in Section 3.4.
From $\tau\circ\Phi_A\circ R=\Phi_S$ we conclude that $\tau\left(b_{2n-1}\right)=ch_n$. 
\begin{lem}\label{borelnotzero} Let $G/K$ be a symmetric space of noncompact type, of odd dimension $d=2n-1$, $\rho:G\rightarrow GL\left(N,{\Bbb C}\right)$ a
representation.
             %   A representation $\rho:G\rightarrow GL\left(N,{\Bbb C}\right)$ has
             Then
                 $$\rho_*b_{2n-1}
\not=0
                \in H^{2n-1}_c \left(G;{\Bbb R}\right)\Leftrightarrow
\overline{\rho}_u^*
                b_{2n-1}\not=0\in H^{2n-1}\left(G_u/K;{\Bbb R}\right)\Leftrightarrow
<b_{2n-1}, \left(\rho_u\right)_*\left[G_u/K\right]>\not=0.$$\end{lem}

                \begin{pf} This follows
from naturality of $D_G$
and from $H_d\left(G_u/K;{\Bbb R}\right)\cong{\Bbb R}$.\end{pf}\\

	It will be clear from the context whether
		we consider the Borel classes 
as elements of $H^*\left(u\left(N\right)\right)
\simeq
H^*\left(U\left(N\right);{\Bbb R}\right)$ or as the (under the van Est
isomorphism) corresponding elements of $H_c^*\left(GL\left(N,{\Bbb C}
                \right);{\Bbb R}\right)$. 

Stabilization $H^*\left(U\left(N+1\right);{\Bbb R}\right)\rightarrow H^*\left(U\left(N\right);
{\Bbb R}\right)$ preserves $b_{2n-1}$, 
thus $b_{2n-1}$ may also be considered as an element of $H^{2n-1}\left(U;
{\Bbb R}\right)\cong 
H^{2n-1}_c\left(GL\left({\Bbb C}\right);{\Bbb R}\right)$.

%Let $n\in{\Bbb N}$.
%We consider $K_{2n-1}\left(
%{\Bbb C}\right)\otimes{\Bbb Q}$       
%as a subgroup of      
%$H_{2n-1}\left(BGL\left({\Bbb C}\right);{\Bbb Q}\right)$, as in Section 2.2. 
%		The {\bf Borel regulator} 
In \hyperref[Thm2]{Theorem \ref*{Thm2}} and \hyperref[Thm3]{Theorem \ref*{Thm3}} we will for subrings $A\subset {\Bbb C}$ consider the homomorphism
$$K_{2n-1}\left(A\right)\otimes{\Bbb Q}\rightarrow{\Bbb R}$$ 
defined by application of the isomorphism
$$K_{2n-1}\left(
A\right)\otimes{\Bbb Q}\cong PH_{2n-1}\left(\mid BGL\left(A\right)\mid^+;{\Bbb Q}\right)\cong
PH_{2n-1}\left(BGL\left(A\right);{\Bbb Q}\right)$$ from Section 2.2 and
pairing
with		the Borel class 
		      $$b_{2n-1}\in H_c^{2n-1}\left(GL\left({\Bbb C}\right);{\Bbb R}\right).$$

\subsection{Projection $H_*^{simp}\left(BGL\left(\overline{\Bbb Q}\right);{\Bbb Q}\right)\rightarrow K_*\left(\overline{\Bbb Q}\right)\otimes{\Bbb Q}$}

Let $A\subset{\Bbb C}$ be a subring and $G=GL\left(A\right)$. Let $I=H_{simp}^{*\ge 1}\left(BG;{\Bbb Q}\right)$ 
be the augmentation ideal of $H^*_{simp}\left(BG;{\Bbb Q}\right)$ and $D=I^2$ the
subspace of decomposable cohomology classes. 

Let $PH_*^{simp}\left(BG;{\Bbb Q}\right)$ be the subspace of primitive elements in homology.
It is easy to check that $c\left(h\right)=0$ for all $c\in D, h\in PH_*^{simp}\left(BG;{\Bbb Q}\right)$.
By \cite[Proposition 3.10]{mm} $I/D$ is 
the dual of $PH_*^{simp}\left(BG;{\Bbb Q}\right)$, which implies
$$D=\left\{c\in I: c\left(h\right)=0\ \forall\ h\in PH_*^{simp}\left(BG;{\Bbb Q}\right)\right\}
$$
%\mbox{\ and\ }
$$PH_*^{simp}\left(BG;{\Bbb Q}\right)=\left\{h\in H_*^{simp}\left(BG;{\Bbb Q}\right): c\left(h\right)=0\ \forall\ c\in D\right\}.$$

%In Section 2.4 
%We have the Borel classes $b_{2n-1}\in H_c^{2n-1}\left(GL\left({\Bbb C}\right);{\Bbb R}\right)$ and, for
%each subring $A\subset {\Bbb C}$, the classes $comp\left(b_{2n-1}\right)\in 
%H^{2n-1}\left(BGL\left(A\right);{\Bbb R}\right)$.

\begin{lem}\label{proj} Let $A\subset{\Bbb C}$ be a subring. Assume that $comp\left(b_{2n-1}\right)\in H^{2n-1}_{simp}\left(BGL\left(A\right);{\Bbb R}\right)$
is not decomposable: $comp\left(b_{2n-1}\right)\not\in D$.
Then there exists a projection $$pr_{2n-1}:H_{2n-1}^{simp}\left(BGL\left(A\right);{\Bbb Q}\right)\rightarrow PH_{2n-1}^{simp}\left(BGL\left(A\right);{\Bbb Q}\right)$$
such that for all $h\in H_{2n-1}^{simp}\left(BGL\left(A\right);{\Bbb Q}\right)$
$$comp\left(b_{2n-1}\right)\left(pr_{2n-1}\left(h\right)\right)=comp\left(b_{2n-1}\right)\left(h\right).$$
\end{lem}

\begin{pf}
Denote $G=GL\left(A\right)$. 
We consider $comp\left(b_{2n-1}\right)\in H^{2n-1}_{simp}\left(BG;{\Bbb Q}\right)$ as a linear map $comp\left(b_{2n-1}\right): 
H_{2n-1}^{simp}\left(BG;{\Bbb Q}\right)\rightarrow {\Bbb Q}.$
We have $$PH_{2n-1}^{simp}\left(BG;{\Bbb Q}\right)=\left\{h\in H_{2n-1}^{simp}\left(BG;{\Bbb Q}\right): c\left(h\right)=0\ \ \forall c\in D\right\}.$$

Since $comp\left(b_{2n-1}\right)\not\in D$ there exists some $e_0\in PH_{2n-1}^{simp}\left(BG;{\Bbb Q}\right)$ with $$comp\left(b_{2n-1}\right)
\left(e_0\right)\not=0.$$
We extend $\left\{e_0\right\}$ to a basis $\left\{e_j:j\in J_P\right\}$ of $PH_{2n-1}^{simp}\left(BG;{\Bbb Q}\right)$ 
and then to a basis $\left\{e_j:j\in J\right\}$ of $H_{2n-1}^{simp}\left(BG;{\Bbb Q}\right)$, for some index sets $\left\{0\right\}\subset J_P\subset J$. 

Since $comp\left(b_{2n+1}
\right)\left(e_0\right)\not=0$, we have that $\left\{e_j^\prime:j\in J\right\}$ defined by
$$e_0^\prime:=e_0, e_j^\prime:=comp\left(b_{2n-1}\right)\left(e_0\right)e_j-comp\left(b_{2n-1}\right)\left(e_j\right)e_0 \mbox{\ for\ }j\in J-\left\{0\right\}$$
is another basis 
of $H_{2n-1}^{simp}\left(BG;{\Bbb Q}\right)$, and $\left\{e_j^\prime:j\in J_P\right\}$ is another basis of
$PH_{2n-1}^{simp}\left(BG;{\Bbb Q}\right)$. 

%By construction, $\left\{e_j^\prime: j\in J-\left\{0\right\}\right\}$ is a basis of $ker\left(comp\left(b_{2n-1}\right)\right)$. Thus, if we 
Let $S\subset H_{2n-1}^{simp}\left(BG;{\Bbb Q}\right)$ be the subspace spanned by $\left\{e_j^\prime: j\not\in J_P\right\}$, then $S\subset ker\left(comp\left(b_{2n-1}\right)\right)$ and we have a decomposition
$$H_{2n-1}^{simp}\left(BG;{\Bbb Q}\right)=PH_{2n-1}^{simp}\left(BG;{\Bbb Q}\right)\oplus S.$$

We use this decomposition to define the projection 
$pr_{2n-1}:H_{2n-1}^{simp}\left(BG;{\Bbb Q}\right)\rightarrow PH_{2n-1}^{simp}\left(BG;{\Bbb Q}\right)$ by $pr_{2n-1}\left(p+s\right)=p$ for $p\in PH_{2n-1}^{simp}\left(BG;{\Bbb Q}\right)$ and $s\in S$. $S\subset ker\left(comp\left(b_{2n-1}\right)\right)$ implies $comp\left(b_{2n-1}\right)\left(pr_{2n-1}\left(p+s\right)\right)=comp\left(b_{2n-1}\right)\left(p+s\right)$.\end{pf}\\

To decide whether the Borel class is indecomposable we apply\footnotemark\footnotetext[2]{We remark that in the already interesting case $A={\Bbb C}$ one can prove indecomposability of the Borel class without using Borel's K-theory computation. 

First, $H_c^*\left(GL\left(N,{\Bbb C}\right);{\Bbb Q}\right)=
\Lambda_{\Bbb Q}\left(b_1,b_3,b_5,\ldots,b_{2N-1}\right)$ implies that $b_{2n-1}$ is not decomposable in
$H_c^*\left(GL\left(N,{\Bbb C}\right);{\Bbb Q}\right)$ for any $N$.
Next, by homology stability of the linear group (\cite[page 77]{bu}), inclusion induces an isomorphism $H^{2n-1}\left(BG;{\Bbb Q}\right)=
H^{2n-1}\left(BGL\left(N,{\Bbb C}\right);{\Bbb Q}\right)$ if 
%$2n-1\le\frac{N-1}{2}$, that is if 
$N\ge 4n+3$.

By Borel's Theorem (see \cite[Theorem 9.6]{bo3}), for each arithmetic subgroup $\Gamma\subset SL\left(N,{\Bbb C}\right)$ we have an isomorphism $j^*:H^{2n-1}_{simp}
\left(B\Gamma;{\Bbb Q}\right)\rightarrow H^{2n-1}_c\left(BSL\left(N,{\Bbb C}\right);{\Bbb Q}\right)$ whenever 
%$2n-1\le\frac{N}{4}$, that is if 
$N\ge 8n+4$. This isomorphism is constructed via the van Est isomorphism, that 
is by integration of forms over simplices. In particular, if $h\in H^*_{simp}\left(
BSL\left(N,{\Bbb C}\right);{\Bbb Q}\right)$ and $i:\Gamma\rightarrow SL\left(N,{\Bbb C}\right)$ is the inclusion, then $comp\left(j^*i^*h\right)=h$.

Now 
%we prove $comp\left(b_{2n-1}\right)\not\in D$ by contradiction. Assume $comp\left(b_{2n-1}\right)$ were decomposable, that is 
if $comp\left(b_{2n-1}\right)=xy$, 
%where $x,y\in I$ are cohomology classes of degree $\ge 1$. Fix some $N\ge 8n+4>4n+3$. Then 
then $b_{2n-1}=j^*i^*comp\left(b_{2n-1}\right)=\left(j^*i^*x\right)\left(j^*i^*y\right)$ is decomposable in $H_c^*\left(GL\left(N,{\Bbb C}\right);{\Bbb Q}\right)$ for for all $h\in H_{2n-1}^{simp}\left(BGL\left({\overline{\Bbb Q}}\right);{\Bbb Q}\right)$, giving a contradiction.}
    Borel's computation of K-theory of integer rings in number fields in \cite{bo3}.

Let $O_F$ be the ring of integers in a number field $F$, which has $r_1$ real and $2r_2$ complex embeddings. Borel proves that the Borel regulator, applied to the different embeddings of $SL\left(O_F\right)$, yields an isomorphism between $PH_{2n-1}^{simp}\left(BSL\left(O_F\right);{\Bbb Z}\right)$ and ${\Bbb Z}^{r_1+r_2}$ resp.\ ${\Bbb Z}^{r_2}$
if $n$ is even resp.\ odd. 
Since decomposable cohomology classes vanish on primitive homology classes, this implies in particular:\\

{\em If $A=O_F$ for a number field F, then $b_{2n-1}$ is not decomposable for even $n$.

If moreover $F$ is not totally real, then $b_{2n-1}$ is not decomposable for all $n$.}\\
\\
In particular, we can apply \hyperref[proj]{Lemma \ref*{proj}} to $A=O_F$, and therefore also to
%If $A_1\subset A_2\subset {\Bbb C}$ are subrings and the Borel class is not decomposable for $A_1$, then of course it is also not decomposable for $A_2$.
%Thus we can actually apply \hyperref[proj]{Lemma \ref*{proj}} to 
all rings $A$ with $O_F\subset A\subset{\Bbb C}$,
in particular to $A={\overline{\Bbb Q}}$:

\begin{cor}\label{projq} 
For all $n$, there exists a projection $$pr_{2n-1}:H_{2n-1}^{simp}\left(BGL\left({\overline{\Bbb Q}}
\right);{\Bbb Q}\right)\rightarrow PH_{2n-1}^{simp}\left(BGL\left({\overline{\Bbb Q}}\right);{\Bbb Q}\right)=K_{2n-1}\left({\overline{\Bbb Q}}\right)\otimes{\Bbb Q}$$
such that for all $h\in H_{2n-1}^{simp}\left(BGL\left({\overline{\Bbb Q}}\right);{\Bbb Q}\right)$
$$comp\left(b_{2n-1}\right)\left(pr_{2n-1}\left(h\right)\right)=comp\left(b_{2n-1}\right)\left(h\right).$$
\end{cor}

\subsection{Compact locally symmetric spaces and K-theory}

In this subsection, we finally
show that to each representation of nontrivial Borel class,
and each compact, oriented, locally symmetric space of noncompact type $M$ we can find a
		nontrivial element $\gamma\left(M\right)\in K_*\left(
\overline{\Bbb Q}\right)\otimes{\Bbb Q}$. 
		
\begin{thm}\label{Thm2} For each symmetric space $G/K$ of noncompact type and odd dimension $d$, and
each representation $\rho:G\rightarrow GL\left(N,{\Bbb C}\right)$ with $\rho^*b_{d}\not=0$, there is some $c_\rho\not=0$ such that 
%the following holds:

to each compact, oriented, locally symmetric space $M=\Gamma\backslash G/K$, with $\rho\left(\Gamma\right)\subset 
GL\left(N,A\right)$ for a subring $A\subset {\Bbb C}$ satisfying the conclusion of \hyperref[proj]{Lemma \ref*{proj}},
		there exists an element $$\gamma\left(M\right)\in K_{d}\left( A\right)\otimes{\Bbb Q}$$
		%such that the homomorphism  $K_{d}\left(A\right)\otimes{\Bbb Q}
		%\rightarrow{\Bbb R}$ defined in Section 2.4 by pairing with the Borel class fulfills 
with $<b_{d},\gamma\left(M\right)>=c_\rho vol\left(M\right).$\end{thm}

		\begin{pf}
%		We recall fromhhh \hyperref[Thm1]{Theorem \ref*{Thm1}} we considered  
%		$Bj_*EM^{-1}\left[M\right]\in H_d^{simp}\left(BG;{\Bbb Z}\right)$.
%	Applying $B\rho_*$ we get
%	an element 
%	$$B\left(\rho j\right)_*EM^{-1}\left[M\right]\in H_d^{simp}\left(BGL\left(N,
%{\Bbb C}\right);{\Bbb Z}
%		\right).$$
%Since $B\left(\rho j\right)$ maps $B\Gamma$ to $BGL\left(N,A\right)$, and 
Using the projection $pr_d$ from \hyperref[proj]{Lemma \ref*{proj}}, we obtain as in Section 2.2 
%$$
%\left(H\rho\right)_*\left[M\right]=
%B\left(\rho j\right)_*EM^{-1}\left[M\right]\in H_d^{simp}\left(BGL\left(N,
%A\right);{\Bbb Z}
%                \right).$$
%By projecting (the image of $\left(B\rho\right)_*j_*\left[M\right]$ in 
%$ H_n\left(GL\left(A\right);{\Bbb Q}
%                \right)$) along the kernel of the Borel regulator,
%we obtain an element 
%	$pr_*\left(B\rho\right)_*j_*\left[M\right]\in P_n\left(GL\left(A\right);{\Bbb Q}
%                \right)=K_n\left(A\right)\otimes{\Bbb Q}$.
%By \hyperref[proj]{Lemma \ref*{proj}} we have a projection $pr_{d}$.
%In Section 2.2, we defined $$
$$\gamma\left(M\right):=Q_{*}pr_{*}B\left(\rho j\right)_{*}EM^{-1}\left[M\right]\in 
K_{d}\left(A\right)\otimes{\Bbb Q}.$$

%By assumption $\rho^*b_{d}\not=0$. Since 
$\rho_*b_d\not=0$ together with $H_c^{d}\left(G;{\Bbb R}\right)=H^d\left(G_u/K;{\Bbb R}\right)={\Bbb R}$ implies
$\rho^*b_{d}=c_\rho
v_{d}$ for some $c_\rho\not=0$.
Using \hyperref[proj]{Lemma \ref*{proj}} we get 
$$
%r_{d}\left(\gamma\left(M\right)\right)=
<b_{d},\gamma\left(
M\right)>=
comp\left(b_{d}\right)\left(pr_{*}B\left(\rho j\right)_{*}EM^{-1}\left[M\right]\right)=$$
$$comp\left(b_{d}\right)\left(B\left(\rho j\right)_{*}EM^{-1}\left[M\right]\right)=comp\left(\rho^*b_{d}\right)
\left(Bj_{*}EM^{-1}\left[M\right]\right)$$
$$=c_\rho <v_{d},Bj_{*}EM^{-1}\left[M\right]>= c_\rho vol\left(M\right),$$ where 
the last equality is true by \hyperref[Thm1]{Theorem \ref*{Thm1}}.
\end{pf}\\

\begin{cor}\label{constant} For each symmetric space $G/K$ of noncompact type and odd dimension
$d=2n-1$, and to
each representation $\rho:G\rightarrow GL\left(N,{\Bbb C}\right)$ with $\rho^*b
_{2n-1}\not=0$, there exists a constant $c_\rho\not=0$, such that the following holds:
to each compact, oriented, locally symmetric space $M=\Gamma\backslash G/K$
there exists an element 
$$\gamma\left(M\right)\in K_{2n-1}\left(
\overline{\Bbb Q}\right)\otimes{\Bbb Q}$$
with 
%$r_{2n-1}:K_{2n-1}\left(\overline{\Bbb Q}
%\right)\otimes{\Bbb Q}
%\rightarrow{\Bbb R}$ fulfills 
$<b_{2n-1},\gamma\left(M\right)>
=c_\rho vol\left(M\right).$\end{cor}
\begin{pf} 
%$G$ is a linear semisimple Lie group without compact factors. 
$dim\left(G/K\right)=2n-1$ implies that $G$ is not locally isomorphic to $SL\left(2,{\Bbb R}\right)$, thus we can apply 
%$because $dim\left(SL\left(2,{\Bbb R}\right)/SO\left(2\right)\right)=2$. Thus $G$ satisfies the assumptions of
Weil's rigidity theorem, which yields a $g\in G$ with
$g\Gamma g^{-1}\in G\left(\overline{\Bbb Q}\right)$.
%Thus, upon replacing $\Gamma$ by $g\Gamma g^{-1}$, 
Hence $M$ is of the form $M=\Gamma\backslash G/K$ with $\Gamma\subset G\left(
\overline{\Bbb Q}\right)$.

Each
representation $\rho:G\rightarrow GL\left(N,{\Bbb C}\right)$ is isomorphic to a representation $\rho^\prime$ such that
$G\left(\overline{\Bbb Q}\right)$ is mapped
to $GL\left(N,\overline{\Bbb Q}\right)$. This follows from the classification of irreducible representations of Lie groups, see \cite{fh}.

By \hyperref[projq]{Corollary \ref*{projq}} we can then apply \hyperref[Thm2]{Theorem \ref*{Thm2}} to $A=\overline{\Bbb Q}$.\end{pf}
% satisfies the conclusion of \hyperref[proj]{Lemma \ref*{proj}} and
%we can apply \hyperref[Thm2]{Theorem \ref*{Thm2}}.\end{pf}

\begin{cor}\label{independent} Let $G/K$ be a symmetric space of noncompact type and $\rho:G\rightarrow
GL\left(N,{\Bbb C}\right)$ a representation with $
\rho^*b_{2n-1}\not=0$, for $2n-1=dim\left(G/K\right)$.
                Then compact, oriented,
locally symmetric spaces $\Gamma\backslash G/K$ of rationally independent volumes
yield rationally independent elements in $K_{2n-1}\left(\overline{\Bbb Q}\right)
\otimes{\Bbb Q}$.\end{cor}

\noindent
{\bf Remark}: {\em In \cite{gon} it was claimed that for (2n-1)-dimensional compact hyperbolic manifolds
one can construct an element $\gamma\left(M\right)\in K_{2n-1}
\left(\overline{\Bbb Q}\right)\otimes{\Bbb Q}$ such
that $<b_{2n-1},\gamma\left(M\right)>=vol\left(M\right)$. However,
since $\rho^*b_{2n-1}$ is a rational cohomology class, $c_\rho$ is rational if and only if $v_{2n-1}$ is a rational 
cohomology class, and this is equivalent to $vol\left(M\right)=<v_{2n-1},
\left[M\right]>\in{\Bbb Q}$. Since, conjecturally, all hyperbolic manifolds 
have irrational volumes, one can probably not get rid of the factor $c_\rho$ in \hyperref[Thm2]{Theorem \ref*{Thm2}}.}\\
\\
%If $M$ is defined over some subring $A\subset\overline{\Bbb Q}$, then we do get an element in $K_n\left(A\right)\otimes{\Bbb Q}$. The proof is an obvious modification. 
In conclusion, we are left with the problem of finding representations of nontrivial Borel class, which will be solved in Section 3.\\

{\bf Compact examples} can e.g.\ be obtained
by Borel's construction of locally symmetric spaces in \cite{bo2b}. A very special case is the construction of
arithmetic hyperbolic manifolds using quadratic forms (cf.\ the textbook \cite[Chapter E.3]{bp}).

%Let $n$ be odd.
Let $u\in{\bf R}$ be an algebraic integer such that all roots of its minimal polynomial have multiplicity 1 and are real and negative (except 
possibly $u$). Assume moreover that $\left(0,\ldots,0\right)$ is the only integer solution of $x_1^2+\ldots
+x_{2n-1}^2-ux_{2n}^2=0$. Let $\widehat{\Gamma}\subset GL\left(2n,{\Bbb Z}\left[u\right]\right)$
be the group of maps preserving $x_1^2+\ldots
+x_{2n-1}^2-ux_{2n}^2$. It is isomorphic to a discrete cocompact subgroup of $SO\left(2n-1,1;{\Bbb Z}\left[u\right]\right)\subset SO\left(2n-1,1;{\Bbb R}\right)$. 
By Selberg's lemma, it contains a torsionfree cocompact subgroup $\Gamma\subset SO\left(2n-1,1;{\Bbb Z}\left[u\right]\right)$. With the computations in Section 3 below one concludes:
If $n$ is even, then the compact manifold $M:=\Gamma\backslash{\Bbb H}^{2n-1}$ (and, for example, a half-spinor representation)
gives a nontrivial element
$\gamma\left(M\right)\in K_{2n-1}\left({\Bbb Z}\left[u\right]\right)\otimes {\Bbb Q}$. If $n$ is odd, then \hyperref[projq]{Corollary \ref*{projq}} can not be applied to
${\Bbb Z}\left[u\right]$ but to $\overline{\Bbb Q}$,
one gets at least a nontrivial element 
$\gamma\left(M\right)\in K_{2n-1}\left(\overline{\Bbb Q}\right)\otimes {\Bbb Q}$.\\

Matthey-Pitsch-Scherer constructed in \cite{mps}
a somewhat stronger invariant for stably parallelisable manifolds: given an embedding $M^d\rightarrow {\Bbb R}^n$
with trivial normal bundle $\nu M$ and a regular neighborhood $U$ they consider the composition
                $${\Bbb S}^n\rightarrow \overline{U}/\partial U\rightarrow
                \overline{U}/\partial U\wedge M_+=Th\left(\nu M\right)\wedge M_+=
                \Sigma^{n-d}M_+\wedge M_+\rightarrow{\Bbb S}^{n-d}\wedge M_+$$
                giving an element $\gamma\left(M\right)\in \pi_d^s\left(M\right)$.

		\section{Existence of representations of nontrivial Borel class}

\subsection{Trace criterion}
		\begin{lem}\label{cartan} Let $G/K$ be a symmetric space of noncompact type, of dimension $2n-1$. Let $\tq\subset\p$ be a Cartan subalgebra of $\g$.
%, let $G_u/K$ be 
%		its compact dual, and $\tq$ the Cartan subalgebra of the Lie algebra $\g_u$. 

Then for a representation $\rho:G\rightarrow GL\left(N,{\Bbb C}\right)$ and its dual $\rho_u:G_u\rightarrow U\left(N\right)\times U\left(N\right)$
the following are equivalent:\\
i) $\rho$ has nonvanishing
		Borel class $\rho^*b_{2n-1}\not=0\in H_c^{2n-1}\left(G;{\Bbb R}\right),$\\
ii) $Tr\left(\left(D_e\rho_u\left(
                it\right)\right)^n\right)\not=0$ for some $t\in\tq$,\\
iii) $Tr\left(\left(D_e\rho\left(
                t\right)\right)^n\right)\not=0$ for some $t\in\tq$.\end{lem}

		\begin{pf} As in \hyperref[dualhom]{Definition \ref*{dualhom}} we
consider the dual representation $\rho_u:\left(G_u,K\right)\rightarrow \left(U\left(N\right)\times U\left(N\right),U\left(N\right)\right)$
and the smooth map
$\overline{\rho_u}:G_u/K\rightarrow U\left(N\right)\times U\left(N\right)/U\left(
N\right)\simeq U\left(N\right).$ Since $\rho_u$ sends $K$ to the first factor of $U\left(N\right)\times U\left(N\right)$ we have $\pi_2\rho_u=\overline{\rho_u} p$, where $\pi_2:U\left(N\right)\times U\left(N\right)\rightarrow U\left(N\right)$
is the projection to the second factor and $p:G_u\rightarrow G_u/K$ projection to the quotient.

%We have the commutative diagram
%$$\begin{xy}
%\xymatrix{        H_c^*\left(GL\left(N,{\Bbb C}\right);{\Bbb R}\right)
%        \ar[r]^{\rho^*}&H_c^*\left(G;{\Bbb R}\right)\\
%        H^*\left(U\left(N\right);{\Bbb R}\right) \ar[r]^{\overline{\rho_u}^*} \ar[u]_{\cong}&
%        H^*\left(G_u/K;{\Bbb R}\right)\ar[u]_{\cong}}
%        \end{xy}$$
%which implies that 
By \hyperref[borelnotzero]{Lemma \ref*{borelnotzero}} we have that $\rho^*b_{2n-1}\not=0\in H_c^{2n-1}\left(G;{\Bbb R}\right)$ if and only if
$$\overline{\rho_u}^*b_{2n-1}\not=0\in H^{2n-1}\left(G_u/K\right).$$
%The projection $p:G_u\rightarrow G_u/K$ induces an injective map 
Averaging of differential forms over the compact group $K$ shows that $p^*:H^*
\left(G_u/K\right)\rightarrow
H^*\left(G_u\right)$ is injective.
%, because a left inverse to $p^*$
%is given by averaging differential
%forms over the compact group $K$. 
Hence, $\overline{\rho_u}^*b_{2n-1}\not=0$ 
if and only if its image in $H^{2n-1}
\left(G_u\right)$ does not vanish. 
The latter equals $$\left(\pi_2 \rho_u\right)^*
b_{2n-1},$$
because $\pi_2\rho_u=\overline{\rho_u} p$.
%where $\pi_2:U\left(N\right)\times U\left(N\right)
%\rightarrow U\left(N\right)$ is the projection to the second factor. 

%Consider the isomorphisms $$\Phi_A:I_A^{2n-1}\left(G_u\right)\rightarrow
%H^{2n-1}\left(G_u\right),
%\Phi_S: I_S^n\left(G_u\right)\rightarrow H^{2n}\left(BG_u\right)$$
%and the homomorphism $$
%R:I_S^n\left(G_u\right)\rightarrow I_A^{2n-1}\left(
%G_u\right)$$  
%from Section 2.4. According to \cite{ca}, the image of $\Phi_A\circ R$ are the 
%transgressive elements and one has $$\tau\circ \Phi_A\circ R =\Phi_S$$ for the universal 
%transgression map $\tau$. In particular 
Using the notation and facts from Section 2.4.3 we have
%$b_{2n-1}=\Phi_A\left(R\left(Tr_n\right)\right)$, which by 
%thus $\left(\pi_2\rho_u\right)^*b_{2n-1}=\left(\pi_2\rho_u\right)^* R\left(C_n\right)$.
%naturality 
%of the transgression map 
%implies that
$$\left(\pi_2\rho_u\right)^*b_{2n-1}=\left(\pi_2\rho_u\right)^*\Phi_A\left(R\left(Tr_n\right)\right)=
\Phi_A\left(R\left(\left(\pi_2\rho_u\right)^*Tr_{n}\right)\right),$$
%thus $\left(\pi_2\rho_u\right)^*b_{2n-1}\not=0$ 
%implies 
and $$\left(\pi_2\rho_u\right)^*ch_n=\left(\pi_2\rho_u\right)^*\Phi_S\left(Tr_n\right)=\Phi_S\left(\left(\pi_2\rho_u\right)^*Tr_n\right).$$
%=0\in H^{2n}\left(BG_u\right).$$
$\Phi_A$ and $\Phi_S$ are isomorphisms, moreover $\tau\circ\Phi_A\circ R=\Phi_S$ 
implies injectivity of $R$. Hence $\left(\pi_2\rho_u\right)^*ch_n\not=0$ if and only if $\left(\pi_2\rho_u\right)^*b_{2n-1}\not=0$.

%		It is known that $H^*\left(BG_u\right)$ is isomorphic to the algebra $S^{G_u}_*\left(\g_u\right)$ of $G_u$-invariant polynomials on $\g_u$ (see \cite{bo2}).
%Recall that $ch_n=\Phi_S\left(Tr_n\right)$ and that
From the definition of $Tr_n$ we see that $$\left(\pi_2\rho_u\right)^*Tr_n\left(A_1,\ldots,A_n\right)=\frac{1}{\left(2\pi i\right)^n}\frac{1}{n!}\sum_{\sigma\in S_n}Tr\left(D_e\left(\pi_2\rho_u\right) A_{\sigma\left(1\right)}\ldots D_e\left(\pi_2\rho_u\right)A_{\sigma\left(n\right)}\right),$$ 
%undee the 
%isomorphism $H^{2n}\left(BU\left(N\right)\right)\simeq
%I^n_S\left(u\left(N\right)\right)$. 
An easy exercise in multilinear algebra shows that a {\em symmetric} polynomial $P\left(x_1,\ldots,x_n\right)$ is nontrivial if and only if there is some $x$ with $P\left(x,x,\ldots,x\right)\not=0$. Hence it is sufficient to check that the invariant polynomial $$Tr\left(\left(D_e\left(\pi_2\rho_u\right)\left(.\right)\right)^n\right)$$ is not trivial on $\g_u$.

Let $\tq_u$ be the Cartan subalgebra of $\g_u$, which corresponds to $\tq$ under the canonical bijection $\kk\oplus\p\simeq\kk\oplus i\p$.
There is an action of the Weyl group $W$ on $\tq_u$, we denote its space of
invariant polynomials by $S_*^W\left(\tq_u
                \right)$.
By a theorem of Chevalley (see \cite{bou}), restriction induces 
an isomorphism $$S_*^{G_u}\left(\g_u\right)\cong
S_*^W\left(\tq_u
		\right).$$ In particular, it suffices 
		to check that $Tr\left(\left(D_e\left(\pi_2\rho_u\right)\left(.\right)\right)^n\right)$ is not trivial on $\tq_u$.

By assumption the Cartan algebra $\tq$ is contained in $\p$. (This can actually always be achieved by a suitable conjugation.)
%Since the conclusion of \hyperref[cartan]{Lemma \ref*{cartan}} is invariant under conjugation, we can without loss of generality assume that $\tq\subset \p$ and 
Thus $\tq_u\subset i\p$.
This implies that, for $t\in\tq_u$, $D_e\rho_u\left(t\right)$ belongs to the second factor of $u\left(N\right)\oplus u\left(N\right)$, and thus 
$D_e\left(\pi_2\rho_u\right)\left(t\right)=D_e\rho_u\left(t\right)$ for 
$t\in \tq_u$, which proves the equivalence of i) and ii). 
Finally we note that, for $t\in\p$, $Tr\left(\left(D_e\rho\left(t\right)\right)^n\right)$
and $Tr\left(\left(D_e\rho_u\left(it\right)\right)^n\right)$ coincide up to a power of $i$. The equivalence of ii) and iii) follows.
		\end{pf}
\begin{cor}\label{even}
Let $G/K$ be a symmetric space of noncompact type. If $d:=dim\left(G/K\right)\equiv 3\ mod\ 4$, then every nontrivial
representation $\rho:G\rightarrow GL\left(N,{\Bbb C}\right)$ has nonvanishing
                Borel class $\rho^*b_{d}\not=0\in H_c^{d}\left(G;{\Bbb R}\right)$.\end{cor}
\begin{pf} We apply \hyperref[cartan]{Lemma \ref*{cartan}} with $d=2n-1$, that is $n$
is even.

For each $t\in \tq$ we have that $$D_e \rho_u\left(it\right)\in u\left(N\right)\oplus u\left(N\right)$$ has purely imaginary eigenvalues, since matrices in $u\left(N\right)\oplus u\left(N\right)$ are skew-symmetric. Hence, if $\rho$ is nontrivial (and thus $D_e\rho_u\not\equiv 0$), the eigenvalues of 
$\left(D_e \rho_u\left(it\right)\right)^n$ are either all positive (if $n\equiv 0\ mod\ 4$) or all negative (if $n\equiv 2\ mod\ 4$). In either case $Tr\left( \left(D_e \rho_u\left(it\right)\right)^n\right)\not=0$.\end{pf}\\

		\subsection{Borel class of Lie algebra representations}
\subsubsection{Preliminaries}
		Let $\g$ be a semisimple Lie algebra and $
	R\left(\g\right)$ its (real) 
	representation ring, with addition $\oplus$ and multiplication $\otimes$. Let $\tq$ be a Cartan subalgebra of $\g$.
	%(We consider represenations $\pi:\g\rightarrow gl\left(N,{\Bbb C}\right)$ which are not necessarily complex linear.) 

In this section we consider, for $n\in{\Bbb N}$,
the map $\beta_{2n-1}:R\left(\g\right)\rightarrow {\Bbb C}\left[\tq\right]$
given by $$\beta_{2n-1}\left(\pi\right)\left(t\right)=Tr\left(\pi\left(t\right)^n\right).$$
One has $\beta_{2n-1}\left(\pi_1\oplus\pi_2\right)=\beta_{2n-1}\left(\pi_1\right)+\beta_{2n-1}
\left(\pi_2\right)$ 
%holds for representations $\pi_1,\pi_2$. Therefore $\beta_{2n-1}$ is uniquely determined by its values 
%for irreducible representations. Moreover, $
and $\beta_{2n-1}\left(\pi_1\otimes\pi_2\right)=\beta_{2n-1}\left(\pi_1\right)\beta_{2n-1}\left(\pi_2\right)$ for representations $\pi_1,\pi_2$.\\
By \hyperref[cartan]{Lemma \ref*{cartan}} a representation $\rho:G\rightarrow GL\left(N,{\Bbb C}
\right)$ has {\em nontrivial Borel class} $\rho^*b_{2n-1}\not=0\in H_c^{2n-1}\left(G;{\Bbb R}\right)$
if and only if $Tr\left(D_e\rho\left(A\right)^n\right)\not =0$ for some $A\in\tq$, in other words if and only if the associated Lie algebra representation $\pi=D_e\rho:\g\rightarrow gl\left(N,{\Bbb C}\right)$ satisfies
$$\beta_{2n+1}\left(\pi\right)\not=0\in{\Bbb C}\left[\tq\right].$$ 
In this section we will investigate for which fundamental representations of Lie algebras the latter condition is satisfied.\\

In the following subsections we will 
%discuss {\bf complex-linear representations}, that is we will 
consider complex
simple Lie algebras $\g$
and the ring $R_{\Bbb C}\left(\g\right)\subset
R\left(\g\right)$ of their ${\Bbb C}$-linear representations.
The general picture can be reduced to that of ${\Bbb C}$-linear representations in view of the following observations.

{\bf Noncomplex Lie algebras.} Let $\pi:\g\rightarrow gl\left(N,{\Bbb C}\right)$ be an ${\Bbb R}$-linear representation of
a simple Lie-algebra $\g$ which is not a complex Lie algebra. Then $\g\otimes{\Bbb C}$ is a simple complex Lie algebra and $\pi$ is the restriction of some ${\Bbb C}$-linear representation $\g\otimes{\Bbb C}\rightarrow gl\left(N,{\Bbb C}\right)$. Let $\tq$ be a Cartan subalgebra of $\g$. Then it is obvious that
an element $t\in\tq\otimes{\Bbb C}$ with $$Tr\left(\pi\left(t\right)^n\right)\not=0$$ exists if and only if such an element exists in $\tq$. Thus $\pi$ has nontrivial Borel class if and only if the ${\Bbb C}$-linear representation $\pi\otimes{\Bbb C}$ has nontrivial Borel class.
%, and
%we can use the results for ${\Bbb C}$-linear representations.

{\bf ${\Bbb R}$-linear representations of complex Lie algebras.} If $\g$ is a simple complex
Lie algebra, then each ${\Bbb R}$-linear representation
$\pi:\g\rightarrow gl\left(N,{\Bbb C}\right)$ is of the
form $\pi=\pi_1\otimes\overline{\pi_2}$ for ${\Bbb C}$-linear representations $\pi_1,\pi_2$. The equality
$$Tr\left(\pi\left(t\right)^n\right)=Tr\left(\pi_1\left(t\right)^n\right) Tr\left(\overline{\pi_2}\left(t\right)^n\right).$$
implies that ${\Bbb R}$-linear representations with $b_{2n-1}\left(\pi\right)$ exist only if there
are ${\Bbb C}$-linear ones.\\

In the sequel we will go through the fundamental representations of simple Lie algebras and discuss whether their Borel class is nontrivial. The results will be subsumed in Section 3.3 in the proof of \hyperref[Thm3]{Theorem \ref*{Thm3}}. For faster reading we are going to highlighten the exceptional cases that will occur in the proof of \hyperref[Thm3]{Theorem \ref*{Thm3}}.

\subsubsection{$\g=sl\left(l+1,{\Bbb C}\right)$}
Let $V={\Bbb C}^{l+1}$ be the standard representation, with basis $e_1,\ldots,e_{l+1}$. 
	Then $$R_{\Bbb C}\left(\g\right)=
	{\Bbb Z}\left[A_1,\ldots,A_l\right]$$ with $A_k$ 
	the induced representation on $\Lambda^kV$, cf.\ \cite[p.377]{fh}. In particular, irreducible representations occur as representations of dominant weight in tensor products of the fundamental representations $A_1,\ldots,A_l$.
%correspond to monomials with coefficient 1, i.e.\ to tensor products of $A_k$'s.
		We compute $\beta_{2n-1}$ on the fundamental representations $A_k, k=1,\ldots,l$. 

A basis of $\Lambda^kV$ is given by 
$$\left\{e_{i_1}\wedge
		\ldots\wedge e_{i_k}:1\le i_1<\ldots<i_k\le l+1\right\}.$$
		As Cartan-subalgebra 
	we may choose the diagonal matrices
$$\tq=\left\{diag\left(h_1,\ldots,h_l,h_{l+1}\right): h_1+\ldots+h_{l+1}=0\right\}.$$
	$diag\left(h_1,\ldots,h_l,h_{l+1}\right)$ acts on 
		$e_{i_1}\wedge\ldots\wedge 
	e_{i_k}$ by multiplication with $h_{i_1}+\ldots +h_{i_k}$. 
		Hence $$\beta_{2n-1}\left(A_k\right)\left(\begin{array}{cccc}h_1&0&\ldots&0\\0&h_2&\ldots&0\\.&.&\ldots&.\\0&0&\ldots&h_{l+1}\end{array}\right)=\sum_{1\le
		i_1<\ldots<i_k\le l+1}\left(h_{i_1}+\ldots +h_{i_k}\right)^n.$$
If {\bf $k=l=1$} then 
$h_1^n+h_2^n$ is a multiple of $h_1+h_2=0$ if and only if $n$ is odd. Thus for $l=1$ we have $\beta_{2n-1}\left(A_1\right)\not=0$ if $n$ is even and {\bf $\beta_{2n-1}\left(A_1\right)=0$ if $n$ is odd}.

If $k=1$ and $l\ge 2$, then
$\sum_{i=1}^{l+1}h_i^{n}$ does not vanish for example for $h_1=2,h_2=-1,h_3=\ldots=h_l=0,h_{l+1}=-1$. 

If {\bf $2\le k \le l$ and $n=1$}, then $$\sum_{1\le
                i_1<\ldots<i_k\le l+1}\left(h_{i_1}+\ldots +h_{i_k}\right)=\left(\begin{array}{c}l\\
k-1\end{array}\right)\left(h_1+\ldots+h_{l-1}\right)=0,$$ 
thus {\bf $\beta_1\left(A_k\right)=0$}. 

If $2\le k\le l$ and $n>1$, then $\beta_{2n-1}\left(A_k\right)\not=0$. Indeed, nontriviality 
%of
%$\sum_{1\le
%                i_1<\ldots<i_k\le l+1}\left(h_{i_1}+\ldots +h_{i_k}\right)^n$ 
can be seen for example
by considering again the diagonal matrix $\left(2,-1,0,\ldots,0,-1\right)\in\tq$, for which we obtain $$\sum_{1\le
                i_1<\ldots<i_k\le l+1}\left(h_{i_1}+\ldots +h_{i_k}\right)^n= \left(2^n-1\right)\left(
        \left(\begin{array}{c}l-2\\k-1\end{array}\right)-
        \left(\begin{array}{c}l-1\\k-1\end{array}\right)\right)<0.$$

		%For $n=1$, the sum is a
		%multiple of $h_1+\ldots+h_{l+1}=0$. For $l=k=1$ and $n$ odd, we have that $b_n\left(A_1\right)=h_1^n+h_2^n$ is a
		%multiple of $h_1+h_2=0$.
		%In all other cases, i.e.\ for $l\ge 2, n>1$ or $l=k=1, n$ even, the sum is not divisible by 
		%$h_1+\ldots+h_{l+1}$ and thus not trivial. This is obvious for even $n$. In the case of odd $n>1$ and $l\ge 2$, it follows for example from the computation $\beta_{2n-1}\left(A_k\right)diag\left(2,-1,-1,0,\ldots,0\right)=\left(2^n-1\right)\left(
%	\left(\begin{array}{c}l-2\\k-1\end{array}\right)-
%	\left(\begin{array}{c}l-1\\k-1\end{array}\right)\right)\not=0$.
%	Thus, for the fundamental representation $A_k$ of $sl\left(l+1,{\Bbb C}\right)$ we 
%have $\beta_{2n-1}\left(A_k\right)$ if and only if either $l\ge 2, n>1$ or $l=k=1, n$ even.
\noindent
{\bf Conclusion:} {\em The exceptional cases with $\beta_{2n-1}\left(A_k\right)=0$ occur for \\
- $k=l=1$, $n$ odd,\\
- $2\le k\le l$, $n=1$.}

\subsubsection{         $\g=so\left(2l,{\Bbb C}\right)$}
 Let $V={\Bbb C}^{2l}$ with ${\Bbb C}$-basis $e_1,\ldots,e_l,f_1,\ldots,f_l$.
Let $Q$ be the quadratic form given by $Q\left(e_i,f_i\right)=Q\left(f_i,e_i\right)=1$ for $i=1,\ldots,l$, 
$Q\left(e_i,f_j\right)=Q\left(f_i,e_j\right)=0$ for $i\not=j$ and
$Q\left(e_i,e_j\right)=Q\left(f_i,f_j\right)=0$ for all $i,j=1,\ldots,l$.

Following \cite[p.268 ff.]{fh} we consider $so\left(2l,{\Bbb C}\right)$ as the skew-symmetric matrices with respect to the quadratic
form $Q:V\times V\rightarrow {\Bbb C}$. (All quadratic forms are equivalent over ${\Bbb C}$ under a suitable change of base, the corresponding Lie groups $SO\left(Q\right)\subset GL\left(N,{\Bbb C}\right)$ are conjugate, thus it is sufficient to consider the Lie algebra $so\left(Q\right)$ with respect to this quadratic form $Q$.)

Let $D_1:so\left(2l,{\Bbb C}\right)\rightarrow gl\left(V\right)$ be the standard representation.

                Then $$R_{\Bbb C}\left(\g\right)={\Bbb Z}\left[D_1,\ldots,D_{l-2},S^+,S^-\right]$$
                with
$D_k:so\left(2l,{\Bbb C}\right)\rightarrow gl\left(\Lambda^kV\right)$ the representation induced from $D_1$
on $\Lambda^kV$,
and $S^\pm$ the half-spinor
                representations.

                As a Cartan-subalgebra we may choose the diagonal matrices
                $$\tq=\left\{diag\left(h_1,\ldots,h_l,-h_1,\ldots,-h_l\right): h_1,\ldots,h_l\in{\Bbb C}\right\}.$$

First we look at $\beta_{2n-1}\left(D_k\right)$ for the fundamental representations $D_k$.

 A basis of $\Lambda^kV$ is given by
$$\left\{e_{i_1}\wedge\ldots\wedge e_{i_p}\wedge
                f_{j_1}\wedge\ldots\wedge f_{j_{k-p}}:
        0\le p\le k, 1\le i_1<\ldots<i_p\le l,
                1\le j_1<\ldots <j_{k-p}\le l\right\}.$$

                $diag\left(h_1,\ldots,h_l,-h_1,\ldots,-h_l\right)$ acts
        on $e_{i_1}\wedge\ldots\wedge e_{i_p}\wedge
                f_{j_1}\wedge\ldots\wedge f_{j_{k-p}}$
        by multiplication with $h_{i_1}+\ldots+h_{i_p}-h_{j_1}-\ldots-h_{j_{k-p}}$. Hence
                $$\beta_{2n-1}\left(D_k\right)\left(\begin{array}{cccccc}h_1&0&\ldots&0&0&\ldots\\0&h_2&\ldots&0&0&\ldots\\
                \ldots& & &\ldots& &\\
                0&0&\ldots&-h_1&0&\ldots\\0&0&\ldots&0&-h_2&\ldots\\
        \ldots& & &\ldots& &\\
        \end{array}\right)$$
$$=\sum_{1\le
                i_1<\ldots<i_p\le l, 1\le j_1<\ldots <j_{k-p}\le l}
                \left(h_{i_1}+\ldots+h_{i_p}-h_{j_1}-\ldots-h_{j_{k-p}}\right)^n.$$
        If $n$ is even, then we get a nonvanishing polynomial. This follows from \hyperref[even]{Corollary \ref*{even}} or
more explicitly for example from $$\beta_{2n-1}\left(D_k\right)\left(diag\left(1,0,\ldots,0,-1,0,\ldots,0\right)\right)>0.$$

%        For odd $n$ we claim $\beta_{2n-1}\left(D_k\right)=0$ for all $k$. This can be seen as follows.
%Let us consider the involution $B\in gl\left(2l,{\Bbb C}\right)$ given by $B\left(e_i\right)=f_i, B\left(f_i\right)=e_i$ for $i=1,\ldots,l$. Clearly $$\left\{Be_{i_1}\wedge\ldots\wedge Be_{i_p}\wedge
%                Bf_{j_1}\wedge\ldots\wedge Bf_{j_{k-p}}:
%        0\le p\le k, 1\le i_1<\ldots<i_p\le l,
%                1\le j_1<\ldots <j_{k-p}\le l\right\}$$
%is a
%and $D_k\left(H\right)^n \left(

%The involution $B$ preserves the quadatic form $Q$, thus 
%it preserves the trace, that is $Tr\left(D_k\left(H\right)^n\right)=\sum_{i=1}^l Q\left(D_k\left(H\right)^ne_i,e_i\right)
%+Q\left(D_k\left(H\right)^nf_i,f_i\right)=\sum_{i=1}^l \sum_{i=1}^l Q\left(D_k\left(H\right)^nBe_i,Be_i\right)
%+Q\left(D_k\left(H\right)^nBf_i,Bf_i\right).$
%Let $H=diag\left(h_1,\ldots,h_l,-h_1,\ldots,-h_l\right)\in \tq$. Then $H\left(Be_i\right)  we have $B^{-1}HB=diag\left(-h_1,\ldots,-h_l,h_1,\ldots,h_l\right) \in\tq$ and the above formula implies for $n$ odd:
%$$\beta_{2n-1}\left(D_k\right)\left(B^{-1}HB\right)=-\beta_{2n-1}\left(D_k\right)\left(H\right).$$
%On the other hand
%$\beta_{2n-1}\left(D_k\right)\left(B^{-1}HB\right)=Tr\left(D_k\left(B^{-1}HB\right)^n\right)$
If $n$ is odd, then the permutation, which transposes $i_r$
and $j_r$ simultaneously for all $r$, multiplies the sum by $-1$, but on the other
hand preserves the sum. Thus {\bf $\beta_{2n-1}\left(D_k\right)=0$ if $n$ is odd}.\\

Next we look at $\beta_{2n-1}\left(S^\pm\right)$ for the half-spinor representations $S^\pm$.

Let $TV=\oplus_{k=0}^mV^{\otimes k}$ be the tensor algebra of $V$ and let $Cl\left(Q\right)=TV/I\left(Q\right)$
be the Clifford algebra of $Q$, where $I\left(Q\right)$ is the
ideal generated by all $v\otimes v+Q\left(v,v\right)1$
with $v\in V$. The grading of $\oplus_{k=0}^mV^{\otimes k}$ induces a well-defined ${\Bbb Z}/2{\Bbb Z}$-grading $Cl\left(Q\right)=Cl\left(Q\right)^{even}\oplus Cl\left(Q\right)^{odd}$ on the Clifford algebra.

%Let $so\left(Q\right)$ be the Lie algebra of skew-adjoint matrices with respect to $Q$. All quadratic forms over ${\Bbb C}$ are equivalent under a suitable change of base, thus there is an isomorphism $I:so\left(m,{\Bbb C}\right)\rightarrow so\left(Q\right)$ and it will be sufficient to describe (following \cite[p.268 ff.]{fh})
%the representations $S^\pm=D_e\rho^\pm \circ I^{-1}:so\left(Q\right)\rightarrow GL\left(N,{\Bbb C}\right)$ and to prove $Tr\left(S^\pm\left(t\right)^n\right)\not=0$ for some $t\in\tq^Q$, where $\tq^Q$ is a Cartan subalgebra of $so\left(Q\right)$.
%%and let $Cl\left(Q\right)=Cl\left(Q\right)^{even}\oplus Cl\left(Q\right)^{odd}=TV/I\left(Q\right)$ be the Clifford algebra of $Q$ with the grading induced from the grading of the tensor algebra $TV=\oplus_{k=0}^mV^{\otimes k}$, where $I\left(Q\right)$ is the generated by all $v\otimes v+Q\left(v,v\right)1, v\in V$.
% As a
%        Cartan-subalgebra $\tq^Q\subset so\left(Q\right)$ we choose (following \cite[p.270]{fh}) the subalgebra of diagonal matrices
%$$\tq^Q=\left\{diag\left(h_1,\ldots,h_l,-h_1,\ldots,-h_l\right): h_1,\ldots, h_l\in{\Bbb C}\right\}.$$
Denote by $E_{ij}$ the elementary matrix with entry $1$ at position $\left(i,j\right)$ and entries $0$ else. Then
$$\left\{A_i:=E_{i,i}-E_{l+i,l+i}, i=1,\ldots,l\right\}$$ is a basis of $\tq$.

By \cite[pp.303-305]{fh}, there is an injective homomorphism $$\iota: so\left(
        Q\right)\rightarrow Cl\left(Q\right)^{even}$$ which maps, in particular, $
        A_i$ to $\frac{1}{2}\left(e_i\otimes f_i-1\right)$.

Let $W$ be the ${\Bbb C}$-subspace of $V$ spanned by $e_1,\ldots,e_l$.
%, $W^\prime$ the subspace spanned by $e_{n+1},\ldots,e_{2n}$.

From the proof of \cite[Lemma 20.9]{fh} we have a homomorphism $$\Phi:Cl\left(Q\right)\rightarrow gl\left(\Lambda^* W\right)$$
with  $$\Phi\left(e_i\right)\left(v_1\wedge\ldots\wedge v_k\right)=e_i\wedge v_1\wedge\ldots\wedge v_k$$
$$\Phi\left(f_i\right)\left(v_1\wedge\ldots\wedge v_k\right)=\sum_{j=1}^k\left(-1\right)^{j-1}2Q\left(v_j,f_i\right)v_1\wedge\ldots \widehat{v_j}\ldots \wedge v_k$$ for all $v_1\wedge\ldots\wedge v_k\in \Lambda^*W$ and
%sends $v\in W$ to $e_i\wedge v$ and $e_{n+i}$
%       sends $v\in W$ to $2v-2Q\left(v,f_i\right)e_i$, for
$i=1,\ldots,l$,
%(This follows from the proof of \cite[Lemma 20.9]{fh}.)
%       This action
%       extends in the obvious way to an action of $Cl\left(Q\right)$ on $\Lambda^* W$.
%       In particular $\Phi\left(\frac{1}{2}\left(e_ie_{n+i}-1\right)\right)\left(v\right)=
%       e_i\wedge v-\frac{1}{2}v$ for all $v\in W$,
which implies
$$\Phi\left(\frac{1}{2}\left(e_i\otimes f_i-1\right)\right)\left(e_{i_1}\wedge\ldots\wedge e_{i_k}\right)=
\frac{1}{2}e_{i_1}\wedge\ldots\wedge e_{i_k}$$ if $i\in\left\{i_1,\ldots,i_k\right\}$
and $$\Phi\left(\frac{1}{2}\left(e_i\otimes f_i-1\right)\right)\left(e_{i_1}\wedge\ldots\wedge e_{i_k}\right)=-\frac{1}{2}e_{i_1}\wedge\ldots\wedge e_{i_k}$$ if $i\not
\in\left\{i_1,\ldots,i_k\right\}$.

By \cite[p.305]{fh},
restriction of $\Phi$ to $Cl\left(Q\right)^{even}$
gives
        rise to an isomorphism $$\Phi^{even}:Cl\left(Q\right)^{even}\rightarrow
        End\left(\Lambda^{even}W\right)\oplus End\left(\Lambda^{odd}W\right).$$ Let $\pi_1,\pi_2$ be the projections from $
End\left(\Lambda^{even}W\right)\oplus End\left(\Lambda^{odd}W\right)$
to the first resp.\ second summand.
The induced homomorphisms
$$S^+:=\pi_1 \Phi^{even} \iota :so\left(Q\right)\rightarrow End\left(\Lambda^{even}W\right)$$
$$S^-:=\pi_2 \Phi^{even} \iota:so\left(Q\right)\rightarrow End\left(\Lambda^{odd}W\right)$$
give the
positive resp.\ negative half-spinor
        representations that we are going to consider.

 Thus $$
        S^\pm\left(A_i\right)\left(
e_{i_1}\wedge\ldots \wedge e_{i_k}\right)=
        \frac{1}{2}
e_{i_1}\wedge\ldots \wedge e_{i_k}$$ if $i\in\left\{i_1,\ldots,i_k\right\}$ and
        $$
 S^\pm\left(A_i\right)\left(
e_{i_1}\wedge\ldots \wedge e_{i_k}\right)=          -\frac{1}{2}
e_{i_1}\wedge\ldots \wedge e_{i_k}$$ if $i\not\in\left\{i_1,\ldots,i_k\right\}$.

        For the positive half-spinor representation  $S^+$ and any $n\in{\Bbb N}$ we obtain
$$\beta_{2n-1}\left(S^+\right)\left(\begin{array}{cccccc}h_1&0&\ldots&0&0&\ldots\\0&h_2&\ldots&0&0&\ldots\\
                \ldots& & &\ldots& &\\
                0&0&\ldots&-h_1&0&\ldots\\0&0&\ldots&0&-h_2&\ldots\\
        \ldots& & &\ldots& &\\
        \end{array}\right)$$        
$$=\frac{1}{2^n}\sum_{0\le k\le l, k\ even}\sum_{
\mid I\mid=k}\left(\sum_{i\in I}h_i-\sum_{j\not\in I}h_j\right)^n.$$

If $n$ is even, then $\beta_{2n-1}\left(S^+\right)\not=0$ follows from \hyperref[even]{Lemma \ref*{even}}. 

If {\bf $n$ is odd and $l$ is even}, then for each $I$ with $k=\mid I\mid$ even we have $I^\prime:=\left\{1,\ldots,l\right\}-I$ with $k^\prime=\mid I^\prime\mid$ even and $
\left(\sum_{i\in I}h_i-\sum_{j\not\in I}h_j\right)^n$ cancels against 
$
\left(\sum_{i\in I^\prime}h_i-\sum_{j\not\in I^\prime}h_j\right)^n$. Thus all summands cancel and {\bf $\beta_{2n-1}\left(S^+\right)=0$}.

We prove that the polynomial is nontrivial for all $n\ge l$ with $n\equiv l\ mod\ 2$, in particular if $n$ and $l$ are both odd.
It suffices to show that for example the coefficient of $h_1^{n-l+1}h_2\ldots h_n$ is not zero. 
First we observe that the coefficient of $h_1^{n-l+1}h_2\ldots h_n$ in $\left(\sum_{i\in I}h_i-\sum_{j\not\in I}h_j\right)^n$ 
is $\frac{n!}{\left(n-l+1\right)!}\left(-1\right)^{n-k}$ if $1\in I$ resp.\ 
$\frac{n!}{\left(n-l+1\right)!}\left(-1\right)^{l-k}$ if $1\not\in I$. Thus the coefficient of $h_1^{n-l+1}h_2\ldots h_n$ in $\sum_{
\mid I\mid=k}\left(\sum_{i\in I}h_i-\sum_{j\not\in I}h_j\right)^n$ is $$\frac{n!}{\left(n-l+1\right)!} \left(\left(\begin{array}{c}l-1\\
k-1\end{array}\right)\left(-1\right)^{n-k}+\left(\begin{array}{c}l-1\\
k\end{array}\right)\left(-1\right)^{l-k}\right).$$ All summands have the same sign because of $n\equiv l\ mod\ 2$. Thus $\beta_{2n-1}\left(S^+\right)\not=0$.

%$$Tr\left(
%        S^+\left(A_i\right)^n\right)=\sum_{0\le k\le l, k\ even}
%\left(  \sum_{i\in\left\{i_1,\ldots,i_k\right\}} \left(\frac{1}{2}\right)^n
%        +\sum_{i\not\in\left\{i_1,\ldots,i_k\right\}} \left(-\frac{1}{2}\right)^n
%        \right)$$
%        $$=\sum_{0\le k\le l,
%k\ even}\left( \left(\begin{array}{c}l-1\\
%        k-1\end{array}\right)\left(\frac{1}{2}\right)^n +
%         \left(\begin{array}{c}l-1\\
%        k\end{array}\right)\left(-\frac{1}{2}\right)^n\right).$$
%        In particular, we have $Tr\left(
%        S^+\left(A_i\right)^n\right)=0$ for $n$ odd, $l\ge 2$, and $Tr\left(
%        S^+\left(A_i\right)^n\right)>0$ for $n$ even. \\

        For the negative half-spinor representation $S^-$ and any $n\in{\Bbb N}$
we obtain
        $$\beta_{2n-1}\left(
        S^-\right)
\left(\begin{array}{cccccc}h_1&0&\ldots&0&0&\ldots\\0&h_2&\ldots&0&0&\ldots\\
                \ldots& & &\ldots& &\\
                0&0&\ldots&-h_1&0&\ldots\\0&0&\ldots&0&-h_2&\ldots\\
        \ldots& & &\ldots& &\\
        \end{array}\right)$$
$$=\frac{1}{2^n}\sum_{0\le k\le l, k\ odd}\sum_{
\mid I\mid=k}\left(\sum_{i\in I}h_i-\sum_{j\not\in I}h_j\right)^n.$$
If $n$ is even, then $\beta_{2n-1}\left(S^-\right)\not=0$ by \hyperref[even]{Lemma \ref*{even}}. 

If $n$ is odd, then the same argument as in the computation of $\beta_{2n-1}\left(D_k\right)=0$ shows that $$\beta_{2n-1}\left(S^+\right)+\beta_{2n-1}\left(S^-\right)=0,$$ thus $ \beta_{2n-1}\left(S^+\right)\not =0$ implies $\beta_{2n-1}\left(S^-\right)\not =0$ if $n\ge l$ and $n\equiv l\ mod\ 2$.

Conversely, if {\bf $n$ is odd and $l$ is even}, then
{\bf $\beta_{2n-1}\left(S^-\right)=0$}.\\
\\
{\bf Conclusion:} {\em The cases with $\beta_{2n-1}\left(\pi\right)=0$ are precisely:\\
- $\pi=D_k$ ($1\le k\le l-2$), $n$ odd,\\
- $\pi=S^\pm$, $l$ even, $n$ odd.}

\subsubsection{		$\g=so\left(2l+1,{\Bbb C}\right)$}
Let $V={\Bbb C}^{2l+1}$ with ${\Bbb C}$-basis $e_1,\ldots,e_l,f_1,\ldots,f_l,g$,
        and $Q$ the quadratic form given by $Q\left(g,g\right)=1,
        Q\left(e_i,f_i\right)=Q\left(f_i,e_i\right)=1$ for $i=1,\ldots,l$,
        and $Q\left(.,.\right)=0$ for all other pairs of basis vectors.

Following \cite[p.268 ff.]{fh} we consider $so\left(2l+1,{\Bbb C}\right)$ as the skew-symmetric matrices with respect to the quadratic
form $Q:V\times V\rightarrow {\Bbb C}$.
		Let $C_1:so\left(2l+1,{\Bbb C}\right)\rightarrow gl\left(V\right)$ be the standard representation.
		Then $$R_{\Bbb C}\left(\g\right)={\Bbb Z}\left[C_1,\ldots,C_{l-1},S\right]$$
		with $C_k:so\left(2l+1,{\Bbb C}\right)\rightarrow gl\left(\Lambda^kV\right)$ the representation induced from $C_1$
on $\Lambda^kV$, and $S$ the spinor representation.

		As a Cartan-subalgebra we may choose the diagonal matrices
$$\tq=\left\{diag\left(h_1,\ldots,h_l,-h_1,\ldots,-h_l,0\right): h_1,\ldots,h_l\in{\Bbb C}\right\}.$$
		Then the computation of $\beta_{2n-1}$ on $C_k$ 
	is exactly the same as for $so\left(2l,{\Bbb C}\right)$ and $D_k$, in particular $\beta_{2n-1}\left(C_k\right)\not=0$ for $n$ even
	and {\bf $\beta_{2n-1}\left(C_k\right)=0$ for $n$ odd}. 

We look at $\beta_{2n-1}\left(S\right)$ for the spinor representation $S$.
As in the case of $so\left(2l,{\Bbb C}\right)$,
        we have $\iota: so\left(                                        Q\right)
        \rightarrow Cl\left(Q\right)^{even}$ with $
\iota\left(E_{i,i}-E_{l+i,l+i}\right)=\frac{1}{2}\left(e_i\otimes f_i-1\right).$

          Let $W$ be the ${\Bbb C}$-subspace of $V$ spanned by $e_1,\ldots,e_{l}$.
%, $W^\prime$ the subspace spanned by $e_{n+1},\ldots,e_{2n}$.
        It follows from the proof of \cite[Lemma 20.16]{fh} that $Cl\left(Q\right)$ acts on $\Lambda^*W$ as follows:
        the action of $e_i$ resp.\ $f_i$, for $i=1,\ldots,l$ is defined as in the case of $so\left(2l,{\Bbb C}\right)$, and 
$g$ acts as
        multiplication by 1 on $\Lambda^{even}W$ and as multiplication by -1 on $\Lambda^{odd}W$.
        In particular, we have again that $\frac{1}{2}\left(e_i\otimes f_i-1\right)$ acts by
        sending
        $e_{i_1}\wedge\ldots\wedge e_{i_k}$ to
                $\frac{1}{2}e_{i_1}\wedge\ldots\wedge e_{i_k}$ if $i\in\left\{i_1,\ldots,i_k\right\}$ resp.\ to $-
\frac{1}{2}e_{i_1}\wedge\ldots\wedge e_{i_k}$ if $i\not\in\left\{i_1,\ldots,i_k\right\}$.

        This action gives
        rise to an isomorphism $Cl\left(Q\right)^{even}\cong
        End\left(\Lambda W\right)$ (see \cite[p.306]{fh}).
        The induced action of $so
        \left(Q\right)$ on $\Lambda
        W$
        is the spinor
        representation $S$.

%        As a
%        Cartan-algebra for $so\left(Q\right)$
%we choose the algebra of diagonal matrices
%$$\tq^Q=\left\{diag\left(h_1,\ldots,h_{l},-h_1,\ldots,-h_{l},0\right), h_1,\ldots,h_{l}\in{\Bbb C}\right\}.$$
        Let $\left\{A_i:i=1,\ldots,{l}\right\}$ be a basis of $\tq$, where $$A_i=E_{i,i}-E_{l+i,l+i}.$$
        $A_i$
        acts on $e_{i_1}\wedge\ldots \wedge e_{i_k}$ by multiplication with
        $\frac{1}{2}$ if $i\in\left\{i_1,\ldots,i_k\right\}$ and
        by multiplication with $-\frac{1}{2}$ if $i\not\in\left\{i_1,\ldots,i_k\right\}$.
        Thus we obtain for any $n\in {\Bbb N}$:
$$\beta_{2n-1}\left(
        S\right)
\left(\begin{array}{cccccc}h_1&0&\ldots&0&0&\ldots\\0&h_2&\ldots&0&0&\ldots\\
                \ldots& & &\ldots& &\\
                0&0&\ldots&-h_1&0&\ldots\\0&0&\ldots&0&-h_2&\ldots\\
        \ldots& & &\ldots& &\\
        \end{array}\right)$$
$$=\frac{1}{2^n}\sum_{0\le k\le l}\sum_{
\mid I\mid=k}\left(\sum_{i\in I}h_i-\sum_{j\not\in I}h_j\right)^n.$$

        Thus, by the same argument as for $D_k$ and $C_k$, 
{\bf $\beta_{2n-1}\left(
        S\right)
        =0$ for $n$ odd} and 
$\beta_{2n-1}\left(
        S\right)
        \not=0$ for $n$ even.\\
\\
{\bf Conclusion:} {\em The cases with $\beta_{2n-1}\left(\pi\right)=0$ are precisely:\\
- $\pi=C_k$ ($1\le k\le l-1$), $n$ odd,\\
- $\pi=S$, $n$ odd.}

 \subsubsection{$\g=sp\left(l,{\Bbb C}\right)$}
                Let $V={\Bbb C}^{2l}$ with basis $\left\{e_1,\ldots,e_l,f_1,\ldots,f_l\right\}$. Consider the symplectic form $Q:V\times V\rightarrow{\Bbb R}$ given by $Q\left(e_i,f_i\right)=1=-Q\left(f_i,e_i\right)$ for $i=1,\ldots,l$, and $Q\left(.,.\right)=0$ for each other pair of basis vectors. Let $Sp\left(l,{\Bbb C}\right)$ be the Lie group of linear maps preserving this symplectic form. Then its lie algebra
$sp\left(l,{\Bbb C}\right)$ consists of matrices $\left(\begin{array}{cc}A&B\\
        C&D\end{array}\right)$, such that the $l$x$l$-blocks $A,B,C,D$
        satisfy $B^T=B,C^T=C,A^T=-D$.
                As a Cartan-subalgebra we may choose the diagonal matrices
$$\tq=\left\{diag\left(h_1,\ldots,h_l,-h_1,\ldots,-h_l\right): h_1,\ldots,h_l\in{\Bbb C}\right\}.$$

                Then $$R_{\Bbb C}\left(\g\right)={\Bbb Z}\left[B_1,\ldots,B_l\right]$$
                where by \cite[p.377]{fh} the fundamental representations $B_k$ are the induced representations of 
$sp\left(l,{\Bbb C}\right)$ on $ker\left(\phi_k:\Lambda^kV\rightarrow \Lambda^{k-2}V\right)$ for $k=1,\ldots,l$, where $\phi_k$ is the contraction using $Q$ defined in \cite[p.260]{fh}
by
$$\phi_k\left(v_1\wedge\ldots\wedge v_k\right)=\sum_{i<j}Q\left(v_i,v_j\right)\left(-1\right)^{i+j-1} v_1\wedge\ldots\wedge\hat{v}_i\wedge\ldots\wedge\hat{v}_j\wedge\ldots\wedge v_k.$$
We consider $\beta_{2n-1}$ for the fundamental representations $B_k, k=1,\ldots,l$.

If $n$ is even, then $\beta_{2n-1}\left(B_k\right)\not=0$ follows from \hyperref[even]{Corollary \ref*{even}}.

We claim that {\bf $\beta_{2n-1}\left(B_k\right)=0$ if $n$ is odd}. This can be seen as follows.
Consider the involution $B\in gl\left(2l,{\Bbb C}\right)$ given by $B\left(e_i\right)=f_i, B\left(f_i\right)=-e_i$ for $i=1,\ldots,l$. It induces an involution on $\Lambda^kV$. 
$B$ preserves the symplectic form $Q$, thus we have 
$$\phi_k\left(Bv_1\wedge \ldots\wedge Bv_k\right)=\sum_{i<j}Q\left(Bv_i,B_j\right)\left(-1\right)^{i+j-1} Bv_1\wedge\ldots\wedge\hat{Bv}_i\wedge\ldots\wedge\hat{Bv}_j\wedge\ldots\wedge Bv_k$$
$$=\sum_{i<j}Q\left(v_i,v_j\right)\left(-1\right)^{i+j-1} Bv_1\wedge\ldots\wedge\hat{Bv}_i\wedge\ldots\wedge\hat{Bv}_j\wedge\ldots\wedge Bv_k=B\left(\phi_k\left(v_1\wedge\ldots\wedge v_k\right)\right),$$
in particular $B$ maps $ker\left(\phi_k\right)$ to itself. If $\left\{b_1,\ldots,b_{dim\left(ker\left(\phi_k\right)\right)}\right\}$
is a basis of $dim\left(ker\left(\phi_k\right)\right)$, then
 $\left\{Bb_1,\ldots,Bb_{dim\left(ker\left(\phi_k\right)\right)}\right\}$
is a basis of $dim\left(ker\left(\phi_k\right)\right)$.

Let $<.,.>$ be the standard scalar product on ${\Bbb C}^{2l}$ such that $\left\{e_1,\ldots,e_l,f_1,\ldots,f_l\right\}$ is an orthonormal basis. We note that $B$ preserves this scalar product. Thus, if $\left\{b_1,\ldots,b_{dim\left(ker\left(\phi_k\right)\right)}\right\}$ is an orthonormal basis of $ker\left(\phi_k\right)\subset{\Bbb C}^{2l}$, then 
$\left\{Bb_1,\ldots,Bb_{dim\left(ker\left(\phi_k\right)\right)}\right\}$ is an orthonormal basis of 
$ker\left(\phi_k\right)$ as well and we have
$$Tr\left(B_k\left(H\right)^n\right)=\sum_{i=1}^{dim\left(ker\left(\phi_k\right)\right)} <B_k\left(H\right)^n b_i,b_i>=
\sum_{i=1}^{dim\left(ker\left(\phi_k\right)\right)} <B_k\left(H\right)^n Bb_i,Bb_i>$$
for each $H\in sp\left(l,{\Bbb C}\right)$.

On the other hand, for $H=diag\left(h_1,\ldots,h_l,-h_1,\ldots,-h_l\right)\in\tq\subset sp\left(l,{\Bbb C}\right)$
and $n$ odd we have $$<B_k\left(H\right)^ne_i,e_i>=h_i^n, <B_k\left(H\right)^nf_i,f_i>=-h_i^n$$
for $i=1,\ldots,n$, which implies $$<B_k\left(H\right)^ne_i,e_i>= -<B_k\left(H\right)^nBe_i,Be_i>, <B_k\left(H\right)^nf_i,f_i>= -<B_k\left(H\right)^nBf_i,Bf_i>.$$
From bilinearity of the scalar product we conclude 
$$<B_k\left(H\right)^nv,v>= -<B_k\left(H\right)^nBv,Bv>$$
for all $v\in{\Bbb C}^{2l}$, in particular for $v=b_1,\ldots,b_{dim\left(ker\left(\phi_k\right)\right)}\in ker\left(\phi_k\right)$.
Thus $$\sum_{i=1}^{dim\left(ker\left(\phi_k\right)\right)} <B_k\left(H\right)^n b_i,b_i>=
-\sum_{i=1}^{dim\left(ker\left(\phi_k\right)\right)} <B_k\left(H\right)^n Bb_i,Bb_i>,$$
which implies $$Tr\left(B_k\left(H\right)^n\right)=\sum_{i=1}^{dim\left(ker\left(\phi_k\right)\right)} <B_k\left(H\right)^n b_i,b_i>=0.$$

%Clearly $$\left\{Be_{i_1}\wedge\ldots\wedge Be_{i_p}\wedge
%                Bf_{j_1}\wedge\ldots\wedge Bf_{j_{k-p}}:
%        0\le p\le k, 1\le i_1<\ldots<i_p\le l,
%                1\le j_1<\ldots <j_{k-p}\le l\right\}$$
%is a
%and $D_k\left(H\right)^n \left(

%The involution $B$ preserves the quadatic form $Q$, thus
%it preserves the trace, that is $Tr\left(D_k\left(H\right)^n\right)=\sum_{i=1}^l Q\left(D_k\left(H\right)^ne_i,e_i\right)
%+Q\left(D_k\left(H\right)^nf_i,f_i\right)=\sum_{i=1}^l \sum_{i=1}^l Q\left(D_k\left(H\right)^nBe_i,Be_i\right)
%+Q\left(D_k\left(H\right)^nBf_i,Bf_i\right).$
%Let $H=diag\left(h_1,\ldots,h_l,-h_1,\ldots,-h_l\right)\in \tq$. Then $H\left(Be_i\right)  we have $B^{-1}HB=diag\left(-h_1,\ldots,-h_l,h_1,\ldots,h_l\right) \in\tq$ and the above formula implies for $n$ odd:
%$$\beta_{2n-1}\left(D_k\right)\left(B^{-1}HB\right)=-\beta_{2n-1}\left(D_k\right)\left(H\right).$$
%On the other hand
%$\beta_{2n-1}\left(D_k\right)\left(B^{-1}HB\right)=Tr\left(D_k\left(B^{-1}HB\right)^n\right)$
\noindent
{\bf Conclusion:} {\em The cases with $\beta_{2n-1}\left(B_k\right)=0$ are precisely:\\
- $1\le k\le l$, $n$ odd.}

		\subsubsection{Exceptional Lie groups}
For the applications of \hyperref[Thm2]{Theorem \ref*{Thm2}} and \hyperref[Thm3]{Theorem \ref*{Thm3}} we will have to 
consider only odd-dimensional manifolds and therefore we are only interested in Lie groups which admit a symmetric space of odd dimension. 
The only exceptional Lie group admitting an odd-dimensional symmetric 
space is $E_7$ with $dim\left(E_7/E_7\left({\Bbb R}\right)\right)=163$. 
The fact that $163\equiv\ 3\ mod\ 4$ implies by \hyperref[even]{Corollary \ref*{even}} that $\rho^*b_{163}\not=0$ holds for each irreducible representation $\rho$. 

For completeness we also show, at least for a specific representation, that
$\rho^*b_{2n-1}\not=0$ holds for each $n\ge 6$.
Namely, we consider the representation $\rho:E_7\rightarrow GL\left(56,{\Bbb C}\right)$, which has been constructed in \cite[Corollary 8.2]{ad}, and we are going to show that this representation satisfies $\rho^*b_{2n-1}\not=0$ for each $n\ge 6$, in particular for $n=82$.

	By \cite[Chapter 7/8]{ad} there is a monomorphism $Spin\left(12\right)\times SU\left(2\right)/{\Bbb Z}_2\rightarrow E_7$
	and the Cartan-subalgebra of the Lie algebra $e_7$ coincides with the Cartan-subalgebra $t$ of $spin\left(12\right)\oplus su\left(2\right)$. According to
	\cite[Corollary 8.2]{ad}, the restriction of $\rho$ to $Spin\left(12\right)\times SU\left(2\right)$ is $\lambda_{12}^1\otimes\lambda_1\oplus S^-\otimes 1$, where $\lambda_{12}^1$ resp.\ $\lambda_1$ are the standard representations and $S^-$ is the
	negative spinor representation. 

For even $n$, we know that $\rho^*b_{2n-1}\not=0$.
If $n$ is odd then, for the derivative $\pi_1$ of the standard representation 
	$\lambda_1$ of $SU\left(2\right)$ we have $Tr\left(\pi_1\left(h\right)^n\right)
=0$, whenever $h\in\tq\cap su\left(2\right)$ belongs to the Cartan-subalgebra of 
$su\left(2\right)$, because the latter
are the diagonal 2x2-matrices of trace 0. Thus the first direct summand $\lambda_{12}^1\otimes\lambda_1$ does not contribute to $Tr\left(\pi\left(h\right)^n\right)$. Hence, for $h=\left(h_{spin},h_{su}\right)\in\tq
\subset spin\left(12\right)\oplus su\left(2\right)$,
we have $Tr\left(\pi\left(h\right)^n\right)=Tr\left(S^-\left(h_{spin}\right)^n\right)$. But the nontriviality of the latter has already been shown in 
Section 3.2.3.

\subsection{Conclusion}

In this section, we discuss, for which symmetric spaces $G/K$ (irreducible,
of noncompact type, of dimension $2n-1$) and which representations $\rho:G\rightarrow GL\left(N,{\Bbb C}\right)$ the inequality $\rho^*b_{2n-1}\not=0$ holds. 

%\begin{df}\label{nontrivial2} We say that a Lie algebra representation $\pi:\g\rightarrow gl\left(N,{\Bbb C}\right)$
%has nontrivial Borel class if $\beta_{2n-1}\left(\pi\right)\not\equiv 0$, for $\beta_{2n-1}:R\left(\g\right)\rightarrow{\Bbb C}\left[\tq\right]$ defined in Section 3.2. \end{df}
%\begin{pro}\label{liealgebra} Let $\rho:G\rightarrow GL\left(N,{\Bbb C}\right)$ be a representation of
%a Lie group $G$, and $\pi:\g\rightarrow gl\left(N,{\Bbb C}\right)$ the associated Lie algebra
%representation $\pi=D_e\rho$. Then $\rho$ has nontrivial Borel class if and only if $\pi$ has nontrivial Borel class.\end{pro}
%\begin{pf} This is precisely the statement of \hyperref[cartan]{Lemma \ref*{cartan}}.\end{pf}

\begin{thm}
The following is a complete list of irreducible symmetric spaces $G/K$ of noncompact type and fundamental representations $\rho:G\rightarrow GL\left(N,{\Bbb C}\right)$ with $\rho^*b_{2n-1}\not=0$ for $2n-1:=dim\left(G/K\right)$.\\
\begin{tabular}[ht]{|l|c|}
  \hline
Symmetric Space & Representation\\
  \hline\hline
 $SL_l\left({\Bbb R}\right)/SO_l, l\equiv 0,3,4,7\ mod\ 8$& any fundamental representation\\
 $SL_l\left({\Bbb C}\right)/SU_l, l\equiv 0\ mod\ 2$& any fundamental representation\\
 $SL_{2l}\left({\Bbb H}\right)/Sp_l, l\equiv 0\ mod\ 2$&  any fundamental representation\\
 $Spin_{p,q}/\left(Spin_p\times Spin_q\right), p,q\equiv 1\ mod\ 2, p\not\equiv q\ mod\ 4$&
any fundamental representation\\
 $Spin_{p,q}/\left(Spin_p\times Spin_q\right), p,q\equiv 1\ mod\ 2, p\equiv q\ mod\ 4$&
positive and negative half-spinor representation\\
% $Spin_l\left({\Bbb C}\right)/Spin_l, l\equiv 3\ mod\ 4$&
%the spinor representation and its conjugate\\
 $SO_l\left({\Bbb C}\right)/SO_l, l\equiv 3\ mod\ 4$& any fundamental representation\\
 $Sp_l\left({\Bbb C}\right)/Sp_l, l\equiv 1\ mod\ 4$& any fundamental representation\\
 $E_7\left({\Bbb C}\right)/E_7$& any fundamental representation\\

\hline
\end{tabular}\\
\end{thm}
\begin{pf}
By \hyperref[cartan]{Lemma \ref*{cartan}}
it suffices to check whether $\beta_{2n-1}\left(\pi\right)\not\equiv 0$, where $\pi$ is the Lie algebra representation induced by $\rho$. Thus we can use the results from Section 3.2.

We use the classification of symmetric spaces as it can 
be read off Table 4 in \cite[p.229 ff.]{ov}. Of course, we are only 
interested in symmetric spaces of odd dimension. 
A simple inspection shows that all odd-dimensional irreducible symmetric spaces of noncompact type are given by the following list:\\
\\
\begin{tabular}[ht]{|l|c|}
  \hline
Symmetric Space & Dimension\\
  \hline\hline
  $SL_l\left({\Bbb R}\right)/SO_l, l\equiv 0,3,4,7\ mod\ 8$& $\frac{1}{2}\left(l-1\right)\left(l+2\right)$\\
$SL_l\left({\Bbb C}\right)/SU_l, l\equiv 0\ mod\ 2$ & $l^2-1$\\
$SL_{2l}\left({\Bbb H}\right)/Sp_l, l\equiv 0\ mod\ 2$ & $\left(l-1\right)\left(2l+1\right)$\\
$Spin_{p,q}/\left(Spin_p\times Spin_q\right), p,q\equiv 1\ mod\ 2$ & $pq$\\
$SO_l\left({\Bbb C}\right)/SO_l, l\equiv 2,3\ mod\ 4$& $\frac{1}{2}l\left(l-1\right)$\\
$Sp_l\left({\Bbb C}\right)/Sp_l, l\equiv 1\ mod\ 2$& $l\left(2l+1\right)$\\
$E_7\left({\Bbb C}\right)/E_7$& 163\\

  \hline
\end{tabular}\\

If $2n-1\equiv 3\mbox{\ mod\ }4$, then $n$ is even and by
\hyperref[even]{Corollary \ref*{even}} all representations $\rho:G\rightarrow GL\left(N,{\Bbb C}\right)$
satisfy $\rho^*b_{2n-1}\not=0$. This applies to  
the following symmetric spaces:\\
\\
\begin{tabular}[ht]{|l|c|}
  \hline
Symmetric Space & Condition\\
  \hline\hline

$SL_l\left({\Bbb R}\right)/SO_l$&$ l\equiv 0,7\ mod\ 8$\\
$SL_l\left({\Bbb C}\right)/SU_l$&$ l\equiv 0\ mod\ 2$\\
$SL_{2l}\left({\Bbb H}\right)/Sp_l$&$ l\equiv 0\ mod\ 4$\\
$Spin_{p,q}/\left(Spin_p\times Spin_q\right)$&$ p,q\equiv 1\ mod\ 2, p\not\equiv q\ mod\ 4$\\
$SO_l\left({\Bbb C}\right)/SO_l$&$ l\equiv 3\ mod\ 4$\\
$Sp_l\left({\Bbb C}\right)/Sp_l$&$ l\equiv 1\ mod\ 4$\\
$E_7\left({\Bbb C}\right)/E_7$&\\
 \hline
\end{tabular}\\
\\
Next we look at the irreducible
locally symmetric spaces of
dimension $\equiv 1\ mod\ 4$.\\
\\
For those $G/K$, whose Lie algebra $\g$ is {\em not} a complex Lie algebra
(this concerns the first 3 cases), we can, as observed in Section 
3.2.1, directly apply the results for the respective complexifications. Thus we have to check whether $\beta_{2n-1}\left(\rho_{\Bbb C}\right)\not=0$.\\
- For $SL_l\left({\Bbb R}\right)/SO_l, l\equiv 3,4\ mod\ 8$, every
fundamental representation $\rho$ satisfies $\rho^*b_{2n-1}\not=0$. (Indeed $l\ge 3,
n\ge 3$, thus we are not in one of the exceptional cases from Section 3.2.2.)\\
- For $SL_{2l}\left({\Bbb H}\right)/Sp_l, l\equiv 2\ mod\ 4$, every 
fundamental representation $\rho$ satisfies $\rho^*b_{2n-1}\not=0$. (Indeed the complexification of $sl_{2l}\left({\Bbb H}\right)$
is $sl_{4l}\left({\Bbb C}\right)$. We have $4l\ge 8$ and $n\ge 3$,
thus we are not in one of the exceptional cases from Section 3.2.2.)\\
- For $Spin_{p,q}/\left(Spin_p\times Spin_q\right), p,q\equiv 1\ mod\ 2, p\equiv q\ mod\ 4$, the
positive and negative half-spinor representations are the only fundamental representations $\rho$ satisfying $\rho^*b_{2n-1}\not=0$. (The assumptions imply that the complexification is $so\left(2l,{\Bbb C}\right)$ with $l$ odd, because of $2l=p+q\equiv 2\ mod\ 4$. In particular $n\equiv l\ mod\ 2$ and we are not in the exceptional case of Section 3.2.3.)\\
\\
For those $G/K$ whose Lie algebra $\g$ is a complex Lie algebra, we use the fact that each ${\Bbb R}$-linear representation is of the form $\rho_1\otimes\overline{\rho_2}$. We get:\\
- For $SO_l\left({\Bbb C}\right)/SO_l, l\equiv 3\ mod\ 4$, 
 we have $l\equiv n\ mod\ 2$ and
by Section 3.2.4 no fundamental representation $\rho$ satisfies $\rho^*b_{2n-1}\not=0$.\\
%the spinor representation and its conjugate 
%are the only fundamental representations $\rho$ satisfying $\rho^*b_{2n-1}\not=0$,\\
- For $Sp_l\left({\Bbb C}\right)/Sp_l, l\equiv 1\ mod\ 4$, 
by Section 3.2.5 no fundamental representation $\rho$ satisfies $\rho^*b_{2n-1}\not=0$.
\end{pf}\\
\\
{\bf Example (Goncharov)}: Consider hyperbolic space ${\Bbb H}^d=Spin_{d,1}/\left(Spin_d\times Spin_1\right).$
It was shown in \cite{gon} that the half-spinor representations have nontrivial Borel class if $d$ is odd. 
%The question was raised (\cite[page 587]{gon})
%whether these are the only fundamental representations of
%$Spin_{n,1}$ with this property. As a
%special case of 
The above results show that for $d\equiv 3\ mod\ 4$ each 
irreducible representation has nontrivial Borel class, but
for $d\equiv 1\ mod\ 4$ the positive and half-negative spinor representation
 are the only fundamental representations with this property.\\
For $d=3$ we will however compute in Section 3.4 that the invariants coming from irreducible representations, albeit all distinct, are rational multiples of each other.

		\subsection{Some clues on computation}

So far we have been using \hyperref[cartan]{Lemma \ref*{cartan}} to decide when $\rho^*b_{2n-1}\not=0$.
%, which is in view of \hyperref[cartan]{Lemma \ref*{cartan}} easier than computing $\rho^*b_{2n-1}$.
In this subsection we will give some clues to
the actual computation of $\rho^*b_{2n-1}$.
%Its results are not needed 
%for the remainder of the paper.

Recall that $H^*\left(\g_u,\kk\right)$ is the cohomology of the complex of $G_u$-invariant forms on $G_u/K$.
For a $d$-dimensional compact symmetric space $G_u/K$, there is an isomorphism $H^d\left(\g_u,\kk\right)\simeq H^d\left(G_u/K;{\Bbb R}\right)\simeq{\Bbb R}$ given by integration 
over $\left[G_u/K\right]$. Moreover, a $G_u$-invariant $d$-form is uniquely determined 
by its value on an orthonormal (for the metric given by the negative of the Killing form) basis $X_1,\ldots,X_d$ of $T_{\left[e\right]}G_u/K\simeq i\p$. By definition, the volume form takes the value 1 on each orthonormal basis. On the other hand, the volume form represents $vol\left(G_u/K\right)\left[G_u/K\right]^v\in H^d\left(G_u/K\right)$, where $\left[G_u/K\right]^v$ means the dual of the fundamental class. Thus we have:

		%For each Lie-algebra-cocycle $P\in C^n\left(\g_u,\kk\right)$, we denote
		%by $\omega_P\in \Omega^n\left(G_u/K\right)$ the corresponding $G_u$-invariant
		%differential form. Then we have the following obvious observation. ($\left[\omega_P
		%\right]$ denotes the cohomology class of $\omega_P$, and
		%$\left[G_u/K\right]^{v}\in H^n\left(G_u/K,{\Bbb R}\right)$ denotes the dual of the fundamental
		%class $\left[G_u/K\right]$. The Riemannian metric is given by the Killing form $B$.)

		\begin{lem}\label{basis} Let $G_u/K$ be a compact symmetric space of dimension $d$, $\omega\in C^d\left(\g_u,\kk\right)$ a $G_u$-invariant $d$-form and $X_1,\ldots,X_d$ an orthonormal
basis of $i\p$.
%		Then, for each $P\in C^n\left(\g_u,\kk\right)$, we
%		have $
Then $$\left[\omega\right]=\omega\left(X_1,\ldots,X_d\right) vol\left(G_u/K\right)\left[G_u/K\right]^{v}\in H^d\left(G_u/K;{\Bbb R}\right).$$  \end{lem}

%		\begin{cor} $\left[\omega_P\right]\not=0$ iff $P\left(X_1,\ldots,X_n\right)\not=0$ 
%		for some (hence any) basis of $i\p$.\end{cor}

%We will apply this to 
The Borel class $b_{2n-1}\in H^{2n-1}
\left(u\left(N\right)
%\left(u\left(N\right)\oplus u\left(N\right),u\left(N\right)
\right)$ is represented by 
%the relative Lie-algebra-cocycle
$b_{2n-1}\left(X_1,\ldots,X_{2n-1}\right)=$
$$\frac{1}{\left(2\pi i\right)^n}\frac{\left(-1\right)^n\left(n-1\right)!}{\left(2n-1\right)!}\sum_{\sigma\in S_{2n-1}}\left(-1\right)^\sigma Tr\left(X_{\sigma\left(1\right)}\left[X_{\sigma\left(2\right)},X_{\sigma\left(3\right)}\right]\ldots
\left[X_{\sigma\left(2n-2\right)},X_{\sigma\left(2n-1\right)}\right]\right).$$ 

Under the identification $H^{2n-1}\left(u\left(N\right)\oplus u\left(N\right),u\left(N\right)\right)\simeq 
H^{2n-1}
\left(u\left(N\right)
%\left(u\left(N\right)\oplus u\left(N\right),u\left(N\right)
\right)$ this gives $b_{2n-1}^{u\oplus u}$ represented by 
$b_{2n-1}^{u\oplus u}\left(Y_1,\ldots,Y_{2n-1}\right)
=b_{2n-1}\left(X_1,\ldots,X_{2n-1}\right)$. Here $Y_1,\ldots,Y_{2n-1}\in u\left(N\right)\oplus u\left(N\right)$ and $X_1:=\pi_2\left(Y_1\right),\ldots, X_{2n-1}:=\pi_2\left(Y_n\right)\in u\left(N\right)$, where $\pi_2$ is the projection to the second summand of $u\left(N\right)\oplus u\left(N\right)$. 

After the identification $H^{2n-1}\left(gl\left(N,{\Bbb C}\right),
u\left(N\right)\right)\simeq 
H^{2n-1}\left(u\left(N\right)\oplus u\left(N\right),
u\left(N\right)\right)$ 
this gives $b_{2n-1}^{gl}$ represented 
by $b_{2n-1}^{gl}\left(y_1,\ldots,y_{2n-1}\right)=i^{2n-1}b_{2n-1}
\left(\frac{x_1}{i},\ldots,\frac{x_{2n-1}}{i}\right)$. Here $x_1=\pi\left(y_1\right),\ldots,x_{2n-1}=
\pi\left(y_{2n-1}\right)$, where $\pi:gl\left(N,{\Bbb C}\right)\rightarrow iu\left(N\right)$
is the projection associated to $gl\left(N,{\Bbb C}\right)=u\left(N\right)\oplus iu\left(N\right)$.

{\bf Borel element.} 
%One can use the canonical isomorphism $C^{2n-1}\left(u\left(N\right)\otimes u\left(N\right),u\left(N\right)\right)\cong
%C^{2n-1}\left(gl\left(N,{\Bbb C}\right),u\left(N\right)\right)$ to consider $b_{2n-1}$ as an element in $\left(gl\left(N,{\Bbb C}\right),u\left(N\right)\right)$.
In the notation of \cite[Section 9.7]{bu} is $b_{2n-1}^{gl}=\frac{1}{\left(2\pi i\right)^n}\Phi_{2n-1}$.
The {\em Borel element} $Bo_n\in C^*\left(gl\left(N,{\Bbb C}\right), u\left(N\right); {\Bbb R}\left(n-1\right)\right)$ is
given in \cite[Section 9.7]{bu} by $Bo_n\left(\wedge_{j=1}^{2n-1}y_j\right)=\Phi_{2n-1}\left(\wedge_{j=1}^{2n-1}
\left(\overline{y}_j^t+y_j\right)\right))$ and this is then used to define the Borel regulator.

For $x_j\in iu\left(N\right)$ we have $\overline{x}_j^t+x_j=2x_j$, hence 
%the formula simplifies to
$Bo_n\left(\wedge_{j=1}^{2n-1}x_j\right)
=\left(2\pi i\right)^n 2^{2n-1}b_{2n-1}\left(x_1,\ldots,x_{2n-1}\right)$. \\

%Under the canonical isomorphism $C^{2n-1}\left(u\left(N\right)\otimes u\left(N\right),u\left(N\right)\right)\cong
%C^{2n-1}\left(gl\left(N,{\Bbb C}\right),u\left(N\right)\right)$ tWe
%\noindent
{\bf Example: Hyperbolic 3-manifolds.}
%Let $G/K=SL\left(2,{\Bbb C}\right)/SU\left(2\right)$. 
%Then $$i\p=\left\{iA\in Mat\left(2,{\Bbb C}\right): Tr\left(A\right)=0, A=\overline{A}^T\right\}=\left\{B\in Mat\left(2,{\Bbb C}\right): Tr\left(B\right)=0, B=-\overline{B}^T\right\}.$$
		A Killing form-orthonormal 
%(with respect to the negative Killing form) 
basis of $T_{\left[e\right]}SL\left(2,{\Bbb C}\right)/SU\left(2\right)=\left\{B\in Mat\left(2,{\Bbb C}\right): Tr\left(B\right)=0, B=\overline{B}^T\right\}$ is 
%given by 
$\left\{\frac{1}{2\sqrt{2}}H,
	\frac{1}{2\sqrt{2}}X, \frac{1}{2\sqrt{2}}Y\right\}$, with
	$$H=\left(\begin{array}{cc}1&0\\0&-1\end{array}\right), X=\left(\begin{array}{cc}0&1\\1&0\end{array}\right),
		Y=\left(\begin{array}{cc}0&i\\-i&0\end{array}\right).$$
		%We have $$\left[H,X\right]=-2Y,\left[H,Y\right]=-2X,\left[X,Y\right]=2H.$$
		%Thus, for each representation
For a representation $\rho:Sl\left(2,{\Bbb C}\right)\rightarrow GL\left(m+1,{\Bbb C}\right)$ 
%with associated
%Lie algebra representation $\pi$
and $\pi=D_e\rho$
%:sl\left(2,{\Bbb C}\right)\rightarrow Mat\left(m+1,{\Bbb C}\right)$ 
we have
		$$\rho^*b_3^{gl}\left(H,X,Y\right)=\frac{i}{\left(2\pi i\right)^2}\frac{1}{6}\left\{2Tr\left(\pi iH\left[\pi iX,\pi iY
		\right]\right)+
2Tr\left(\pi iX\left[\pi iY,\pi iH\right]\right)
+2Tr\left(\pi iY\left[\pi iH,\pi iX
                \right]\right)\right\}$$
		$$=-\frac{1}{6\pi^2}Tr\left(\left(\pi iH\right)^2\right)-\frac{1}{6\pi^2}Tr\left(\left(\pi iX\right)^2\right)
-\frac{1}{6\pi^2}Tr\left(\left(\pi iY\right)^2\right).$$
%		By the classification of irreducible representations of $sl\left(2,{\Bbb C}
%		\right)$, 
Each $m+1$-dimensional irreducible representation
		is equivalent to $\pi_m$ given
		by $$\pi_m\left(iH\right)=\left(\begin{array}{ccccc}im&0&0&...&0\\
		0&i\left(m-2\right)&0&...&0\\
%		0&0&i\left(m-4\right)&...&0\\
		.&.&.&...&.\\
		0&0&.&...&-im\end{array}\right), $$
$$		\pi_m\left(iX\right)=
		\left(\begin{array}{ccccc}0&-i&0&...&0\\
		-im&0&-2i&...&0\\
		0&-i\left(m-1\right)&0&...&0\\
		.&.&.&...&-im\\
		0&0&0&..-i&0\end{array}\right),
\pi_m\left(iY\right)=
		\left(\begin{array}{ccccc}0&1&0&...&0\\
		-m&0&2&...&0\\
		0&-\left(m-1\right)&0&...&0\\
		.&.&.&...&m\\
		0&0&0&..-1&0\end{array}\right).$$
		The diagonal
		entries of
		$\pi_m\left(iH\right)^2$ are $$\left(
		-m^2,-\left(m-2\right)^2,\ldots,0,\ldots,-\left(m-2\right)^2,-m^2\right),$$
and
the diagonal entries of
		$\pi_m\left(iX\right)^2$ and
		$\pi_m\left(iY\right)^2$ are both equal to
$$\left(-m,-m-2\left(m-1\right),-2\left(m-1\right)-3\left(m-2\right), ...\right).$$
In particular
		$Tr\left(\pi_m\left(iX\right)^2\right)=Tr\left(\pi_m\left(iY\right)^2\right)$.
After scaling with $\frac{1}{2\sqrt{2}}$ we conclude
		$$\rho_m^*b_3^{gl}\left(\frac{1}{2\sqrt{2}}H,\frac{1}{2\sqrt{2}}X,
	\frac{1}{2\sqrt{2}}Y\right)=
%-\frac{1}{96\sqrt{2}\pi^2}Tr\left(\left(\pi_m H\right)^2\right)-\frac{1}{48\sqrt{2}\pi^2}Tr\left(\left(\pi_m X\right)^2\right)$$
%$$=
\frac{i}{96\sqrt{2}\pi^2}\sum_{k=0}^m\left(m-2k\right)^2+\frac{i}{48\sqrt{2}\pi^2}
\sum_{k=0}^m
k\left(m-k+1\right)+\left(k+1\right)\left(m-k\right).$$
%Since the Borel regulator corresponds to $2^{2n-1}\left(2\pi i\right)^n\rho^*b_{2n-1}=-32\pi^2\rho_1^*b_3$ 

{\bf Comparison of Borel element and hyperbolic volume form.} 

If $m=1$, that is for the inclusion $\rho_1=j:SL\left(2,{\Bbb C}\right)\subset GL\left(2,{\Bbb C}\right)$,
we get $$j^*b_3^{gl}\left(\frac{1}{2\sqrt{2}}H,\frac{1}{2\sqrt{2}}X,
        \frac{1}{2\sqrt{2}}Y\right)=\frac{i}{16\sqrt{2}\pi^2}, j^*b_3^{u\oplus u}\left(\frac{i}{2\sqrt{2}}H,\frac{i}{2\sqrt{2}}X,
        \frac{i}{2\sqrt{2}}Y\right)=\frac{1}{16\sqrt{2}\pi^2}.$$ Thus by \hyperref[basis]{Lemma \ref*{basis}} 
$j^*b_3^{u\oplus u}\in H^3\left(\g_u,\kk\right)\simeq H^3\left(S^3\right)$ represents $\frac{1}{16\sqrt{2}\pi^2} vol\left(S^3\right)\left[S^3\right]$. 

%One the other hand, since 
Explicit computation shows that the projection $SL\left(2,{\Bbb C}\right)\rightarrow {\Bbb H}^3$ maps $\frac{1}{2\sqrt{2}}H,\frac{1}{2\sqrt{2}}X,\frac{1}{2\sqrt{2}}Y$ to vectors of hyperbolic length $\frac{1}{\sqrt{2}}$. Thus
{\bf the hyperbolic metric is given by one half of the Killing form}: lenghts are multiplied by $\frac{1}{\sqrt{2}}$, volumes by $\frac{1}{2\sqrt{2}}$. 
By \hyperref[basis]{Lemma \ref*{basis}} this means that
the isomorphism $H^3\left(sl\left(2,{\Bbb C}\right),su\left(2\right)\right)\simeq H^3\left(S^3\right)$
%frac{1}{2\sqrt{2}}$.
%we have that evaluation of the hyperbolic volume form on the orthonormal basis gives $\frac{1}{\sqrt{2}^3}$, thus its cohomology class represents 
sends the class $\left[dvol\right]$ of the hyperbolic volume form to $\frac{1}{2\sqrt{2}}vol\left(S^3\right)\left[S^3\right]$. In particular $j^*b_3^{gl}\in C^3\left(sl\left(2,{\Bbb C}\right),su\left(2\right)\right)$ is $\frac{1}{8\pi^2}$ times the class
of the hyperbolic volume form.

For the Borel element we have $Bo_2=-32\pi^2b_3^{gl}$, it follows that the                                                                        Borel element is $-16$ times the hyperbolic volume.
%Thus an  orthonormal basis of $\p=
%T_{\left[e\right]}{\Bbb H}^3$ with respect to the hyperbolic metric is given by $\left\{-\frac{i}{2}H,-\frac{i}{2}X,-\frac{i}{2}Y\right\}$. By definition, evaluation of the volume form on an orthonormal basis yields 1. Since evaluation of $\rho_1^*b_3$ on this orthonormal basis yields $\frac{1}{8\pi^2}$, we conclude that
%the Borel class gives $\frac{1}{8\pi^2}$ times the hyperbolic volume form, when both are considered as elements of $C^3\left(sl\left(2,{\Bbb C}\right),su\left(2\right)\right)$. 
%{\bf Borel element.} The Borel element $Bo_n\in C^*\left(gl\left(N,{\Bbb C}\right),u\left(N\right);{\Bbb R}\left(n-1\right)\right)$
%is defined in \cite[Section 9.7]{bu} by $$Bo_n\left(x_1,\ldots,x_{2n-1}\right)
%=\left(2\pi i\right)^n b_{2n-1}\left(\overline{x}_1^t+x_1,\ldots,\overline{x}_{2n-1}^t+x_{2n-1}\right).$$
%For $x_j\in i\p$ we have $\overline{x}_j^t+x_j=2x_j$, hence the formula simplifies to
%$Bo_n\left(x_1,\ldots,x_{2n-1}\right)
%=\left(2\pi i\right)^n 2^{2n-1}b_{2n-1}\left(x_1,\ldots,x_{2n-1}\right)$. In the case at hand this means 
%Moreover the Borel regulator corresponds to $2^{2n-1}\left(2\pi i\right)^n\rho^*b_{2n-1}=
%$Bo_2=-32\pi^2b_3$. It follows that the
%Borel element is $-16$ times the hyperbolic volume.
(\cite{ds} and \cite{ny} compute the imaginary part of the Borel regulator to be $\frac{1}{2\pi^2}$ times the hyperbolic volume, but they are using a different definition.)\\
%while \cite{ny} defines the imaginary part of the Borel regulator 
%to be $\frac{1}{2\pi^2}$ times the hyperbolic volume.)\\

{\bf Example: $SL\left(3,{\Bbb R}\right)/SO\left(3\right)$.}  
Let $\rho:SL\left(3,{\Bbb R}\right)\rightarrow GL\left(3,{\Bbb C}\right)$
be the inclusion. Since $SL\left(3,{\Bbb R}\right)/SO\left(3\right)$
is 5-dimensional, we wish to compute $\rho^*b_5$. 
Let $$H_1=\left(\begin{array}{ccc}i&0&0\\0&-i&0\\0&0&0\end{array}\right),
X_1=\left(\begin{array}{ccc}0&i&0\\i&0&0\\0&0&0\end{array}\right),
Y_1=\left(\begin{array}{ccc}0&-1&0\\1&0&0\\0&0&0\end{array}\right).$$
We will use the convention that, for $A\in\left\{H,X,Y\right\}$ if $A_1$ is defined (in a given basis), then $A_2$ is obtained via the base change $e_1\rightarrow e_2, e_2\rightarrow e_3, e_3\rightarrow e_1$ and $A_3$ is obtained
via the base change $e_1\rightarrow e_3,
e_3\rightarrow e_2, e_2\rightarrow e_1$.

We have $\left[H_1,H_2\right]=0, \left[H_1,X_1\right]=2Y_1, \left[H_1,X_2\right]=-Y_2,\left[H_1,X_3\right]=-Y_3,\left[X_1,X_2\right]=iY_3$ and more relations are obtained out of these ones by base changes.

A basis of $i\p$ is given by $H_1,H_2,X_1,X_2,X_3$. The formula for $\rho^*b_5
\left(H_1,H_2,X_1,X_2,X_3\right)$ contains 120 summands. (24 of them contain $\left[H_1,H_2\right]=0$
or $\left[H_2,H_1\right]=0$.)

Each summand appears four times because, for example, $H_1\left[H_2,X_1\right]\left[X_2,X_3\right]$ also shows up as $
-H_1\left[X_1,H_2\right]\left[X_2,X_3\right], 
-H_1\left[H_2,X_1\right]\left[X_3,X_2\right]$ and $
H_1\left[X_1,H_2\right]\left[X_3,X_2\right]$. Thus one has to add 30 summands (6 of them zero), and multiply their sum by $4$.

We note that all summands of the form $H_1\left[H_2,.\right]\left[.,.\right]$
give after base change corresponding elements of the form $H_2\left[H_1,.\right]\left[.,.\right]$, which are summed with the opposite sign. Thus these terms cancel each other. The same cancellation occurs between summands of the form
$X_2\left[.,.\right]\left[.,.\right]$ and $X_3\left[.,.\right]\left[.,.\right]$. Thus we only have to sum up summands of the form $X_1\left[.,.\right]\left[.,.\right]$ and we get 
$$\left(2\pi i\right)^35!\rho^*b_5\left(H_1,H_2,X_1,X_2,X_3\right)=$$
$$4 Tr\left(X_1\left[H_1,H_2\right]\left[X_2,X_3\right]\right)
+4 Tr\left(X_1
\left[X_2,X_3\right]\left[H_1,H_2\right]\right)+$$
$$+4 Tr
\left(X_1
\left[H_1,X_2\right]\left[X_3,H_2\right]\right)
+
4 Tr\left(X_1
\left[X_3,H_2\right]\left[H_1,X_2\right]\right)+$$
$$+
4 Tr \left(X_1
\left[H_1,X_3\right]\left[H_2,X_2\right]\right)+
4 Tr \left(X_1
\left[H_2,X_2\right]\left[H_1,X_3\right]\right)$$
$$=0+0+4 Tr\left(X_1Y_2Y_3\right)+ 4 Tr\left(X_1Y_3Y_2\right)
+4 Tr\left(-2X_1Y_3Y_2\right) +4 Tr\left(-2 X_1Y_2Y_3\right)$$
$$=
0+0+4i+4i-8i-8i=-8i.$$
We note that $H_1,H_2,X_1,X_2,X_3$ are pairwise orthogonal and have norm $2\sqrt{3}$. Dividing each of them by $2\sqrt{3}$ gives an 
orthonormal
basis, on which evaluation of $\rho^*b_5$ gives 
$$\rho^*b_5\left(\frac{1}{2\sqrt{3}}H_1,\frac{1}{2\sqrt{3}}H_2,\frac{1}{2\sqrt{3}}X_1,
\frac{1}{2\sqrt{3}}X_2,\frac{1}{2\sqrt{3}}X_3\right)=
\frac{1}{\left(2\sqrt{3}\right)^5}\frac{1}{5!}
\frac{1}{\left(2\pi i\right)^3}\left(-8i\right)=\frac{1}{34560\sqrt{3}\pi^3}.$$
The Borel element is 
%$2^{2n-1}\left(2\pi i\right)^nb_{2n-1}=
$-256\pi^3ib_5$, thus its value on the orthonormal basis is $-\frac{i}{135\sqrt{3}}$.

\section{The cusped case}

{\bf Outline.} In Section 2 we defined $\gamma\left(M\right)$ for closed manifolds $M=\Gamma\backslash G/K$, using the image of the fundamental class $\left[M\right]\in H_*
\left(M\right)\cong H_*\left(B\Gamma\right)$
in $H_*\left(BG\right)$ for the 
construction, and the volume 
cocycle in $C_{simp}^*\left(BG\right)$ for the proof of the desired nontriviality
properties. In this section we would like to give an analogous construction for cusped manifolds.

For this we would like
to map the fundamental class $$\left[M,\partial M\right]\in H_*\left(DCone\left(\cup_i\partial_i M\rightarrow M\right)\right)\cong
H_*\left(DCone\left(\cup_i B\Gamma_i \rightarrow B\Gamma\right)\right)$$ (all
notions are defined in Section 4.2) to the homology of some completion $BG^{comp}$ of $BG$, where the completion should be chosen such that the
volume class extends to $BG^{comp}$.

The completion $BG^{comp}$ that we define in Section 4.2.2 will be chosen such that for each $c\in\partial_\infty G/K$ a cone over $BG$ is added. There is a natural extension of the volume class to this set and the addition of all points at infinity will leave us the flexibility to remember the geometry of cusps.

Now, if $G/K$ has rank one, then each path-component $\partial_i M$ corresponds to a cusp $c_i\in\partial_\infty G/K$ and this will allow us to define an image of 
$\left[M,\partial M\right]$ in $H_*\left(BG^{comp}\right)$. 

Of course this does not apply to $SL\left(N,{\Bbb C}\right)/SU\left(N\right)$, which has rank $N-1$, but if $\rho:\left(G,K\right)
\rightarrow \left(SL\left(N,{\Bbb C}\right),SU\left(N\right)\right)$ is a representation for a rank one space $G/K$, then we get a well-defined image of $c_i$ in $\partial_\infty
SL\left(N,{\Bbb C}\right)/SU\left(N\right)$ and can thus define the image of $\left[M,\partial M\right]$ in $H_*\left(BSL\left(N,{\Bbb C}\right)^{comp}\right)$.

In Section 4.4 we will show that this homology class has a preimage in $H_*\left(BSL\left(N,{\Bbb C}\right)\right)$. This then finally allows to generalize the Goncharov 
construction (Theorem 4).

\subsection{Preparations}
Let $G$ be a connected, semisimple Lie group with maximal compact subgroup $K$.
Thus $G/K$ is a symmetric space of noncompact type. Throughout Section 4 we will make
the assumption $rank\left(G/K\right)=1$.
%$G$ acts on $\widetilde{M}=G/K$ and thus on the ideal boundary $\partial_\infty \widetilde{M}$. 
%Let $G=KAN$ be an Iwasawa decomposition of $G$. We say that a group $P\subset G$ is parabolic if there is some $x\in
%\partial_\infty \widetilde{M}$ with $P\subset\left\{g\in G: gx=x\right\}$. This is equivalent to the condition that there exists a $g\in G$ with $gPg^{-1}\subset AN$.

We will consider a manifold $M$ with boundary $\partial M$ such that $int\left(M\right)=M-\partial M$ is a finite-volume locally symmetric space of noncompact type
of rank one. This means 
%that there is a symmetric space $G/K$ of noncompact type
%of rank one and a discrete subgroup $\Gamma\subset G$ such that 
$$int\left(M\right)=\Gamma\backslash G/K$$ for a (not necessarily cocompact) lattice $\Gamma\subset G$.
%In this section we will assume that $\Gamma\backslash G/K$ has finite volume but need not be compact.

%It maps $K$ into some maximally compact subgroup $K^\prime\subset GL\left(N,{\bf C}\right)$.
We note that connected,
semisimple Lie groups are perfect, hence each representation $\rho:G\rightarrow GL\left(N,{\Bbb C}\right)$
has image in
$SL\left(N,{\Bbb C}\right)$. Further we will assume that $\rho$ maps $K$ to $SU\left(N\right)$, which can be achieved upon conjugation.
% and maps $K$ to $SU\left(N\right)$. 
%We denote $\widetilde{X}$ the symmetric space $\widetilde{X}=SL\left(N,{\bf C}\right)/K^\prime$.
%The induced map $$\rho:G/K\rightarrow SL\left(N,{\Bbb C}\right)/SU\left(N\right)$$ 
%maps geodesics in $G/K$ 
%to geodesics in $SL\left(N,{\Bbb C}\right)/SU\left(N\right)$ and does not increase distances.
%Thus one gets a map of ideal boundaries $$\partial_\infty \widetilde{M}\rightarrow
%\partial_\infty \left(SL\left(N,{\Bbb C}\right)/
%SU\left(N\right)\right)$$ which is $\rho$-equivariant for the actions 
%of $G$ resp.\ $SL\left(N,{\Bbb C}\right)$. In particular, parabolic subgroups are mapped to parabolic subgroups.
%We recall that
%any maximal unipotent subgroup $B\subset SL\left(N,{\Bbb C}\right)$
%is conjugate to the group $B_0$ of upper triangular matrices with all diagonal entries equal 1. 
%If $\Gamma\subset G$ is a discrete subgroup of finite covolume, then we can assume $\Gamma\subset G\left(\overline{\Bbb Q}\right)$ by Weil rigidity. Note that, for a semisimple Lie group $G$, each representation $\rho:G\rightarrow SL\left(N,{\Bbb C}\right)$ maps $G\left(\overline{\Bbb Q}\right)$ to $SL\left(N,\overline{\Bbb Q}\right)$, and for each subring $\ff\subset{\Bbb C}$, maps $G\left(\ff\right)$ to $SL\left(\ff\right)$. 
\subsubsection{ Negative curvature and visibility manifolds}

If $int\left(M\right)=\Gamma\backslash G/K$ is a locally symmetric space of noncompact type of rank one, then its sectional curvature $sec$
is bounded between two negative constants, after scaling with a constant factor one has
$$-4\le sec \le -1.$$ 
In particular, by \cite[page 440]{ebe}, the universal covering $\widetilde{int\left(M\right)}=G/K$ 
is a 'visibility manifold' in the sense of \cite{ebe}. 

The structure of finite-volume quotients of visibility manifolds has been described in \cite{ebe}. The following Lemma collects those results from the proof of \cite[Theorem 3.1]{ebe} that we will frequently use in this paper.
(We denote by $\partial_\infty\widetilde{int\left(M\right)}=\partial_\infty\left(G/K\right)$ the ideal boundary 
of $\widetilde{int\left(M\right)}=G/K$, that is the set of equivalence classes of geodesic rays, where rays are equivalent if they
are asymptotic, see
%. The union of $\widetilde{int\left(M\right)}$ with its ideal boundary carries a well-known topology defined
\cite[Section 1]{ebe}.)
\begin{lem}\label{visibility}
Let $\widetilde{N}$ be a simply connected, complete Riemannian manifold, $\Gamma$ be a discrete group of isometries of $\widetilde{N}$ and $N=\widetilde{N}/\Gamma$.

If 
$\widetilde{N}$ is a visibility manifold (\cite{ebe}) of nonpositive sectional curvature and $N$ has finite volume, then each end of
$N$
% $int\left(M\right)=M-\partial M$ 
has a neighborhood $E$
homeomorphic to $U_c/P_c$, where $c\in\partial_\infty\widetilde{N}$, $U_c$ is a horoball centered at $c$ and $P_c\subset \Gamma$ is a discrete group of parabolic isometries fixing $c$.

In particular, if $N$ has finitely many ends, then there are end neighborhoods $E_1,\ldots,E_s$ such that $K=N-\cup_{i=1}^s E_i$ 
is compact and for $i=1,\ldots,s$ there are homeomorphisms of pairs $\left(E_i,\partial \overline{E}_i\right)\rightarrow \left(U_{c_i}/P{c_i},L_{c_i}/P_{c_i}\right)$, where $c_i\in\partial_\infty \widetilde{N}$ and $L_{c_i}$ is 
the horosphere centered at $c_i$ which bounds the horoball $U_{c_i}$.\end{lem}
%\begin{pf} This is shown in the proof of \cite[Theorem 3.1]{ebe}.\end{pf}
\begin{cor}\label{homeo}
If $M$ is a compact manifold with boundary, $\partial_1M,\ldots,\partial_s M$ are the connected components of $\partial M$, and $N:=int\left(M\right)=M-\partial M$ carries a Riemannian metric of finite volume such that $\widetilde{N}$ is a visibility manifold, then, 
with the notation of \hyperref[visibility]{Lemma \ref*{visibility}}, we have a homeomorphism of tuples
$$\left(M,\partial_1 M,\ldots,\partial_sM\right)\rightarrow \left(\left(\widetilde{N}-\cup_{i=1}^s U_{c_i}\right)/\Gamma,L_{c_1}/P_{c_1},\ldots, L_{c_s}/P_{c_s}\right).$$\end{cor} 
\begin{pf} By the proof of \cite[Theorem 3.1]{ebe}, the neighborhood $E_i$ is {\em Riemannian collared}, which implies in particular the existence of a diffeomorphism $E_i\cong\partial \overline{E}_i\times\left(0,\infty\right)$.
The claim follows.\end{pf}\\
\\
%and if $E_i\simeq \partial_iM\times\left(0,\infty\right)$ a neighborhood of 
%the end of $int\left(M\right)=M-\partial M$ that 
%corresponds to $\partial_iM$  has a neighborhood $E_i$ \end{lem}
%If $\partial M$ has connected components $\partial_1M,\ldots,\partial _sM$ and if the end of $int\left(M\right)$ corresponding to $\partial_iM$ has a neighborhood $E_i=U_{c_i}/P_{c_i}$ with $c_i\in\partial_\infty\widetilde{int
%\left(M\right)}$, then 
We will
say that $\Gamma c_i\subset \partial_\infty\widetilde{N}$ is {\em the set of parabolic fixed points corresponding to $\partial_iM$}.

It is at this point where we need the assumption $rank\left(G/K\right)=1$. In the higher 
rank case it is not true that there is a unique $\Gamma$-orbit of 
parabolic fixed points $\Gamma c_i\subset \partial_\infty \left(G/K\right)$ 
associated to a boundary component $\partial_iM$. The isomorphism $\pi_1M\cong \Gamma$ does {\em not} send $\pi_1\partial_iM$ to a subgroup of some $Fix\left(c_i\right)$, if $rank\left(G/K\right)\ge 2$.\\
%, where $\Gamma$ is $\pi_1M=\pi_1int\left(M\right)$ acting by deck transformations on $\widetilde{int\left(M\right)}$ and thus 
%on $\partial_\infty\widetilde{int\left(M\right)}$.\\
\\
{\bf $\pi_1$-injective boundary.}
In the proof of \hyperref[preimage]{Proposition \ref*{preimage}} and \hyperref[Thm3]{Theorem \ref*{Thm3}} we will use that $\pi_1\partial_iM\rightarrow\pi_1M$ is
injective for each path-component $\partial_iM$ of $\partial M$. We are going to explain how this fact follows from well-known properties of visibility manifolds.
\begin{cor}\label{piinjective} 
Under the assumptions of \hyperref[homeo]{Corollary \ref*{homeo}} we have
%Let $M$ be a compact manifold with boundary $\partial M$.
%If $int\left(M\right)$ admits a Riemannian metric of finite volume such that the universal covering $\widetilde{M}$ with the pull-back metric is a visibility manifold of nonpositive sectional curvature,
that $\pi_1\partial_iM\rightarrow\pi_1M$ is
injective for each path-component $\partial_iM$ of $\partial M$.\end{cor}
\begin{pf}
%By covering theory, the universal covering $\widetilde{M}$ is the set of homotopy classes (rel.\ boundary) of paths $\gamma:\left[0,1\right]\rightarrow M$ with $\gamma\left(0\right)=x_0$, for a fixed point $x_0\in M$. If $\pi_1\partial_iM\rightarrow \pi_1M$ were not injective, then $\widetilde{\partial_iM}\rightarrow \widetilde{M}$ would not be injective. By covering theory, $\widetilde{M}$ is the, up to homeomorphism unique, simply connected covering of $M$.
From \hyperref[homeo]{Corollary \ref*{homeo}} we get a commutative diagram
\[\begin{xy}
\xymatrix{\partial_iM \ar[r] \ar[d] & M \ar[d]\\
L_{c_i}/P_{c_i}\ar[r] & \left(\widetilde{N}-\cup_{i=1}^s U_{c_i}\right)/\Gamma}
\end{xy}\]
where the vertical arrows are homeomorphisms, thus inducing isomorphisms $\pi_1\partial_iM\rightarrow P_{c_i}$ and $\pi_1M\rightarrow \Gamma$, and the horizontal arrows are induced by inclusions. 
If $P_{c_i}\rightarrow 
\Gamma$ were not injective, then the lift of $\iota:L_{c_i}/P_{c_i}\rightarrow \left(\widetilde{N}-\cup_{i=1}^s U_{c_i}\right)/\Gamma$ 
to the universal coverings would not be injective. However the lift of $\iota$ is the inclusion $\tilde{\iota}:L_{c_i}\rightarrow \widetilde{N}-\cup_{i=1}^s U_{c_i}$.
%$M=\widetilde{M}/\Gamma$ for a group of isometries $\Gamma\subset Isom\left(\widetilde{M}\right)$ with $\Gamma\simeq\pi_1M$.
%For $i=1,\ldots,s$ let $E_i$ be the end of $int\left(M\right)$ that corresponds
%to $\partial_iM$.
%By the proof of \cite[Theorem 3.1]{ebe} we have that each end $E_i$ of $int\left(M\right)$ has a neighborhood
%of the form $N_i\times\left(0,\infty\right)$ such that $N_i\times\left\{0\right\}=L_{c_i}/P_{c_i}$ 
%for a simply connected horosphere $L_{c_i}$ (centered at some $c_i\in\partial_\infty\widetilde{int\left(M\right)}$)
%and a subgroup $P_{c_i}\subset\Gamma$ that acts freely and properly discontinuously on $L_{c_i}$. The horosphere $L_{c_i}$ bounds the horoball $U_{c_i}$ given by \hyperref[visibility]{Lemma \ref*{visibility}}.
%This implies that there is a homeomorphism of tuples $$\left(M,\partial_1 M,\ldots,\partial_s M\right)\rightarrow \left(\left(\widetilde{int\left(M\right)}-\cup_{i=1}^s \Gamma U_{c_i}\right)/\Gamma, L_{c_1}/P_{c_1},\ldots, L_{c_s}/P_{c_s}\right)$$
%homotopy equivalence of pairs 
%In particular, the boundary component $\partial_iM$ that corresponds to this end is homeomorphic to $N$ and its fundamental group is isomorphic to $P_c$. We have $\widetilde{N}=L_c\subset \widetilde{M}$, which implies that $\pi_1N\rightarrow \pi_1M$ is injective because, by covering theory, $\widetilde{N}$ resp.\ $\widetilde{M}$ are
%the set of homotopy classes (rel.\ $\left\{0,1\right\}$) of paths $\gamma:\left[0,1\right]\rightarrow N$ resp.\ $M$ with $\gamma\left(0\right)=x_0$, for a fixed point $x_0\in N\subset M$.
\end{pf}\\
\\
Moreover $M$ and all $\partial_iM$ are aspherical by the Cartan-Hadamard Theorem and by \cite{ebe}.\\
\\
{\bf Identification of $\pi_1\partial_iM$ with a subgroup of $\pi_1M$.}
If $\partial M$ is not connected, then we have to choose different basepoints $x,x_1,\ldots,x_s$ for the definition of
$\pi_1\left(M,x\right),\pi_1\left(\partial_1M,x_1\right),\ldots,
\pi_1\left(\partial_sM,x_s\right)$. We can obtain subgroups $\Gamma_1,\ldots,\Gamma_s\subset\pi_1\left(M,x\right)$ isomorphic to $\pi_1\left(\partial_1M,x_1\right),\ldots,
\pi_1\left(\partial_sM,x_s\right)$, respectively, as follows:
\begin{df}\label{subgroup} Let $M$ be a manifold, $\partial_1M,\ldots,\partial_sM $ the connected components of $\partial
M$, $x\in M,x_1\in\partial_1M,\ldots, x_s\in\partial_sM, \Gamma=\pi_1\left(M,x\right)$.

Fix lifts $\tilde{x},\tilde{x}_1,\ldots,\tilde{x}_s$ of $x,x_1,\ldots,x_s$ to the universal covering $\pi:\widetilde{M}
\rightarrow M$, for $i=1,\ldots,s$ fix pathes 
$\tilde{l}_i:\left[0,1\right]\rightarrow \widetilde{M}$ with 
$\tilde{l}_i\left(0\right)=\tilde{x}$ and $\tilde{l}_i\left(1\right)=\tilde{x}_i$, let $l_i=\pi\circ \tilde{l}_i:\left[0,1\right]\rightarrow M$,
denote $\left[l_i\right]$ its homotopy class rel.\ $\left\{0,1\right\}$ and define
$$\Gamma_i:=\left\{\left[l_i\right]^{-1}*\gamma*\left[l_i\right]:\gamma\in\pi_1\left(\partial_iM,x_i\right)\right\}\subset\Gamma$$
to be the subgroup of $\Gamma$ which corresponds to $\pi_1\left(\partial_iM,x_i\right)$ after conjugation with $\left[l_i\right]$.\end{df}
The subgroup $\Gamma_i$ depends on the chosen lift 
$\tilde{x}_i$ but, for given $\tilde{x},\tilde{x}_i$, not on 
%the chosen path 
$\tilde{l}_i$.

With the homeomorphism from \hyperref[homeo]{Corollary \ref*{homeo}} we obtain $\Gamma_i=P_{c_i}$. We will say that $c_i$ is the {\em cusp associated to $\Gamma_i$}. In particular, $\Gamma_i\subset Fix\left(c_i\right)$.

(The choice of $c_i$ in its $\Gamma$-orbit depends on the chosen lift $\tilde{x}_i$ of $x_i$.)\\
\\
{\bf Compactification of universal covering by cusps.}
In the following Corollary we consider $\widetilde{int\left(M\right)}\bigcup\cup_{i=1}^s
\Gamma c_i$ as a subspace of $\widetilde{int\left(M\right)}\bigcup\partial_\infty \widetilde{int\left(M\right)}$, where the latter has the well-known
topology defined for example in Section 1 of \cite{ebe}.
The definition of the disjoint cone $DCone$ is given in Section 4.2.1 below.
\begin{cor}\label{cusps} Let the assumptions of \hyperref[visibility]{Lemma \ref*{visibility}} hold and let a fixed homeomorphism
$f:int\left(M\right)-\cup_{i=1}^s E_i\rightarrow M$ be given.
Then we have a projection $$\overline{\pi}:
\widetilde{int\left(M\right)}\bigcup\cup_{i=1}^s
\Gamma c_i\rightarrow DCone\left(\cup_{i=1}^s\partial_i M\rightarrow M\right) $$
such that $\overline{\pi}
\mid_{\widetilde{int\left(M\right)}- \cup_{i=1}^s \Gamma U_{c_i}}:\widetilde{int\left(M\right)}- \cup_{i=1}^s
\Gamma U_{c_i}\rightarrow
int\left(M\right)-\cup_{i=1}^s E_i$ is the restriction of the universal covering $\pi:\widetilde{int\left(M\right)}\rightarrow int\left(M\right)$,
$\overline{\pi}\mid_{\Gamma U_{c_i}}:
\Gamma U_{c_i}\rightarrow
 E_i \cup Cone\left(\partial_i M\right)-C_i$ is a covering with deck group $\Gamma$ and $\overline{\pi}$ maps
$\Gamma c_i$ to $C_i$ for $i=1,\ldots,s$, where $C_i$ is the cone point of $Cone\left(\partial_iM\right)$.\end{cor}
\begin{pf} Each boundary component $\partial_i M$ corresponds to an end (with neighborhood $E_i$) of $int\left(M\right)$ and thus by 
\hyperref[visibility]{Lemma \ref*{visibility}} to a unique $\Gamma$-orbit $\Gamma c_i$ with $c_i\in\partial_\infty\widetilde{int\left(M\right)}$ such that $E_i=U_{c_i}/P_{c_i}$. Let $\overline{E}_i$ be
the one-point compactification of $\overline{E}_i$, denote $C_i^+$ be the compactifying point, and let $M_+$ be the compactification of $M$ obtained by adding $C_1^+,\ldots,C_s^+$ to $M$. (This is homeomorphic to the space $M_+$ which will be considered in Section 4.1.2.)
Then we have 
homeomorphisms $f_0:int\left(M\right)-\cup_{i=1}^s E_i\rightarrow M$ and $f_i:\overline{E}_i\rightarrow Cone\left(\partial_iM\right)$ such that $f_0=f_i$ on $\partial E_i$ for $i=1,\ldots,s$, hence they yield a well-defined homeomorphism $f:M_+\rightarrow 
DCone\left(\cup_{i=1}^s\partial_i M\rightarrow M\right)$ which sends $C_i^+$ to $C_i$, the cone point over $\partial_iM$. 

Moreover, the universal covering $\pi:\widetilde{int\left(M\right)}\rightarrow int\left(M\right)$
sends $\gamma U_{c_i}$ to $E_i$ for each $\gamma\in\Gamma$, thus it can be continuously extended to $\Gamma c_i$ by $\pi\left(\gamma c_i\right)=C_i^+$ for $\gamma\in\Gamma$. 

Composition of $\pi$ with the homeomorphism $f$ yields the desired projection $\overline{\pi}$.\end{pf}\\
%is homeomorphic to $Cone\left(\partial_iM\right)$, by a homeomorphism sending the added point to the cone point $C_i$. Since $int\left(M\right) - \cup_{i=1}^s E_i$ is homeomorphic to $M$, this implies the existence of the projection.\end{pf}\\

Again, by the remark after \hyperref[homeo]{Corollary \ref*{homeo}}, also \hyperref[cusps]{Corollary \ref*{cusps}} requires the assumptions 
of \hyperref[visibility]{Lemma \ref*{visibility}} and would not work if $\widetilde{int\left(M\right)}=G/K$ were a symmetric space with $rank\left(G/K\right)\ge 2$.
(However there is a version of \hyperref[cusps]{Corollary \ref*{cusps}} for locally symmetric spaces of ${\Bbb Q}$-rank 1, which we will exploit in forthcoming work with Inkang Kim.)

\subsubsection{Generalized Cisneros-Molina-Jones construction}

The aim of Section 4 will be to associate a K-theoretic invarant to cusped locally symmetric spaces. 
We mention that, by an argument 
completely analogous to \cite{cj},
one can define an element $\alpha\left(M\right)$ and can associate to each representation  
$\rho:\left(G,K\right)\rightarrow \left(SL\left(N,{\Bbb C}\right),SU\left(N\right)\right)$
the push-forward
$$\left(B\rho\right)_d\left(\alpha\left(M\right)\right)\in H_d\left(B\left(SL\left(N,{\Bbb C}\right),{\mathcal{F}}\left(B
\right)\right)\right),$$
where $B\subset SL\left(N,{\Bbb C}\right)$ is a maximal unipotent subgroup.

In the case of hyperbolic 3-manifolds, Cisneros-Molina and Jones lifted
the invariant $\alpha\left(M\right)$ to $K_3\left({\Bbb C}\right)\otimes{\Bbb Q}$, and proved its nontriviality by relating it to the Bloch invariant. We describe 
now how to do a very similar construction for arbitrary locally 
symmetric spaces of noncompact type with finite volume. 
Unfortunately we did not succeed to evaluate the Borel class on the constructed invariant. This is 
the reason why we will actually pursue another approach, using relative group homology and closer in spirit to \cite{gon}, in the remainder of this section. The construction is however included at this point because its main step, \hyperref[mapR]{Lemma \ref*{mapR}}, will be crucial for the proof of \hyperref[preimage]{Proposition \ref*{preimage}}.

Let $M$ be an aspherical (compact, orientable, connected) $d$-manifold with aspherical boundary, $\ff\subset{\Bbb C}$ a subring
and $\rho: \pi_1M\rightarrow SL\left(\ff\right)$ a representation\footnotemark\footnotetext[3]{Notation: We will denote by ${\Bbb F}\subset{\Bbb C}$ an arbitrary subring (with 1), while $A\subset{\Bbb C}$ will denote a subring satisfying the assumptions of \hyperref[proj]{Lemma \ref*{proj}}.}.

To push forward the fundamental class $\left[M_+\right]\in H_d\left(M_+;{\Bbb Q}\right)$ 
one would like to have a map $R:M_+\rightarrow \mid BSL\left(\ff\right)\mid^+$ such that the following diagram 
%(with $h^M:M\rightarrow \mid B\pi_1M\mid$ the homotopy equivalence from Section 2.2 and $\mid B\rho\mid$ induced by $\rho$)
commutes up to homotopy:
\[\begin{xy}
\xymatrix{M \ar[r]^q \ar[d]_{\mid B\rho\mid h^M} & M_+ \ar[d]_{R}\\
\mid BSL\left(\ff\right)\mid\ar[r]^{incl}& \mid BSL\left(\ff\right)\mid^+}
\end{xy}\]
If this is the case, 
then
one can consider $R_*\left[M_+\right]$ and use the isomorphism $ H_d\left(\mid BSL\left(\ff\right)\mid^+;{\Bbb Q}\right)\cong H_d\left(\mid BSL
\left(\ff\right)\mid;{\Bbb Q}\right)$ 
%(the inverse of Quillen's isomorphism from Section 2.1)
to define 
an element
%$I_d^{-1}R_d\left[M_+\right]\in 
in $H_d\left(\mid BSL\left(\ff\right)\mid;{\Bbb Q}\right)$. 
%(And 
%thus, if the assumptions of \hyperref[proj]{Lemma \ref*{proj}} are satisfied for $A=\ff$, one obtains an element in $K_d\left(\ff\right)\otimes{\Bbb Q}$).

\begin{lem}\label{mapR} Let $M$ be a manifold with boundary such that 
$M$ and the path-components $\partial_1M,\ldots,\partial_sM$ of
$\partial M$ are aspherical. Let $q:M\rightarrow M_+$ be the canonical projection.

Let $\ff\subset {\Bbb C}$ a subring and $\rho:\pi_1M\rightarrow SL\left(N,\ff\right)$ be a representation such that $\rho\left(\pi_1\partial_i M\right)$ is unipotent for $i=1,\ldots,s$. 

Then there exists a
continuous map $R:M_+\rightarrow \mid BSL\left(N,\ff\right)\mid^+$ such that $$R\circ q=incl\circ \mid B\rho\mid \circ h^M,$$
where $incl: \mid BSL\left(N,\ff\right)\mid\rightarrow \mid BSL\left(N,\ff\right)\mid^+$ is the inclusion.
\end{lem}

\begin{pf}
Let $F$ be the homotopy fiber of $\mid BSL\left(N,
\ff\right)\mid \rightarrow \mid BSL\left(N,\ff\right)\mid^+$. 
It is well-known (e.g.\ \cite[page 
336]{cj}) that $\pi_1F$ is isomorphic to the Steinberg group $St\left(N;\ff\right)$. Let $\Phi:St\left(N,\ff\right)\rightarrow SL\left(N,\ff\right)$ be the canonical homomorphism.

By assumption, $\rho$ maps $\pi_1\partial_1 M$ 
into some maximal unipotent subgroup $B\subset SL\left(n,{\Bbb F}\right)$ of parabolic elements.
$B$ is conjugate to $B_0\subset SL\left(n,{\Bbb F}\right)$, the group of upper triangular matrices with all diagonal entries equal to 1. By \cite[Lemma 4.2.3]{ro} 
there exists a homomorphism $\Pi:B_0\rightarrow St\left(N,\ff\right)$ with $\Phi\Pi=id$. Applying conjugations and composing with $\rho$, we get
a homomorphism $\tau:\pi_1\partial_1M\rightarrow St\left(N;
\ff\right)$ such that $\Phi\tau=\rho\mid_{\pi_1\partial_1 M}$.

$\partial_1 M$ is aspherical,
hence $\tau$ is induced by some continuous mapping $g_1:\partial_1 M\rightarrow F$, and the diagram

$$\begin{xy}\xymatrix{
\partial_1 M \ar[r]^{i_1} \ar[d]^{g_1}&
M \ar[d]^{\mid B\rho\mid h^M}\\
F\ar[r]^{j}&
\mid BSL\left(N,{\Bbb F}\right)\mid}
\end{xy}$$

commutes up to some homotopy $H_t$. 
%(Here $\Gamma=\pi_1M$ and $B\rho:B\Gamma\rightarrow BSL\left({\Bbb F}\right)$ is induced by $\rho$.)

This construction can be repeated for all connected components $\partial_1M,\ldots,\partial_sM$ of $\partial M$.
For each $r=1,\ldots,s$ we get
a continuous map $g_r:\partial_rM\rightarrow F$ such that $jg_r\sim \mid
B\rho\mid h^M i_r$. Altogether, we get a continuous map $g:\partial M\rightarrow F$
such that $jg$ is homotopic to
$\mid
B\rho
\mid h^M i$.

By \cite[Lemma 8.1]{cj} this implies the existence of the desired map $R$.
\end{pf}\\

Hence one obtains an element in
%$I_d^{-1}R_d\left[M_+\right]\in 
$H_d\left(\mid BSL\left({\Bbb F}\right)\mid;{\Bbb Q}\right)$.
%\simeq K_*\left({\Bbb F}\right)\otimes{\Bbb Q}$. 
Unfortunately we did not succeed to prove its nontriviality, i.e.\ to evaluate the Borel class.
Therefore we will in the remainder of Section 4 pursue a different approach, closer in spirit to \cite{gon}, but 
surrounding the problem that $\partial M$ may be disconnected. 
%The reason why the generalized Cisneros-Molina-Jones construction was included
%at this point is that \hyperref[mapR]{Lemma \ref*{mapR}} will be important for the proof of \hyperref[preimage]{Proposition \ref*{preimage}}.

We mention that another "basis-point independent" approach might use multicomplexes in the sense of Gromov, but also here we were able to evaluate the Borel class only in the case that there are 2 or less boundary 
components. Also, in the case of hyperbolic 3-manifolds, yet another approach is due to Neumann-Yang \cite{ny}. 
For hyperbolic 3-manifolds of finite volume, Zickert has given in \cite{zi} a direct construction of a fundamental class $\left[M,\partial M\right]\in H_3\left(SL\left(2,{\Bbb C}\right),B_0\right)$, even in the case of possibly disconnected boundary.
It should be interesting to generalize and
compare the different constructions.

\subsection{Cuspidal completion}

\subsubsection{Disjoint cone}

We start with a {\bf notational remark}: the notion of disjoint cone for topological spaces resp.\ simplicial sets. This notion will be
useful for considering the
homology of a group relative to possibly more than one subgroup.

{\bf Disjoint cone of topological spaces.} Let $X$ be a topological space and $A_1,\ldots,A_s\subset X$ a set of (not necessarily disjoint) subspaces. There is a (not necessarily injective) continuous mapping $$i:A_1\dot{\cup}\ldots\dot{\cup}A_s\rightarrow X$$ from the {\bf disjoint} union $A_1\dot{\cup}\ldots\dot{\cup}A_s$ to $X$. 

We define the {\bf disjoint cone} $$DCone\left(\cup_{i=1}^sA_i\rightarrow X\right)$$ 
to be the pushout of the diagram
$$ \begin{xy}
\xymatrix{
A_1\dot{\cup}\ldots\dot{\cup}A_s\ar[r]^i \ar[d]&
X \ar[d]\\
Cone\left(A_1\right)\dot{\cup}\ldots\dot{\cup}Cone\left(A_s\right)\ar[r]&DCone\left(\cup_{i=1}^sA_i\rightarrow X\right)}
\end{xy}$$

If $X$ is a CW-complex and $A_1,\ldots,A_s$ are disjoint sub-CW-complexes, then for $*\ge 2$
$$H_*\left(DCone\left(\cup_{i=1}^sA_i\rightarrow X\right)\right)\cong
H_*\left(Cone\left(\cup_{i=1}^sA_i\rightarrow X\right)\right)=H_*\left(X,\cup_{i=1}^s A_i\right).$$
%A special case is that of a compact manifold $M$ with disconnected boundary $\partial M$, consisting of path-components 
%$\partial_1 M\dot{\cup}\ldots\dot{\cup}\partial_s M$. Then $DCone\left(\cup_{i=1}^s\partial_iM\rightarrow M\right)$ 
%is the space $M_+$ from Section 4.1.
%(In this case, the union of components is a disjoint union. Nonetheless $DCone$ is different from $Cone$.)
{\bf Disjoint cone of simplicial sets.}
We will need the cuspidal completion of 
a classifying space, which fits into the setting of simplicial sets.
(The point of the construction is that it may remember the geometry 
of the cusps of locally symmetric spaces. Thus it will serve as a technical device to handle the cusped case.)

For a simplicial set $\left(B,\partial_B,s_B\right)$ and a symbol $c$, {\em the cone over $B$ with cone point $c$} is the quasi-simplicial set whose $k$-simplices are 

- either k-simplices in $B$,

- or cones over $k-1$-simplices in $B$ with cone point $c$.\\
({\em By convention, the cone point is always the {\bf last} vertex of the cone over a $k-1$-simplex.})

The boundary operator $\partial$ in $Cone\left(B\right)$ is defined by $ \partial\sigma=\partial_B\sigma$ if $\sigma\in B$ and $\partial Cone\left(\tau\right)=Cone\left(\partial_B\tau\right)+\left(-1\right)^{dim\left(\tau\right)+1}\tau$ if $\tau\in B$. 

If $Y$ is a simplicial set and $\left\{B_i:i\in I\right\}$ a family of simplicial subsets indexed over a set $I$, then we define the quasi-simplicial set $DCone\left(\cup_{i\in I}B_i\rightarrow Y\right)$ 
%to be the 
%quasi-simplicial set which is the
as the push-out 
%of the diagram
$$ \begin{xy}
\xymatrix{
\dot{\cup}_{i\in I} B_i\ar[r] \ar[d]&
Y \ar[d]\\
\dot{\cup}_{i\in I} Cone\left(B_i\right)\ar[r]&DCone\left(\cup_{i\in I}B_i\rightarrow Y\right)}
\end{xy}$$

\subsubsection{ Construction of $BG^{comp}$ and $B\Gamma^{comp}$}

%As a technical device to handle the cusped case we introduce the cuspidal completion of $BG$. 
We recall from the beginning of Section 2.1 that $BG$ is the simplicial set realizing the bar construction. 
Thus its $k$-simplices are of the form $\left(g_1,\ldots,g_k\right)$ with $g_1,\ldots,g_k\in G$. We recall that $\partial_\infty\left(G/K\right)$ denotes the ideal boundary of $G/K$. The point of the following definition is that it allows to consider the geometry at each $c\in \partial_\infty\left(G/K\right)$ separately.
(As in the remark after \hyperref[homeo]{Corollary \ref*{homeo}} the definition of $B\Gamma^{comp}$ will assume
%$\widetilde{int\left(M\right)}=G/K$ were a symmetric space with
$rank\left(G/K\right)=1$.)

\begin{df}\label{cuspidal} Let $G/K$ be a symmetric space of noncompact type.
% and $\tilde{x}\in G/K$. 
%Let $n=dim\left(G/K\right)$.
%For each $c\in\partial_\infty G/K$ let $X_c\subset BG$ 
%be the
%subcomplex generated (under the face operation) by 
%$\left(1,p_1,\ldots,p_{n-1}\right)\in \left(BG\right)_{n-1}: vol\left(str\left(\tilde{x},p_1\tilde{x},\ldots,p_{n-1}\tilde{x}\right)\right)<\infty\right\}$.
We define the cuspidal completion $BG^{comp}$
of $BG$ to be 
%the pushout of the following diagram:
%$$ \begin{xy}
%\xymatrix{
%\dot{\cup}_{c\in\partial_\infty G/K} BG\ar[r] \ar[d]&
%BG \ar[d]\\
%\dot{\cup}_{c\in\partial_\infty G/K} Cone\left(BG\right)\ar[r]&BG^{comp}}
%\end{xy}$$
%\end{df}
%In other words, $BG^{comp}$
%is the mapping cone\footnotemark\footnotetext[5]{The mapping cone is a simplicial complex in the natural way:
%k-simplices of $BG^{comp}$ are 
%- either k-simplices in the target $BG$, 
%- or cones (with cone point $c\in\partial_\infty G/K$) over k-1-simplices in $BG$.\\
%More generally, if $X$ is a simplicial complex and $A_1,\ldots,A_s$ are simplicial subcomplexes, then $DCone\left(\cup_{i=1}^sA_i\rightarrow X\right)$ is a simplicial complex whose k-simplices are either k-simplices in $X$ or cones over k-1-simplices in some $A_i$.}
$$DCone \left(\dot{\cup}_{c\in\partial_\infty \left(G/K\right)} BG\rightarrow BG\right).$$  \end{df}
%but with the union $\dot{\cup}_{c\in\partial_\infty G/K} BG$ to be understood 
%{\bf as a disjoint union}.\\
\noindent
{\em Notation}: The cone point of 
%$Cone\left(BG\right)\subset BG^{comp}$ 
corresponding to $c\in\partial_\infty \left(G/K\right)$ will also be denoted by
$c$.

%\begin{df} Let $M$ be a compact manifold with boundary $\partial M=\partial_1M\cup\ldots\cup\partial_sM$ such that $int\left(M\right)=\Gamma\backslash G/K$ is a rank one locally symmetric space of noncompact type and finite volume. 
%%Let $\Gamma\subset G$ be a discrete subgroup such that $\Gamma\backslash G/K$ has finite volume. 
%Let $\Gamma c_i\subset\partial_\infty\left(G/K\right)$ be the set of parabolic fixed points corresponding to $\partial_iM$ by the remark after \hyperref[visibility]{Lemma \ref*{visibility}}.
%Then we define 
%$$B\Gamma^{comp}\subset BG^{comp}$$ to be the 

\begin{df}\label{conem}
Let $M$ be a manifold
with $\pi_1$-injective boundary $\partial M$, let $\partial_1M,\ldots,\partial_sM$ be the connected components of $\partial M$,
%satisfying the assumptions of \hyperref[straightsimplex]{Definition \ref*{straightsimplex}},
fix $x_0\in M$ and $x_i\in\partial_i M$ for $i=1,\ldots,s$, and let $\Gamma_i\subset \Gamma:=\pi_1\left(M,x\right)$ be defined according to
\hyperref[subgroup]{Definition \ref*{subgroup}}.

%Let $B\Gamma$ be the simplicial set as defined in Section 2.1 and $B\Gamma_i$ the naturally embedded simplicial subsets.
Assume that $M$ satisfies the assumptions of \hyperref[homeo]{Corollary \ref*{homeo}}
%carries a Riemannian metric of satisfies the assumptions of there exists a Riemannian metric of nonpositive secional curvature on $int\left(M\right)$ such that
%$\widetilde{int\left(M\right)}$ is a visibility manifold and $int\left/M\right)$ has 
%finite volume.
and let $c_i\in\partial_\infty\widetilde{int\left(M\right)}$ be the cusp associated to
$\Gamma_i$.
%Let $\Gamma c_i\subset\partial_\infty\left(G/K\right)$ be the set of parabolic fixed points corresponding to $\partial_iM$ by the remark after \hyperref[visibility]{Lemma \ref*{visibility}}.
Then we define
$$B\Gamma^{comp}=DCone\left(\cup_{i=1}^s B\Gamma_i\rightarrow B\Gamma\right)$$
%\subset BG^{comp}$$
%\subset
%DCone\left(\cup_{c\in\partial_\infty\left(G/K\right)}BG\rightarrow BG\right)$$ 
to be the
quasi-simplicial set whose $k$-simplices $\tau$
are either of the form $$\tau=\left(\gamma_1,\ldots,\gamma_k\right)$$
with $\gamma_1,\ldots,\gamma_k\in\Gamma$ or for some $i\in\left\{1,\ldots,s\right\}$ of the form
$$\tau=\left(p_1,\ldots,p_{k-1}, c_i\right)$$
with $p_1,\ldots,p_{k-1}\in\Gamma_i$.
% and $c_i\in\partial_\infty\widetilde{int\left(M\right)}$ corresponding to $\partial_iM$.
\end{df}

\noindent
{\em Notation}: The cone point 
of $Cone\left(B\Gamma_i\right)\subset B\Gamma^{comp}$ will be denoted by 
%$c_i$. This notation is suggested by the following observation, where $c_i$ is identified with the cusp 
$c_i\in\partial_\infty \left(G/K\right)$, the cusp associated to 
$\Gamma _i$. This is justified by the following observation.
%corresponding to $c_i\in\partial_\infty \widetilde{int\left(M\right)}$ is also denoted by $c_i$.\\

\begin{obs}\label{trivial}Let $M$ be a compact manifold with boundary $\partial M=\partial_1M\cup\ldots\cup\partial_sM$ such that $int\left(M\right)=\Gamma\backslash G/K$ is a locally symmetric space of noncompact type of rank one with finite volume.
Then $B\Gamma^{comp}\subset BG^{comp}$, where the cone point $c_i$ of $Cone\left(B\Gamma_i\right)$ corresponds to $c_i\in\partial_\infty\left(G/K\right)$ as the cone point of the corresponding copy of $Cone\left(BG\right)$.\end{obs}

{\em Remark}: $B\Gamma^{comp}$, as a subset of $BG^{comp}$, depends on the chosen identification of $\pi_1\left(\partial_iM,x_i\right)$ with a subgroup $\Gamma_i$ of $\Gamma$.

\subsubsection{Volume cocycle} In Section 2.3 we defined the volume cocycle $c\nu_d\in C_{simp}^d\left(BG\right)$ for a symmetric space $G/K$ of noncompact type. In this subsection we will extend $c\nu_d$ to $\overline{c\nu}_d\in C_{simp}^d\left(BG^{comp}\right)$.
  
For the remainder of this section we fix some $\tilde{x}\in G/K$. Let $d=dim\left(G/K\right)$.\\

We define the {\bf volume cocycle} $\overline{c\nu}_d\in C^d_{simp}\left(BG^{comp}\right)$ as follows.

For $\left(g_1,\ldots,g_d\right)\in BG$ we define 
$$\overline{c\nu}_d\left(g_1,\ldots,g_d\right)=algvol\left(str\left(\tilde{x},g_1\tilde{x},\ldots,g_1\ldots g_d\tilde{x}\right)\right)=\int_{str\left(\tilde{x},g_1\tilde{x},\ldots,g_1\ldots g_d\tilde{x}\right)}dvol$$
and for $\left(p_1,\ldots,p_{d-1},c\right)\in Cone\left(BG\right)$ with $c\in\partial_\infty \left(G/K\right)$ we define
$$\overline{c\nu}_d\left(p_1,\ldots,p_{d-1},c\right)=algvol\left(str\left(\tilde{x},p_1\tilde{x},\ldots,p_1\ldots
p_{d-1}\tilde{x},c\right)\right)=\int_{str\left(\tilde{x},p_1\tilde{x},\ldots,p_1\ldots
p_{d-1}\tilde{x},c\right)}dvol.$$
(This is defined because ideal d-simplices in a $d$-dimensional symmetric space $G/K$ of noncompact type have finite volume.)

The computation in Section 2.3 shows that $\delta \overline{c\nu}_d\left(g_1,\ldots,g_{d+1}\right)=0$ for $\left(g_1,\ldots,g_{d+1}\right)\in BG$. Moreover, for $\left(p_1,\ldots,p_d,c\right)\in Cone\left(BG\right)$ with
$c\in\partial_\infty\left(G/K\right)$ we have
$$\delta \overline{c\nu}_d\left(p_1,\ldots,p_d,c\right)= \overline{c\nu}_d
\left(\left(p_2,\ldots,p_d,c\right)+\sum_{i=1}^{d-1}\left(p_1,\ldots,p_ip_{i+1},\ldots,p_{d},c\right)+\left(-1\right)^{d+1}\left(p_1,\ldots,p_d\right)\right)$$
$$=...=\int_{\partial str\left(\tilde{x},p_1\tilde{x},\ldots,p_1\ldots p_d\tilde{x},c\right)}dvol=\int_{ str\left(\tilde{x},p_1\tilde{x},\ldots,p_1\ldots p_d\tilde{x},c\right)} d\left(dvol\right)=0.$$
This proves that $\overline{c\nu}_d$ is a simplicial cocycle on $BG^{comp}$.
Let $\overline{cv}_d=\left[\overline{c\nu}_d\right]\in H^d_{simp}\left(BG^{comp}\right).$ 
%be its cohomology class. 

By construction we have $\overline{c\nu}_d\mid_{BG}=c\nu_d$ and thus $\overline{cv}_d\mid_{BG}=comp\left(v_d\right)$ for the volume class $v_d=\left[\nu_d\right]\in H_c^d\left(G;{\Bbb R}\right)$ defined in Section 2.3.\\
\\
The Borel class $b_d\in H_c^d\left(GL\left({\Bbb C}\right);{\Bbb R}\right)$ defined in Section 2.4 may 
also be considered as a class $b_d\in H_c^d\left(SL\left({\Bbb C}\right);{\Bbb R}\right)$. 
For a representative $\beta_d\in C_c^d\left(SL\left({\Bbb C}\right);{\Bbb R}\right)$ of $b_d$ we define $c\beta_d\in C_{simp}^d\left(BSL\left({\Bbb C}\right);{\Bbb R}\right)$
by 
$$c\beta_d\left(g_1,\ldots,g_d\right):=\beta_d\left(1,g_1,g_1g_2,\ldots,g_1g_2\ldots g_d\right).$$
Then $c\beta_d$ represents $$comp\left(b_d\right)\in H_{simp}^d\left(BSL\left({\Bbb C}\right);{\Bbb R}\right).$$
%for the comparison map $comp$ defined in Section 2.2.
%We will use the same letter for the restriction of $comp\left(b_d\right)$ to $BSL\left(N,{\Bbb C}\right)$, that is we will consider $b_d$ as
%a class $b_d\in H^d\left(BSL\left(N,{\Bbb C}\right);{\Bbb R}\right)$.
\begin{lem}\label{fb} Let $d,N\in{\Bbb N}$ with $d$ odd. Let $G/K$ be a $d$-dimensional symmetric space of noncompact type.
If $\rho:\left(G,K\right)\rightarrow \left(SL\left(N,{\Bbb C}\right),SU\left(N\right)\right)$ is a
representation, then
there exists a quasi-simplicial set $BSL\left(N,{\Bbb C}\right)^{fb}$ with $$BSL\left(N,{\Bbb C}\right)\subset BSL\left(N,{\Bbb C}\right)^{fb}\subset
BSL\left(N,{\Bbb C}\right)^{comp}$$
and a homomorphism 
$$\overline{c\beta}_d: 
C_d^{simp}\left(BSL\left(N,{\Bbb C}\right)^{fb};{\Bbb R}\right)
\rightarrow {\Bbb R},$$
such that \\
i) $ \overline{c\beta}_d\mid_{C_d^{simp}\left(BSL\left(N,{\Bbb C}\right);{\Bbb R}\right)}$
is a cocycle representing $comp\left(b_d\right)$,\\
ii) 
%if $G/K$ is a $d$-dimensional symmetric space of noncompact type and $\rho:G\rightarrow SL\left(N,{\Bbb C}\right)$ a 
%representation, then 
we have 
$$\left(B\rho\right)_d\left(C_d^{simp}\left(BG^{comp};{\Bbb R}\right)\right)\subset C_d^{simp}\left(BSL\left(N,{\Bbb C}\right)^{fb};{\Bbb R}\right)$$ 
and
$\rho^*\overline{c\beta}_d$ represents $c_\rho \overline{cv}_d.$ 
(In particular, $\overline{c\beta}_d$ is well-defined on $\left(B\rho\right)_dH_d\left(BG^{comp};{\Bbb R}\right)$.)\end{lem}
\noindent
Here, $c_\rho\in {\Bbb R}$ is defined by the equality $\rho^*b_d=c_\rho v_d\in H_c^d\left(G;{\Bbb R}\right)$
%, cf.\ the proof of
from \hyperref[Thm2]{Theorem \ref*{Thm2}}.
%holds 
%for each representation $\rho:G\rightarrow SL\left(N,{\Bbb C}\right)$, where $c_\rho$ is defined by Theorem 2.\end{lem}

\begin{pf} %By the van Est Theorem 
%(Section 2.4.2)
%there is an isomorphism $$I:H^d_c\left(SL\left(N,{\Bbb C}\right)\right)
%\rightarrow H^d\left(sl\left(N,{\Bbb C}\right),su\left(N\right)\right),$$ where $H^*
%\left(sl\left(N,{\Bbb C}\right),su\left(N\right)\right)$ 
%is the cohomology of the complex of $SL\left(N,{\Bbb C}\right)$-invariant differential forms on $SL\left(N,{\Bbb C}\right)/SU\left(N\right)$. 
Let $dbol$ be\footnotemark\footnotetext[4]{The reason for the notation '$dbol$' is that $dbol$ relates to the Borel class $b_n$ as $dvol$ relates to $v_n$.

The superscript '$fb$' in $BSL\left(N,{\Bbb C}\right)^{fb}$ stands for 'finite Borel class'.} an
$SL\left(N,{\Bbb C}\right)$-invariant differential form on $SL\left(N,{\Bbb C}\right)/SU\left(N\right)$ representing $\nu_{SL\left(N,{\Bbb C}\right)}\left(b_d\right)$, where $\nu_{SL\left(N,{\Bbb C}\right)}:H^d_c\left(
SL\left(N,{\Bbb C}\right);{\Bbb R}\right)\rightarrow H^d\left(sl\left(N,{\Bbb C}\right),su\left(N\right)\right)$ from Section 2.4.2.
%differential form representing $I\left(b_d\right)$. 
Then a representative $\beta_d$ of $b_d$ is given by
$$\beta_d\left(g_0,g_1,\ldots,g_d\right):=\int_{str\left(g_0\tilde{x},g_1\tilde{x},\ldots,g_d\tilde{x}\right)}dbol$$
for each $\left(g_0,g_1,\ldots,g_d\right)\in \left(SL\left(N,{\Bbb C}\right)\right)^{d+1}$. (This follows from the explicit description of the van Est isomorpism in \cite[Theorem 1.1]{dup}.)

Since the van Est isomorphism is functorial and $\rho^*b_d=c_\rho v_d$,
we have that $\rho^* dbol-c_\rho dvol$
is an exact differential form. Moreover, $\rho^*dbol$ and $dvol$ are $G$-invariant differential forms on $G/K$. Hence they are harmonic and 
$\rho^* dbol-c_\rho dvol$
is an exact harmonic form, thus zero and we conclude 
$$\rho^* dbol=c_\rho dvol.$$

Let $\widetilde{X}=SL\left(N,{\Bbb C}\right)/SU\left(N\right)$. 
Define
$$BSL\left(N,{\Bbb C}\right)^{fb}_d:=
BSL\left(N,{\Bbb C}\right)_d\cup$$
$$\bigcup_{c\in\partial_\infty \widetilde{X}}\left\{\left(p_1,\ldots,p_{d-1},c\right)
\in Cone\left(BSL\left(N,{\Bbb C}\right)\right): \int_{str\left(\tilde{x},p_1\tilde{x},\ldots,p_1\ldots p_{d-1}\tilde{x},c\right)}dbol<\infty\right\}.$$
%(Recall that $\int_{\left(1,p_1,\ldots,p_{d-1},c\right)}p^*dbol=\int_{str\left(\tilde{x},p_1\tilde{x},\ldots,p_1\ldots p_{d-1}\tilde{x},c\right)}dbol$.)\\

This defines the d-simplices of $BSL\left(N,{\Bbb C}\right)^{fb}$ and
we define $BSL\left(N,{\Bbb C}\right)^{fb}$ 
to be the quasi-simplicial set generated by $BSL\left(N,{\Bbb C}\right)^{fb}_d$ under face maps. \\

Define $\overline{c\beta}_d:BSL\left(N,{\Bbb C}\right)^{fb}_d\rightarrow {\Bbb R}$ by 
$$\overline{c\beta}_d\left(g_1,\ldots,g_d\right)=\int_{str\left(\tilde{x},g_1\tilde{x},
\ldots,g_1\ldots g_d\tilde{x}\right)}dbol$$
if $\left(g_1,\ldots,g_d\right)\in BSL\left(N,{\Bbb C}\right)_d$, and
$$\overline{c\beta}_d\left(p_1,\ldots,p_{d-1}\right)=\int_{str\left(\tilde{x},\ldots,p_1\tilde{x},p_1\ldots p_{d-1}
\tilde{x},c\right)}dbol$$
if $\left(p_1,\ldots,p_{d-1}\right)\in BSL\left(N,{\Bbb C}\right)_{d-1}$, $c\in\partial_\infty \widetilde{X}$ and $\left(p_1,\ldots,p_{d-1},c\right)\in BSL\left(N,{\Bbb C}\right)_{d}^{fb}$.
%$\overline{b}_n$ is given by integration of a closed form, hence it vanishes on $\partial BSL\left(N,{\Bbb C}\right)_{n+1}\cap BSL\left(N,{\Bbb C}\right)^{fb}_n$ and
%gives a well-defined homomorphism
%$$\overline{b}_n: im\left(H_n\left(BSL\left(N,{\Bbb C}\right)^{fb};{\Bbb R}\right)\rightarrow H_n\left(BSL\left(N,{\Bbb C}
%\right)^{comp};{\Bbb R}\right)\right)\rightarrow {\Bbb R},$$

By construction, $\overline{c\beta}_d\mid_{C_d^{simp}\left(BSL\left(N,{\Bbb C}\right);{\Bbb R}\right)}$ is a cocyle representing $comp\left(b_d\right)$. 

%Let $G/K$ be an $n$-dimensional symmetric space of noncompact type and $\rho:G\rightarrow SL\left(N,{\Bbb C}\right)$ a
%representation. 
%To prove b), let $z\in C_d^{simp}\left(BG^{comp};{\Bbb R}
%\right)$.
% be a cycle, i.e.\ $dz=0$.
The homomorphism $\rho:G\rightarrow SL\left(N,{\Bbb C}\right)$ extends to a well-defined map $\partial_\infty G/K\rightarrow \partial_\infty\widetilde{X}$, thus we obtain a well-defined simplicial map $B\rho:BG^{comp}\rightarrow BSL\left(N,{\Bbb C}\right)^{comp}$.
From $\rho^* dbol=c_\rho dvol$ we have
%$$\int_{\rho\left(z\right)}dbol=c_\rho\int_z dvol.$$
$\left(B\rho\right)_d\left(C_d^{simp}\left(BG^{comp};{\Bbb R}\right)\right)\subset C_d^{simp}\left(BSL\left(N,{\Bbb C}\right)^{fb};{\Bbb R}\right)$ and
$\rho^*\overline{c\beta}_d$ represents $c_\rho \overline{cv}_d.$ \end{pf}

%For a simplex $\left(1,p_1,\ldots,p_{n-1},c\right)\in Cone\left(BSL\left(N,{\Bbb C}\right)\right)$, we define $\overline{b}_n$ by
%$$\overline{b}_n\left(1,p_1,\ldots,p_{n-1},c\right):=\int_{str\left(\tilde{x},p_1\tilde{x},\ldots,p_{n-1}\tilde{x},c\right)}dbol.$$
%It is easily checked that $d\overline{b}_n=0$ and that $\overline{b}_n=b_n$ on $BSL\left(N,{\Bbb C}\right)$.

%Thus $\overline{b}_n$ extends $b_n$ to a simplicial cocycle on $BSL\left(N,{\Bbb C}\right)^{comp}$.
%Moreover, $\rho^* dbol=c_\rho dvol$ implies $$\rho^*\overline{b}_n=c_\rho \overline{v}_n.$$\end{pf}

%Remark: In particular, if $\rho:G\rightarrow SL\left(N,{\Bbb C}\right)$ is a representation of a semisimple Lie group and $n=dim\left(G/K\right)$, then $\overline{b}_n$ takes finite values on simplices in the image of $\rho$.

\begin{df}\label{subfield} Let ${\Bbb F}\subset{\Bbb C}$ be a subring (with 1) and $G/K$ a symmetric space of noncompact type. Then 
we define 
%$$BSL\left(N,{\Bbb F}\right)^{comp}=DCone\left(\dot{\cup}_{c\in\partial_\infty\left(SL\left(N,{\Bbb C}\right)/SU\left(N\right)\right)}
%BSL\left(N,{\Bbb F}\right)\rightarrow BSL\left(N,{\Bbb F}\right)\right)
%\subset BSL\left(N,{\Bbb C}\right)^{comp}$$
$$BG\left({\Bbb F}\right)^{comp}=DCone\left(\dot{\cup}_{c\in\partial_\infty\left(G/K\right)}
BG\left({\Bbb F}\right)\rightarrow BG\left({\Bbb F}\right)\right)
\subset BG^{comp}.$$
For $G=SL\left(N,{\Bbb C}\right)$ we define $$BSL\left(N,{\Bbb F}\right)^{fb}=BSL\left(N,{\Bbb C}\right)^{fb}\cap BSL\left(N,{\Bbb F}\right)^{comp}.$$
%Moreover we let $$BSL\left({\Bbb F}\right)^{fb}=\cup_{N\in{\Bbb N}}BSL\left(N,{\Bbb F}\right)^{fb}$$ be the increasing union with respect to the obvious inclusion
\end{df}

\subsection{Straightening of interior and ideal simplices}

The purpose of this section is to describe an explicit realization of the isomorphism
$$H_*\left(DCone\left(\cup_{i=1}^s\partial_iM\rightarrow M\right)\right)\cong
H_*^{simp}\left(DCone\left(\cup_{i=1}^sB\Gamma_i\rightarrow B\Gamma\right)\right)$$
for $\Gamma=\pi_1M,\Gamma_1=\pi_1\partial_1M,\ldots,\Gamma_s=\pi_1\partial_sM$, under the assumptions of \hyperref[visibility]{Lemma \ref*{visibility}}, that 
is if $M$ is a finite-volume quotient of a nonpositively curved visibility manifold.

In Section 2.1 we used straightening to define the Eilenberg-MacLane map on genuine simplices. In this section
we will extend the Eilenberg-MacLane map to ideal simplices.\\
%We recall that we defined in Section 4.2.1 the notion of disjoint cone for simplicial sets. We remember that, by definition, the
%cone points are always the last vertex of cones over simplices.
\begin{df}\label{Chat} 
Let $M$ be a compact manifold with boundary, let $\partial_1 M,\ldots,\partial_sM$ be the connected components of $\partial M$. Let $x_0,x_i,\Gamma,\Gamma_i$ be defined according to \hyperref[subgroup]{Definition \ref*{subgroup}}. 
We denote
$$\widehat{C}_*\left(M\right):=C_*^{simp}\left(DCone\left(\cup_{i=1}^sC_*\left(\partial_i M\right)\rightarrow C_*\left(M\right)\right)\right).$$
For $i=1,\ldots,s$ let $C_i$ be the cone point of $Cone\left(C_*\left(\partial_iM\right)\right)$.
A vertex of a simplex in 
$\widehat{C}_*\left(M\right)$ is an ideal vertex, if it is in one of the cone points $C_1,\ldots,C_s$, and an interior
vertex else.
%For fixed $x_0\in M$, $x_i\in\partial_iM$, and a fixed identification of $\pi_1\left(\partial_iM,x_i\right)$ with a subgroup of $\pi_1\left(M,x_0\right)$ (using a path from $x_0$ to $x_i$), 
Then we define $$\widehat{C}_*^{x_0}\left(M\right)\subset\widehat{C}_*\left(M\right)$$
to be the subcomplex freely generated by those simplices for which\\
- either all vertices are in $x_0$,\\
- or the last vertex is an ideal vertex $C_i$,
%, corresponding to some boundary component
%$\partial_i M$, 
all other vertices are in $x_0$, and the homotopy
classes (rel.\ $\left\{0,1\right\}$)
of all edges
between interior vertices belong to $\Gamma_i\subset \pi_1\left(M,x_0\right)$.
\end{df}

By construction, $\widehat{C}_*\left(M\right)$ and $\widehat{C}_*^{x_0}\left(M\right)$
are chain complexes.

From now on we assume that the assumptions of \hyperref[cusps]{Corollary \ref*{cusps}} (and thus the assumptions of \hyperref[visibility]{Lemma \ref*{visibility}}) hold for $N=int\left(M\right)=M-\partial M$. In particular we have the projection 
$\overline{\pi}:
\widetilde{int\left(M\right)}\bigcup\cup_{i=1}^s
\Gamma c_i\rightarrow DCone\left(\cup_{i=1}^s\partial_i M\rightarrow M\right) $ from \hyperref[cusps]{Corollary \ref*{cusps}}.

\begin{df}\label{straightsimplex} 
%If $M$ is a manifold with boundary $\partial M=\partial_1M\cup\ldots
%\cup\partial_sM$ and that $int\left(M\right)$ has a Riemannian metric of finite volume such that $\widetilde{int\left(M\right)}$
%is a visibility manifold of nonpositive sectional curvature and $int\left(M\right)=\Gamma\backslash\widetilde{int\left(M\right)}$ for a discrete group of isometries $\Gamma$, if $\Gamma c_i$ denotes the $\Gamma$-orbit of the $c_i\in\partial_\infty\widetilde{int\left(M\right)}$ which corresponds to $\partial_iM$, and if
%the projection $\pi$ is defined according to the foregoing remark,
Let the assumptions of 
%\hyperref[visibility]{Lemma \ref*{visibility}} 
\hyperref[cusps]{Corollary \ref*{cusps}}
hold. A simplex in $DCone\left(\cup_{i=1}^s\partial_i M\rightarrow M\right)$ is said to be {\bf straight} if some (hence any) lift to $\widetilde{int\left(M\right)}\bigcup\cup_{i=1}^s\Gamma c_i\subset \widetilde{int\left(M\right)}\cup\partial_\infty\widetilde{int\left(M\right)}$ is straight.

In particular a $k$-simplex $\sigma\in \widehat{C}^*\left(M\right)$
%DCone\left(\cup_{i=1}^sC_*\left(\partial_i M\rightarrow M\right)$
is straight if it is either of the form
$$\sigma=\pi\left(str\left(\tilde{x}_0,\tilde{x}_1,\ldots,\tilde{x}_k\right)\right)$$
with $\tilde{x}_0,\tilde{x}_1,\ldots,\tilde{x}_k\in \widetilde{int\left(M\right)}$ or of the form
$$\sigma=\pi\left(str\left(\tilde{x}_0,\tilde{x}_1,\ldots,\tilde{x}_{k-1},\gamma c_i\right)\right)$$
with $\tilde{x}_0,\tilde{x}_1,\ldots,\tilde{x}_{k-1}\in \widetilde{int\left(M\right)},\gamma\in\Gamma,i\in\left\{1,\ldots,s\right\}$.
\end{df}

\begin{df}\label{straight} Let $M$ be a manifold satisfying the assumptions of \hyperref[straightsimplex]{Definition \ref*{straightsimplex}}.
%Assume that $M$ is a manifold with boundary $\partial M$ and that $int\left(M\right)$ has a Riemannian metric such that $\widetilde{int\left(M\right)}$
%is a visibility manifold of nonpositive sectional curvature.
Let $x_0\in M$. Then we define the chain complex
$$\widehat{C}_*^{str,x_0}\left(M\right):={\Bbb Z}\left[\left\{\sigma\in \widehat{C}_*^{x_0}\left(M\right): \sigma\mbox{\ straight}\right\}\right].$$\end{df}
%as the subcomplex generated by the straight simplices.\end{df}
%$\widehat{C}_*^{str,x_0}\left(M\right)$ is a chain complex because faces of straight simplices are straight.

\begin{lem}\label{iso} 
Let $M$ be a compact manifold with boundary, let $\partial_1 M,\ldots,\partial_sM$ be the connected components of $\partial M$. Let $x_0,x_i,\Gamma,\Gamma_i$ be defined according to \hyperref[subgroup]{Definition \ref*{subgroup}}.
Moreover let 
the assumptions of
%\hyperref[visibility]{Lemma \ref*{visibility}}
\hyperref[cusps]{Corollary \ref*{cusps}}
hold.

a) 
Then there is an isomorphism of chain complexes
 $$\Phi:\widehat{C}_*^{str,x_0}\left(M\right)\rightarrow C_*^{simp}\left(B\Gamma^{comp}\right).$$

b) The inclusion $$\widehat{C}_*^{str,x_0}
\left(M\right)\rightarrow \widehat{C}_*\left(M\right)$$
%\rightarrow C_*\left(DCone\left(\cup_{i=1}^s\partial_iM\rightarrow M\right)\right)$$
is a chain homotopy equivalence.

c) The
composition of $\Psi:=\Phi^{-1}$ with the inclusion $\widehat{C}_*^{str,x_0}
\left(M\right)\rightarrow 
%\widehat{C}_*\left(M\right)\rightarrow 
C_*\left(DCone\left(\cup_{i=1}^s\partial_iM\rightarrow M\right)\right)$
induces an isomorphism
$$EM_*:
H_*^{simp}\left(B\Gamma^{comp}\right)\rightarrow   H_*\left(DCone\left(\cup_{i=1}^s\partial_iM\rightarrow M\right)\right).$$
\end{lem}
\begin{pf}
a)
In Section 2.1 we defined a chain isomorphism $\Phi:C_*^{str,x_0}\left(M\right)\rightarrow C_*^{simp}\left(B\Gamma\right)$ by $\Phi\left(\sigma\right)=\left(g_1,\ldots,g_k\right)$, where $\sigma\in
C_k^{str,x_0}\left(M\right)$ is a continuous map $\sigma:\Delta^k\rightarrow M$ with $\sigma\left(w_j\right)=x_0$ for $j=0,\ldots,k$,
and $g_j\in\Gamma=\pi_1\left(M,x_0\right)$ is the homotopy class (rel.\ vertices) of $\sigma\mid_{\gamma_j}$ for $j=1,\ldots,k$.
 Moreover we defined a chain isomorphism $\Psi:C_*^{simp}\left(B\Gamma\right)\rightarrow C_*^{str,x_0}\left(M\right)$ by
$\Psi\left(g_1,\ldots,g_k\right):=\pi\left(str\left(\tilde{x}_0,g_1\tilde{x}_0,g_1g_2\tilde{x}_0,\ldots,g_1
\ldots g_k\tilde{x}_0\right)\right)$ and we proved $\Phi\Psi=id$ and $\Psi\Phi=id$. We will now extend $\Phi$ and $\Psi$ to chain isomorphisms
$$\Phi: \widehat{C}_*^{str,x_0}\left(M\right)\rightarrow 
C_*^{simp}\left(B\Gamma^{comp}\right),$$
$$\Psi:
C_*^{simp}\left(B\Gamma^{comp}\right)
\rightarrow
 \widehat{C}_*^{str,x_0}\left(M\right)$$ and will prove that the extensions are inverse to each other.

Let $\sigma\in \widehat{C}_k^{str,x_0}\left(M\right)$ be a straight $k$-simplex which is not in $C_k^{str,x_0}\left(M\right)$. This means that
%its last vertex is an ideal vertex $C_j$, corresponding to the cone point of the cone over $C_*\left(\partial_jM\right)$ for some $j\in\left\{1,\ldots,s\right\}$.
the lift $\tilde{\sigma}$ of $\sigma$ to 
$$\widetilde{int\left(M\right)}\bigcup\cup_{i=1}^s\Gamma c_i\subset\widetilde{int\left(M\right)}\cup\partial_\infty\widetilde{int\left(M\right)}$$
is 
of the form 
$$\tilde{\sigma}=
\pi\left(str\left(
\gamma_0\tilde{x}_0,\gamma_1\tilde{x}_0,\ldots,\gamma_{k-1}\tilde{x}_0,\gamma c_i\right)\right)$$
for some $i\in\left\{1,\ldots,s\right\}$ and some $\gamma_0,\ldots,\gamma_{k-1},\gamma\in \Gamma_i$.
%a straight simplex with last vertex $v_k$ belonging to $\Gamma c_j$ and the remaining vertices $v_0,\ldots,v_{k-1}$ belonging to
%$\Gamma \tilde{x}_0$ for some lift $\tilde{x}_0\in \widetilde{int\left(M\right)}$ of $x_0$. Let $v_i=\gamma_i\tilde{x}_0$ for $i=0,\ldots,k-1$ and 
We define
$$\Phi\left(\sigma\right)=\left(\gamma_1\gamma_0^{-1},\ldots,\gamma_{k-1}\gamma_{k-2}^{-1}, c_i\right),$$
where $ c_i$ is the cone point of $Cone\left(B\Gamma_i\right)$. 
%(This is well-defined because $\gamma_i\gamma_{i-1}^{-1}$
%does not depend on the choice of a lift $\tilde{x}_0$.)

Conversely, if a simplex 
$\tau\in 
C_*^{simp}\left(B\Gamma^{comp}\right)
$ does not belong to $C_*^{simp}\left(B\Gamma\right)$ 
then $\tau\in Cone\left(B\Gamma_i\right)$ for some $i\in\left\{1,\ldots,s\right\}$, but $\tau\not\in B\Gamma_i$, thus $\tau$
is of the form 
$$\tau=\left(p_1,\ldots,p_{k-1}, c_i\right)\in 
C_*^{simp}\left(DCone\left(\cup_{i=1}^s
%\cup_{\gamma\in \Gamma} 
B\Gamma_i\rightarrow B\Gamma\right)\right)$$ 
for some $i\in\left\{1,\ldots,s\right\}$,
with $p_1,\ldots,p_{k-1}\in \Gamma_i$ and $ c_i$ the cone point of $Cone\left(B\Gamma_i\right)$.
Then we define 
$$\Psi\left(\tau\right)=\pi\left(str\left(\tilde{x}_0,p_1\tilde{x}_0,\ldots,p_1\ldots p_{k-1}\tilde{x}_0, c_i\right)\right)\in \widehat{C}_*^{str,x_0}\left(M\right).$$

From Section 2.1 we have $\Psi\left(\partial \tau\right)=\partial \Psi\left(\tau\right)$ for $\tau\in C_*^{simp}\left(B\Gamma\right)$.
On the other hand, if $\tau  =\left(p_1,\ldots,p_{k-1}, c_i\right)\in 
C_*^{simp}\left(B\Gamma^{comp}\right)$, then a straightforward computation shows
$$\Psi\left(\partial \tau\right)-\partial\Psi\left(\tau\right)=
%$$\Psi \left(\partial\left(p_1,\ldots,p_{k-1}, c_i\right)\right)=$$
%$$\Psi\left(p_2,\ldots,p_{k-1}, c_i\right)+
%\sum_{i=1}^{k-2}\Psi\left(p_1,\ldots,p_ip_{i+1},\ldots,p_{k-1}, c_i\right)+
%\left(-1\right)^{k-1}\Psi\left(p_1,\ldots,p_{k-2}, c_i\right)+\left(-1\right)^k\Psi\left(p_1,\ldots,p_{k-1}\right)$$
\pi\left(str\left(\tilde{x}_0,p_2\tilde{x}_0,\ldots,p_2\ldots p_{k-1}\tilde{x}_0, c_i\right)\right)- \pi\left(str\left(p_1\tilde{x}_0,p_1p_2\tilde{x}_0,\ldots,p_1p_2\ldots p_{k-1}\tilde{x}_0, c_i\right)\right).$$
%$$+
%\sum_{i=1}^{k-2}\left(-1\right)^i\pi\left(str\left(\tilde{x}_0,\ldots,p_1\ldots p_{i-1}\tilde{x}_0,p_1\ldots p_ip_{i+1}\tilde{x}_0,\ldots,p_1\ldots p_{k-1}\tilde{x}_0, c_i\right)\right)$$
%$$+\left(-1\right)^{k-1}\pi\left(str\left(\tilde{x}_0,p_1\tilde{x}_0,\ldots,p_1\ldots p_{k-2}\tilde{x}_0,c_i\right)\right)+
%\left(-1\right)^k\pi\left(str\left(\tilde{x}_0,p_1\tilde{x}_0,\ldots,p_1\ldots p_{k-1}\tilde{x}_0\right)\right)$$
%$$=
%\pi\left(str\left(p_1\tilde{x}_0,p_1p_2\tilde{x}_0,\ldots,p_1p_2\ldots p_{k-1}\tilde{x}_0, c_i\right)\right)$$
%$$+
%\sum_{i=1}^{k-2}\left(-1\right)^i\pi\left(str\left(\tilde{x}_0,\ldots,p_1\ldots p_{i-1}\tilde{x}_0,p_1\ldots p_ip_{i+1}\tilde{x}_0,\ldots,p_1\ldots p_k\tilde{x}_0, c_i\right)\right)$$
%$$+\left(-1\right)^{k-1}\pi\left(str\left(\tilde{x}_0,p_1\tilde{x}_0,\ldots,p_1\ldots p_{k-2}\tilde{x}_0, c_i\right)\right)
%+\left(-1\right)^k \pi\left(str\left(\tilde{x}_0,p_1\tilde{x}_0,\ldots,p_1\ldots p_{k-1}\tilde{x}_0\right)\right)$$
%$$=\pi\left(\partial str\left(\tilde{x}_0,p_1\tilde{x}_0,\ldots,p_1\ldots p_{k-1}\tilde{x}_0, c_i\right)\right)=
%\partial \Psi\left(p_1,\ldots,p_{k-1}, c_i\right),$$ where we have used that 
Thus $p_1\in\Gamma_i\subset Fix\left(c_i\right)$
and $\Gamma$-invariance of $\pi\left(str\left(.\right)\right)$ implies
%$\pi\left(str\left(\tilde{x}_0,p_2\tilde{x}_0,\ldots,p_2
%\ldots p_{k-1}\tilde{x}_0,c_i\right)\right)=\pi\left(str\left(p_1\tilde{x}_0,p_1p_2\tilde{x}_0,\ldots,p_1p_2\ldots p_{k-1}\tilde{x}_0,c_i\right)\right)$ for each deck transformation $p_1\in \Gamma_i$.
$\Psi\left(\partial\tau\right)=\partial\Psi\left(\tau\right)$, that is
$\Psi$ is a chain map.

Clearly $\Phi\left(\pi\left(str\left(\tilde{x}_0,p_1\tilde{x}_0,\ldots,p_1\ldots p_{k-1}\tilde{x}_0,c_i\right)\right)\right)=\left(p_1,\ldots,p_{d-1},c_i\right)$, thus $\Phi\Psi=id$. On the other hand,
a straight simplex $\sigma:\Delta^k\rightarrow M$ with the first $k$ vertices in $x_0$ and the last vertex in $\Gamma c_i$
is uniquely determined by the homotopy classes (rel.\ vertices) of
$p_j=\left[\sigma\mid_{\gamma_j}\right]$ for $j=1,\ldots,k-1$, because its lift to $\widetilde{M}$ must be in the $\Gamma$-orbit of
$str\left(\tilde{x}_0,p_1\tilde{x}_0,\ldots,p_1\ldots p_{k-1}\tilde{x}_0,c_i\right)$. Thus $\Psi\Phi=id$. This shows that $\Psi$ and $\Phi$ are inverse to each other, in particular both are chain isomorphisms.\\
\\
b) We define a chain homotopy $\widehat{C}_*\left(M\right)\rightarrow\widehat{C}^{x_0}_*\left(M\right)$, left-inverse to the inclusion, by induction on the dimension of simplices. First, for each each $v\in C_0\left(\partial_iM\right)$ we fix 
a chain homotopy from $v$ to $x_i$ {\em inside $\partial _iM$}. The fixed path $l_i$ from \hyperref[subgroup]{Definition \ref*{subgroup}} 
provides us with a chain homotopy from $x_i$ to $x_0$. Composition of these two chain homotopies yields a chain homotopy from $v$ to $x_0\in 
\widehat{C}_0^{x_0}\left(M\right)$. If $v\in C_0\left(M\right)-C_0\left(\partial M\right)$, then we fix an arbitrary chain homotopy from $v$ to $x_0$. For the cone points fix the constant chain homotopy.
Now for each 1-simplex $e$ we have a chain homotopy of its vertices into either $x_0$ or one of the cone points. This chain homotopy of $\partial e$ can be extended to a chain homotopy of $e$. If $e$ had vertices in $\partial_iM$, then we observe that the chain homotopy of the vertices consisted 
of two steps. In the first step the vertices were homotoped {\em inside $\partial_i M$} into $x_i$. Thus $e$ can 
be homotoped inside $\partial_iM$ into a loop with vertices in $x_i$, which then represents an element of 
$\pi_1\left(\partial_iM,x_i\right)$. In the second step the vertices were homotoped along the $l_i$, thus $e$ can be homotoped into a loop representing an element of $\Gamma_i$ as defined in \hyperref[subgroup]{Definition \ref*{subgroup}}.
 Thus we have a chain homotopy from $\widehat{C}_1\left(M\right)$ to $\widehat{C}^{x_0}_1\left(M\right)$. 
A standard argument shows that this chain homotopy can be recursively extended to the $\widehat{C}_k\left(M\right)$ for all $k\in{\Bbb N}$.

%Assume that 
%$z=\sum_{i=1}^r\tau_i
%+\sum_{j=1}^p\kappa_j
%\in \widehat{C}_k\left(M\right)$ is a singular chain, where $\tau_i
%\in C_*\left(M\right)$
%and $\kappa_j
%\in Cone\left(C_*\left(\partial M\right)\right)$, and that the chain homotopy is alread defined on the $k-1$-skeleton.
%For $j=1,\ldots,p$ let $\partial_{i_j}M$ be the path component of $\partial M$ corresponding to the cone point which is the
%ideal vertex of $\kappa_j$. We can homotope $z$ such that all interior
%vertices of all $\kappa_j$ are homotoped into the respective $x_{i_j}\in\partial_{i_j}M$ and also all edges between interior vertices
%are homotoped into $\partial_{i_j}M$, the ideal vertices remaining fixed during the homotopy. Next we use the path from $x_i$ to $x_0$ (which realizes the given identification of $\pi_1\left(\partial_iM,x_i\right)$ with a subgroup of $\pi_1\left(M,x_0\right)$ in \hyperref[Chat]{Definition \ref*{Chat}}) and homotope $z$ such that all interior vertices of all $\kappa_j$ are homotoped into $x_0$. Finally we homotope $z$ such that all vertices of all $\tau_i$ are homotoped into $x_0$. The result is a cycle in $\widehat{C}_*^{x_0}\left(M\right)$ which is homotopic (hence chain-homotopic) to $z$.
We then apply the usual straightening procedure (\cite{bp}, Lemma C.4.3) to construct 
a
chain homotopy $\widehat{C}_*^{x_0}\left(M\right)\rightarrow \widehat{C}_*^{str,x_0}\left(M\right)$, left-inverse to the inclusion.\\
%This shows that $\widehat{C}_*^{str,x_0}
%\left(M\right)\rightarrow \widehat{C}_*\left(M\right)$ is surjective in homology.
%To prove injectivity in homology, let $z\in \widehat{C}_*^{str,x_0}
%\left(M\right)$ be a cycle and assume that $z=\partial w$ with $w\in  \widehat{C}_*\left(M\right)$. Since $\partial w$ is straight, we can apply the 
%same straightening procedure to
%homotope $w$, keeping $\partial w=z$ fixed,
%into some chain $str\left(w\right)\in\widehat{C}_*^{str,x_0}
%\left(M\right)$. Then $\partial str\left(w\right)=z $
%implies that $z$ is 0-homologous in $\widehat{C}_*^{str,x_0}
%\left(M\right)$.
\\
c) Comparison of the Mayer-Vietoris sequences shows that
inclusion $\widehat{C}_*\left(M\right)\rightarrow C_*\left(DCone\left(\cup_{i=1}^s\partial_i M\rightarrow M\right)\right)$
is a homology equivalence. Hence c) follows from a) and b).\end{pf} \\

Thus for $d$-dimensional compact, orientable Riemannian manifold $M$ of nonpositive sectional curvature, 
$EM_d^{-1}\left[M,\partial M\right]\in H_d^{simp}\left(DCone\left(\cup_{i=1}^s B\Gamma_i\rightarrow B\Gamma\right)\right)$ is well-defined.

\subsection{Construction of $\overline{\gamma}\left(M\right)$}

The results of the previous sections allow us to consider the image of the fundamental class $\left[M,\partial M\right]$ in $H_d^{simp}\left(DCone\left(\cup_{i=1}^s B\Gamma_i\rightarrow B\Gamma\right)\right)$ and its push-forward, for representations satisfying suitable assumptions, in $H_d^{simp}\left(BSL\left(N,{\Bbb F}\right)^{fb}\right)$. The aim of this subsection will be to show that this element has a (uniquely defined) preimage in $H_d^{simp}\left(BSL\left(N,{\Bbb F}\right)\right)$.
\begin{pro}\label{preimage}
 Let
$M$ be a compact, oriented, connected manifold with boundary components
$\partial_1M,\ldots,\partial_sM$ such that
 $Int\left(M\right)=\Gamma\backslash
         G/K$ is a locally
        symmetric space of noncompact type of rank one with finite volume. 

Fix $x_0\in M$ and $x_i\in\partial_i M$ for $i=1,\ldots,s$, and fix the isomorphisms of $\pi_1\left(\partial_iM,x_i\right)$ with subgroups $\Gamma_i$ of $\Gamma=\pi_1\left(M,x_0\right)$ given by \hyperref[subgroup]{Definition \ref*{subgroup}}.
Assume that, for some subring ${\Bbb F}\subset{\Bbb C}$, we have an inclusion
$$j:\left(\Gamma,\Gamma_i\right)\rightarrow \left(G\left({\Bbb F}\right),\Gamma_i\right).$$ 
Let $$\rho:G\left(\ff\right)\rightarrow SL\left(N,\ff\right)$$ 
be a representation. 
%let $\Gamma_i^\prime:=\rho\left(\Gamma_i\right)$ for 
%$i=1,\ldots,s$.\\
        Denote $$\left[M,\partial M\right]\in H_d\left(DCone\left(\cup_{i=1}^s\partial_i M\rightarrow M\right);{\Bbb Q}\right)$$
%\cong H_d^{simp}\left(
%DCone\left(\cup_{i=1}^s B\Gamma_i\rightarrow B\Gamma\right);{\Bbb Q}\right)$$
        the fundamental class of $M$. Then $$  B\left(\rho j\right)_*EM^{-1}\left[M,\partial 
M\right]\in H_d^{simp}\left(BSL\left(N,{\Bbb F}\right)^{fb}
%DCone\left(\cup_{i=1}^sB\Gamma_i^\prime\rightarrow BSL\left({\Bbb F}\right)
;{\Bbb Q}\right)$$ has a preimage
$$\overline{\gamma}\left(M\right)\in H_d^{simp}\left(BSL\left(N,{\Bbb F}\right);{\Bbb Q}\right).$$
The preimage does not depend on the chosen identification of $\pi_1\partial_iM$ with a subgroup $\Gamma_i\subset \Gamma$. 
\end{pro}

\begin{pf}  Let $\Gamma_i^\prime:=\rho\left(\Gamma_i\right)$ for
$i=1,\ldots,s$.
In a first step we will prove that the desired preimage exists if 
$$B\left(\rho j\right)_{*}EM^{-1}\left[\partial_i M\right]=0\in H_{d-1}^{simp}\left(
B\Gamma_i^\prime;                                                                                                                                                                    {\Bbb Q}\right)$$ for $i=1,\ldots,s$. In the second step we will then prove this equality. 
%It suffices the prove the existence of $\overline{\gamma}\left(M\right)$ for ${\Bbb F}={\Bbb C}$ because the so-defined $\overline{\gamma}\left(M\right)$ can be restricted to $H_d^{simp}\left(BSL\left(N,{\Bbb F}\right);{\Bbb Q}\right)$ for ${\Bbb F}\subset {\Bbb C}$.
%We prove the result for ${\Bbb Z}$-coefficients, the statement for 
%${\Bbb Q}$-coefficients follows by 
%application of the change-of-rings homomorphism. 
%We have $EM_d^{-1}\left[M,\partial M\right]\in H_d^{simp}\left(B\Gamma^{comp}\right)$, thus $Bj_dEM_d^{-1}\left[M,\partial M\right]\in H_d^{simp}\left(BG^{comp}\right)$ is in the image of $H_d^{simp}\left(DCone\left(\cup_{i=1}^s B\Gamma_i\rightarrow BG\right)\right)$ where the cone point of $Cone\left(B\Gamma_i\right)$ is mapped to $c_i$, the cusp associated to $\Gamma_i$ as defined after \hyperref[subgroup]{Definition \ref*{subgroup}}.

{\bf Step 1:} 
%Since $\Gamma_i$, as a group of parabolic isometries with a common fixed point, is unipotent, then 
%also the image of $\Gamma_i^\prime$ in 
%$SL\left(N,{\Bbb C}\right)$ is unipotent
%and therefore there is some $c_i^\prime\in\partial_\infty\left(SL\left(N,{\Bbb C}\right)/SU\left(N\right)\right)$ with $\Gamma_i^\prime\subset Fix\left(c_i^\prime\right)$. Then $B\left(\rho j\right)_dEM_d^{-1}\left[M,\partial M\right]\in H_d^{simp}\left(BSL\left(N,{\Bbb F}\right)^{fb}
%\right)$ is in the image of $H_d^{simp}\left(DCone\left(\cup_{i=1}^s B\Gamma_i^\prime\rightarrow BSL\left(N,{\Bbb F}\right)\right)\right)$, where  
%the cone point of $Cone\left(B\Gamma_i^\prime
%\right)$ is mapped to $c_i^\prime$.
%
%We are going to prove that (each) preimage of $  B\left(\rho j\right)_dEM_d^{-1}\left[M,\partial
%M\right]$ in $H_d^{simp}\left(DCone\left(\cup_{i=1}^sB\Gamma_i^\prime\rightarrow BSL\left(N,{\Bbb F}\right)
%\right);{\Bbb Q}\right)$
%has a preimage $\overline{\gamma}\left(M\right)\in H_d^{simp}\left(BSL\left(N,{\Bbb F}\right);{\Bbb Q}\right).$
Assume that $B\left(\rho j\right)_{*}EM^{-1}\left[\partial_i M\right]=0\in H_{d-1}^{simp}\left(
%BSL\left({\Bbb F}\right)
B\Gamma_i^\prime;
{\Bbb Q}\right)$ for $i=1,\ldots,s$ and consider the commutatve diagram
$$ \begin{xy}
\xymatrix{
Z_d\left(M,\partial M\right)\ar[r]^{c} \ar[d]^{\partial_i}&\widehat{Z}_d\left(M\right)\ar[r]^{\Phi\circ str}&
Z_d^{simp}\left(B\Gamma^{comp}\right)\ar[r]^{Bj_*}& Z_d^{simp}\left(BG\left({\Bbb F}\right)^{comp}\right)\ar[r]^{B\rho_*}&
Z_d^{simp}\left(BSL\left(N,{\Bbb F}\right)^{fb}\right)\\
Z_{d-1}\left(\partial_iM\right)\ar[r]^=&Z_{d-1}\left(\partial_iM\right)\ar[r]^{EM^{-1}}\ar[u]^{Cone}&
Z_{d-1}^{simp}\left(B\Gamma_i\right)\ar[r]^=\ar[u]^{Cone}&
Z_{d-1}^{simp}\left(B\Gamma_i\right)\ar[r]^{B\left(\rho j\right)_{*}}\ar[u]^{Cone}& Z_{d-1}^{simp}\left(B\Gamma_i^\prime\right)\ar[u]^{Cone}}
\end{xy}$$
where $Z_d\left(M,\partial M\right)\subset C_d\left(M,\partial M\right)$ is the subgroup of
relative cycles, and for
a relative cycle $z$ we define $c\left(z\right)=z+Cone\left(\partial z\right)\in\widehat{C}_d\left(M\right)$ and $\partial_iz$ to be the image of $\partial z\in C_{d-1}\left(\partial M\right)$ under the projection from $C_{d-1}\left(\partial M\right)$ to its direct summand 
$C_{d-1}\left(\partial_i M\right)$.

If $z\in C_d\left(M,\partial M\right)$
%C_d^{simp}\left(DCone\left(\cup_{i=1}^s B\Gamma_i\rightarrow B\Gamma\right);{\Bbb Q}\right)$ 
is a relative cycle that represents $\left[M,\partial M\right]$, then $\partial_iz$ represents $\left[\partial_iM\right]$.
%If for $i=1,\ldots,s$ 
By the assumption of Step 1 we have chains $z_i^\prime\in C_d^{simp}\left(
%BSL\left({\Bbb F}\right)
B\Gamma_i^\prime\right)$ with
$$\partial z_i^\prime=B
\left(\rho j\right)_{*}EM^{-1}\left(\partial_i z \right),$$
%then we let $z_i:=\left(EM_d\right)^{-1}\left(z_i^\prime\right)$ and conclude 
then $$B\left(\rho j \right)_*\Phi\left(str\left(c \left(z\right)\right)\right)- \sum_{i=1}^sCone\left(z_i^\prime\right)\in Z_d^{simp}\left(BSL\left(N,\ff\right)^{fb}\right)$$ 
is a genuine cycle in $Z_d^{simp}\left(BSL\left(N,\ff\right)\right)$, whose image in $Z_d\left(
BSL\left(N,\ff\right)^{fb}
%DCone\left(\cup_{i=1}^sB\Gamma_i^\prime\rightarrow BSL\left(N,{\Bbb F}\right)
\right)$ again represents $B\left(\rho j\right)_*EM^{-1}\left[M,\partial
M\right]$. Therefore its homology class
%$B\left(\rho j\right)_d\Phi\left(str\left(c \left(z\right)\right)\right)- \sum_{i=1}^sCone\left(z_i^\prime\right)$ represents 
gives the desired $\overline{\gamma}\left(M\right)$. 

{\bf Step 2:} It remains to 
prove $B\left(\rho j\right)_{*}EM^{-1}\left[\partial_i M\right]=0$. Let $f_i:\partial_i M\rightarrow M$ be the inclusion, $q:M\rightarrow M_+$ the projection.
Thus $qf_i$ is constant.
% and the induced map $\left(qf\right)_*$ in homology
%is trivial. 
Recall that $\Gamma_i\subset G$ consists of parabolic isometries with the same fixed point in $\partial_\infty G/K$ (see \cite[Theorem 3.1]{ebe}), thus
$\Gamma_i$ and hence
$\Gamma_i^\prime:=\rho\left(\Gamma_i\right)$ are unipotent and we 
%Since $\Gamma_i^\prime=\rho\left(\Gamma_i\right)$ is unipotent, we 
can apply \hyperref[mapR]{Lemma \ref*{mapR}} and obtain a continuous map $$R:
M_+\rightarrow \mid BSL\left(N,\ff\right)\mid^+$$ such that $$R\circ q\circ f_i=incl\circ \mid B\left(\rho j\right)\mid\circ h^M\circ f_i.$$
In particular, $
incl\circ \mid B\left(\rho j\right)\mid\circ h^{\partial_iM}=
incl\circ \mid B\left(\rho j\right)\mid\circ h^M\circ f_i
%=incl\circ \mid B\left(\rho j\right)\mid\circ h^{\partial_iM}
:\partial_iM\rightarrow \mid BSL\left(N,\ff\right)\mid^+$
is constant. 

Since $incl\circ\mid B\left(\rho j\right)\mid:\mid B\Gamma_i\mid\rightarrow \mid BSL\left(N,\ff\right)\mid^+$
factors over $\mid B\Gamma_i^\prime\mid^+$ and since $$\mid
B\Gamma_i^\prime\mid\subset \mid BSL\left(N,\ff\right)\mid\subset
\mid BSL\left(N,\ff\right)\mid^+$$ are inclusions (the first by $\Gamma_i^\prime\subset SL\left(N,\ff\right)$, the second by 
%$BSL\left(N,{\Bbb F}\right)\subset BSL\left({\Bbb F}\right)$ and by the
the definition of the plus construction via attaching cells to $\mid BSL\left(N,\ff\right)\mid$), this implies that 
$$incl\circ \mid B\left(\rho j\right)\mid\circ h^{\partial_iM}:\partial_iM\rightarrow \mid B\Gamma_i^\prime\mid^+$$ is constant.
Since $incl_*:H_*\left(\mid B\Gamma_i^\prime\mid;{\Bbb Q}\right)\rightarrow
H_*\left(\mid B\Gamma_i^\prime\mid^+;{\Bbb Q}\right)$ is 
an isomorphism, this implies $$\mid B\left(\rho j\right)\mid_* h^{\partial_iM}_*=0,$$
in particular $$
%B\left(\rho j\right)_{d-1}\left[\partial_i M\right]=
\mid B\left(\rho j\right)\mid_{*} h^{\partial_iM}_{*}\left[\partial_iM\right]=
%\mid B\left(\rho j\right)\mid_{d-1} h^M_{d-1}\left(f_i\right)_{d-1}\left[\partial_iM\right]=
0\in H_{d-1}\left(\mid B\Gamma_i^\prime\mid;
{\Bbb Q}\right)$$ for $i=1,\ldots,s$. 

But
$h^{\partial_iM}_{*}\left[\partial_iM\right]$ is the image of $EM^{-1}\left[\partial_iM\right]$ under the isomorphism $H_{d-1}^{simp}\left(B\Gamma_i;{\Bbb Q}\right)\rightarrow H_{d-1}\left(\mid B\Gamma_i\mid;{\Bbb Q}\right)$ (see Section 2.1), hence $\mid B\left(\rho j\right)\mid_{*} h^{\partial_iM}_{*}\left[\partial_iM\right]$ is the image of 
$$B\left(\rho j\right)_{*}EM^{-1}\left[
\partial_iM\right]$$ under the isomorphism $H_{d-1}^{simp}\left(B\Gamma_i^\prime;{\Bbb Q}\right)\rightarrow H_{d-1}\left(\mid B\Gamma_i^\prime\mid;{\Bbb Q}\right)$. Thus 
$$B\left(\rho j\right)_{*}EM^{-1}\left[
\partial_iM\right]=0.$$

{\bf Welldefinedness:} The construction of $\overline{\gamma}\left(M\right)$ as the homology class represented by
$B\left(\rho j \right)_*\Phi\left(str\left(c \left(z\right)\right)\right)- \sum_{i=1}^sCone\left(z_i^\prime\right)$ involves a choice of chains $z_i^\prime\in C_d^{simp}\left(
B\Gamma_i^\prime\right)$ with
$\partial z_i^\prime=B
\left(\rho j\right)_{*}EM^{-1}\left(\partial_i z \right)$. 

If $z_i^\prime$
and $z_i^{\prime\prime}$ are two such choices,
then $\partial z_i^\prime=\partial z_i^{\prime\prime}$ implies $z_i^\prime-z_i^{\prime\prime}\in Z_d^{simp}\left(B\Gamma_i^
\prime\right)$. Now we have $\Gamma_i=\pi_1\partial_iM$ and $\partial_iM$
is an aspherical $\left(d-1\right)$-manifold, hence 
$cd_{\Bbb Q}\left(\Gamma_i\right)=d-1$. 

We claim that this implies $cd_{\Bbb Q}\left(\Gamma_i^\prime\right)\le d-1$.
By a theorem of Gruenberg one has $cd_{\Bbb Q}\left(\Gamma_i\right)=h\left(\Gamma_i\right)$ for a finitely generated torsion-free nilpotent group
$\Gamma_i$, where $h$ denotes the Hirsch length and $cd_{\Bbb Q}$ the rational
cohomological dimension.
Also $\Gamma_i^\prime=\rho\left(\Gamma_i\right)$ is nilpotent, finitely generated and obviously
$h\left(\Gamma_i^\prime\right)\le h\left(\Gamma_i\right)$. Let $N\subset G$ be the maximal nilpotent (in the sense of the Iwasawa decomposition) group containing
$\Gamma_i$, then $\rho\left(N\right)$ is conjugate into the group 
of upper triangular matrices, in particular it is torsionfree. Thus 
$\Gamma_i^\prime$ is a finitely generated torsion-free nilpotent group, to which Gruenberg's Theorem applies, and we obtain $$cd_{\Bbb Q}\left(\Gamma_i^\prime\right)=h\left(\Gamma_i^\prime\right)\le h\left(\Gamma_i\right)=d-1.$$
%$h\left(\Gamma_i\right)=cd_{\Bbb Q}\left(\Gamma_i\right)=d-1$, where $h$ denotes the Hirsch length and $cd_{\Bbb Q}$ the rational
%cohomological dimension. (Here we use that $\Gamma_i$ is nilpotent because $G/K$ has rank one. A straightforward inductive application of 
%the Hochschild-Serre spectral sequence shows $cd_{\Bbb Q}\le h$. For finitely generated torsion-free nilpotent groups, cohomological dimension and Hirsch length actually agree by a theorem of Gruenberg, hence $cd_{\Bbb Q}\left(\Gamma_i\right)=h\left(\Gamma_i\right)$.)
%Also $\Gamma_i^\prime=\rho\left(\Gamma_i\right)$ is nilpotent and obviously
%$h\left(\Gamma_i^\prime\right)\le h\left(\Gamma_i\right)$, hence $cd_{\Bbb Q}\left(\Gamma_i^\prime\right)\le d-1$
Hence the $d$-cycle $z_i^\prime-z_i^{\prime\prime}$ must be 0-homologous in $C_*^{simp}\left(
B\Gamma_i^\prime\right)$. This implies that $Cone\left(z_i^\prime\right)-Cone\left(z_i^{\prime\prime}\right)$ is 0-homologous in
$C_*^{simp}\left(BSL\left(N,{\Bbb F}\right)^{fb}\right)$ and the homology class of $\overline{\gamma}\left(M\right)$ does not depend on the choice of $z_i^\prime$.

Moreover the construction of $\overline{\gamma}\left(M\right)$
involves a choice of a relative cycle
$z\in Z_d\left(M,\partial M\right)$. 

If $z$ and $z^\prime$ are two relative cycles representing $\left[M,\partial M\right]$,
then $z-z^\prime=\partial w+u$ for some $w\in C_{d-1}\left(M\right),u\in Z_d\left(\partial M\right)$. Again
$H_d\left(\Gamma_i^\prime;{\Bbb Q}\right)=0$ implies that $B
\left(\rho j\right)_*EM^{-1}\left(u\right)$ is a boundary. Hence $B
\left(\rho j\right)_*\Phi\left(str\left(c\left(z-z^\prime\right)\right)\right)$ is a boundary and thus
the homology class $\overline{\gamma}\left(M\right)$ does not
depend on the choice of $z$.

{\bf Independence of $\pi_1\partial_iM \cong \Gamma_i$:} The identification of $\pi_1\partial_iM$ with a
subgroup of $\Gamma$ depends on a path $\tilde{l}_i:\left[0,1\right]\rightarrow \widetilde{M}$. 
For two different paths one obtains subgroups which are conjugate in $\Gamma$. 
Conjugation in $\Gamma$ induces the identity homomorphism in group homology, thus 
the image of $H_*\left(B\Gamma_i\right)$ in 
$H_*\left(B\Gamma\right)$ and the image of 
$H_*\left(B\Gamma_i^\prime\right)$ in
$H_*\left(BSL\left(N,{\Bbb F}\right)^{fb}\right)$ 
do not depend on the chosen identification. In particular $\overline{\gamma}\left(M\right)$ does not depend.
\end{pf}\\
\\
For the special case of hyperbolic manifolds and half-spinor representations, Proposition 2 was proved in \cite[Theorem 2.12]{gon}. The proof in \cite{gon} uses very special properties of 
the half-spinor representations and seems not to generalize to other representations.

\subsection{Evaluation of Borel classes}
In \hyperref[Thm2]{Theorem \ref*{Thm2}} we proved for closed manifolds the equality $<b_{2n-1},\gamma\left(M\right)>=c_\rho vol\left(M\right)$. 
The following theorem will prove the analogous result for the cusped case.
 \begin{thm}\label{Thm3} 
a)  Let
$M$ be a compact, oriented, connected $2n-1$-manifold with boundary components
$\partial_1M,\ldots,\partial_sM$ such that
 $Int\left(M\right)$
is a locally
        symmetric space of noncompact type $Int\left(M\right)=\Gamma\backslash
         G/K$ of rank one with finite volume.
                Let $\rho:\left(G,K\right)\rightarrow \left(SL\left(N,{\Bbb C}\right),SU\left(N\right)\right)$ be a representation and
%such
%that $\rho\left(\pi_1\partial_i M\right) $ is unipotent for $i=1,\ldots, s$.
let $c_\rho$ be defined by \hyperref[Thm2]{Theorem \ref*{Thm2}}. Let 
$$\overline{\gamma}\left(M\right)\in H_{2n-1}\left(BSL\left(N,\overline{\Bbb Q}\right),{\Bbb Q}\right)$$ 
be defined by \hyperref[preimage]{Proposition \ref*{preimage}}, let $\overline{\overline{\gamma}}\left(M\right)$ 
be the image of $\overline{\gamma}\left(M\right)$ in $H_{2n-1}\left(BGL\left(\overline{\Bbb Q}\right),{\Bbb Q}\right)$ and define 
$$\gamma\left(M\right):=pr_{2n-1}
\left(\overline{\overline{\gamma}}\left(M\right)\right)\in PH_{2n-1}\left(BGL\left(\overline{\Bbb Q}\right),{\Bbb Q}\right)
\cong K_{2n-1}\left(\overline{\Bbb Q}\right)\otimes{\Bbb Q},$$ where $pr_{2n-1}$ is defined in \hyperref[projq]{Corollary \ref*{projq}}. Then
$$<b_{2n-1},\gamma\left(M\right)>=c_\rho vol\left(M\right).$$
b) If $\Gamma\subset G\left(A\right)$ for a subring $A\subset{\Bbb C}$ that satisfies the assumption of 
\hyperref[proj]{Lemma \ref*{proj}}
%, if we have an inclusion
%$j:\Gamma\rightarrow G\left(A\right)$ and 
and
if $\rho$ maps $G\left(A\right)$ to $SL\left(N,A\right)$\footnotemark\footnotetext[5]{For a semisimple Lie group $G$,
each representation $\rho:G\rightarrow SL\left(N,{\Bbb C}\right)$ is isomorphic to a representation which maps $G\left(A\right)$ to
$SL\left(N,A\right)$. (This can be read off the classification of representations of semisimple Lie groups, see \cite{fh}.)},
and if
$$\gamma\left(M\right):=pr_{2n-1}\left(\overline{\overline{\gamma}}\left(M\right)\right)\in 
PH_{2n-1}\left(BGL\left(A\right),{\Bbb Q}\right)\cong K_{2n-1}\left(A\right)\otimes{\Bbb Q},$$
where $\overline{\overline{\gamma}}\left(M\right)$
is the image of $\overline{\gamma}\left(M\right)\in H_{2n-1}\left(BSL\left(N,A\right),{\Bbb Q}\right)$ (defined
%$\overline{\gamma}\left(M\right)\in H_{2n-1}\left(BGL\left(A\right),{\Bbb Q}\right)$ is 
%defined 
by \hyperref[preimage]{Proposition \ref*{preimage}}) in $H_{2n-1}\left(BGL\left(A\right),{\Bbb Q}\right)$, and $pr_{2n-1}$ is given by \hyperref[proj]{Lemma \ref*{proj}}, then
$$<b_{2n-1},\gamma\left(M\right)>=c_\rho vol\left(M\right).$$
\end{thm}
        \begin{pf}
Denote $d=2n-1$.

$G$ is a linear semisimple Lie group without compact factors, not locally isomorphic to $SL\left(2,{\Bbb R}\right)$. By Weil rigidity we can assume (upon conjugation) that $\Gamma \subset G\left(\overline{\Bbb Q}\right)$. 
%For a semisimple Lie group $G$, each representation $\rho:G\rightarrow GL\left(N,{\Bbb C}\right)$ is isomorphic to a representation which maps $G\left(A\right)$ to 
%$SL\left(N,A\right)$. (This can be read off the classification of representations of semisimple Lie groups, see \cite{fh}.)
By \hyperref[projq]{Corollary \ref*{projq}}, $A=\overline{\Bbb Q}$ satisfies the assumptions of \hyperref[proj]{Lemma \ref*{proj}}. 
Thus a) is a consequence of b). We are going to prove b).

Let $z\in C_d\left(M,\partial M\right)$ represent $\left[M,\partial M\right]$. Then $\partial z\in C_{d-1}\left(\partial M\right)$ and
$$z+Cone\left(\partial z\right)\in \widehat{C}_d\left(M\right)\subset C_d\left(DCone\left(\cup_{i=1}^s\partial_iM\rightarrow M\right)\right)$$
represents the fundamental class. 
%Let $$\sum_{i=1}^r\tau_i$$ be some triangulation of $\left(M,\partial M\right)$. Let $$\sum_{j=1}^p\kappa_j:=DCone\left(\partial\left(\sum_{i=1}^r \tau_i\right)\right)$$ be the cone over the induced triangulations of $\partial_1 M,\ldots,\partial_sM$, with one cone point for each path-component of $\partial M$. Thus
%$$\sum_{i=1}^r\tau_i+\sum_{j=1}^p\kappa_j$$
%is a triangulation of $DCone\left(\cup_{i=1}^s\partial_i M\rightarrow M\right)$.

From \hyperref[cusps]{Corollary \ref*{cusps}} we get a homeomorphism $$DCone\left(\cup_{l=1}^s\partial_i M\rightarrow M\right)\cong \Gamma\backslash G/K \cup \left\{c_1,\ldots,c_s\right\},$$
where $c_l$ 
%are points ("cusp points")
corresponds to the cone point of $Cone\left(\partial_lM\right)$ for $l=1,\ldots,s$. Thus we can define $algvol\left(\sigma\right)=\int_\sigma dvol$ for $\sigma\in C_*\left(DCone\left(\cup_{i=1}^s\partial_i M\rightarrow M\right)\right)$, where $dvol$ is the volume form for the locally symmetric metric and the cusps $c_l$ are declared to have measure zero.
%evaluate the volume form of the locally symmetric metric on $z+Cone\left(\partial z\right)$, declaring the cone points to have measure zero.
%$\partial_1 M,\ldots,\partial_s M$ of $\partial M$. 
%Note that the inclusion $\Gamma\subset G$ maps (the isomorphic image of)
%each $\Gamma_i$ to a parabolic subgroup of $G$.

By Stokes' Theorem, evaluation of the volume form on $z+Cone\left(\partial z\right)$ does not depend on the chosen representative $z$ of $\left[M,\partial M\right]$. In particular we can, by Whitehead's Theorem, assume that $z$ is given by a triangulation of $\left(M,\partial M\right)$. Then $z+Cone\left(\partial z\right)$ is an ideal triangulation of $M$ and evaluation of the volume form gives the sum of the signed volumes of simplices in that triangulation, that is $vol\left(M\right)$. This shows that
$$algvol\left(z+Cone\left(\partial z\right)\right)=vol\left(M\right).$$

%$\sum_{i=1}^r\tau_i+\sum_{j=1}^p\kappa_j$ is an (ideal) triangulation of $\Gamma\backslash G/K \cup \left\{C_1,\ldots,C_s\right\}$, therefore $$vol\left(M\right)=vol\left(\Gamma\backslash G/K\right)=\sum_{i=1}^ralgvol\left(\tau_i\right)+\sum_{j=1}^palgvol\left(\kappa_j\right),$$ where algvol is (as in Section 2.3) the signed volume computed by integration of the volume form (and defining the volume of $\left\{C_1,\ldots,C_s\right\}$ to be zero). 
Let $x_0,x_i,\Gamma,\Gamma_i$ be defined according to \hyperref[subgroup]{Definition \ref*{subgroup}}.
Let $$str:\widehat{C}_*\left(M\right)\rightarrow \widehat{C}_*^{str,x_0}\left(M\right)$$
be the chain homotopy inverse of the inclusion given by part b) of \hyperref[iso]{Lemma \ref*{iso}}.

Then $str\left(z+Cone\left(\partial z\right)\right)$ is homologous to $z+Cone\left(\partial z\right)$, thus Stokes' Theorem implies $$algvol\left(str\left(z+Cone\left(\partial z\right)\right)\right)=algvol\left(z+Cone\left(\partial z\right)\right)=vol\left(M\right).$$

Let $$z+Cone\left(\partial z\right)=\sum_{i=1}^ra_i\tau_i+\sum_{j=1}^p b_j\kappa_j$$
with $\tau_i\in C_*\left(M\right)$ and $\kappa_j\in \cup_{l=1}^s Cone\left(C_*\left(\partial_lM\right)\right)$ for $i=1,\ldots,r, j=1, \ldots,p$.

Let $w_0,\ldots,w_d$ be the vertices of the standard simplex $\Delta^d$. By the proof
of \hyperref[iso]{Lemma \ref*{iso}}, the isomorphism
$$\Phi:\widehat{C}_*^{str,x_0}
\left(M\right)\rightarrow C_*^{simp}\left(
%DCone\left(\cup_{l=1}^s B\Gamma_l\rightarrow 
B\Gamma^{comp}\right)$$
maps the
interior simplex $str\left(\tau_i\right)$ to $$\left(\gamma_1^i,\ldots,\gamma_d^i\right)\in B\Gamma,$$
where $\gamma_k^i\in\Gamma$ is the homotopy class of the (closed) edge from
$\tau_i\left(w_{k-1}\right)$ to $\tau_i\left(w_k\right)$, \\
and 
the ideal simplex $str\left(\kappa_j\right)$ 
to $$\left(p_1^j,\ldots,p_{d-1}^j,c_{l_j}\right)\in Cone\left(B\Gamma_{l_j}\rightarrow B\Gamma\right),$$
where $\kappa_j\in Cone\left(C_*\left(\partial_{l_j}M\right)\right)$ and 
$c_{l_j}\in\partial_\infty G/K$ is the cusp associated to $\Gamma_{l_j}$ (cf.\ the remark after \hyperref[subgroup]{Definition \ref*{subgroup}}) and $p_k^j$ is the homotopy class of the (closed) edge from
$\kappa_j\left(w_{k-1}\right)$ to $\kappa_j\left(w_k\right)$. 
Thus, in the setting of \hyperref[preimage]{Proposition \ref*{preimage}}, we have that
$$\left(Bj\right)_*EM^{-1}\left[M,\partial M\right]
\in H_d\left(BG\left(A\right)^{comp};{\Bbb Q}\right)$$ is represented by $$\sum_{i=1}^r \left(\gamma_1^i,\ldots,\gamma_d^i\right) +\sum_{j=1}^p \left(p_1^j,\ldots,p_{d-1}^j,c_{l_j}\right).$$
Let $str\left(
\tilde{x},\gamma_1^i\tilde{x},\ldots,\gamma_1^i\ldots\gamma_d^i\tilde{x}\right)$
%\in C_*\left( G/K\right)$ 
be the unique straight simplex with vertices $\tilde{x},\gamma_1^i\tilde{x},\ldots,\gamma_1^i\ldots\gamma_d^i\tilde{x}$, and 
$str\left(\tilde{x},p_1^j\tilde{x},\ldots,p_1^j\ldots p_{d-1}^j\tilde{x},c_{l_j}\right)$ the unique ideal straight simplex with interior vertices $\tilde{x}, p_1^j\tilde{x},\ldots,p_1^j\ldots p_{d-1}^j\tilde{x}$ and ideal vertex $c_{l_j}$.

By construction we have $$
%\pi p\left(1,\gamma_1^i,\ldots,\gamma_d^i\right)=
\overline{\pi}\left(str\left(\tilde{x},\gamma_1^i\tilde{x},\ldots,\gamma_1^i\ldots\gamma_d^i\tilde{x}\right)\right)=str\left(\tau_i\right), $$
$$
%\pi p\left(1,p_1^j,\ldots,p_{d-1}^j,c^j\right)=
\overline{\pi}\left(str\left(\tilde{x},p_1^j\tilde{x},\ldots,p_1^j\ldots p_{d-1}^j\tilde{x},c_{l_j}\right)\right)=str\left(\kappa_j\right).$$
Hence 
%(since $\pi^*$ maps the volume form of $M$ to the volume form of $G/K$)
$$\int_{str\left(
\tilde{x},\gamma_1^i\tilde{x},\ldots,\gamma_1^i\ldots \gamma_d^i\tilde{x}\right)}dvol_{G/K}=\int_{str\left(
\tilde{x},\gamma_1^i\tilde{x},\ldots,\gamma_1^i\ldots \gamma_d^i\tilde{x}\right)}\overline{\pi}^*dvol_M=
\int_{str\left(\tau_i\right)}dvol_M=algvol\left(str\left(\tau_i\right)\right),$$
$$ 
\int_{str\left(\tilde{x},p_1^j\tilde{x},\ldots,p_1^j\ldots p_{d-1}^j\tilde{x},c_{l_j}\right)}
dvol_{G/K}
=\int_{str\left(\tilde{x},p_1^j\tilde{x},\ldots,p_1^j\ldots p_{d-1}^j\tilde{x},c_{l_j}\right)}
\overline{\pi}^*dvol_M=\int_{str\left(\kappa_j\right)}dvol_M=algvol\left(str\left(\kappa_j\right)\right).$$
%is the volume of the straight simplex with vertices 
%$\tilde{x}_0,\gamma_1^i\tilde{x}_0,\ldots,\gamma_n^i\tilde{x}_0$, 
%i.e.\ of the lift of $\tau_i$ to $\widetilde{M}$ with first vertex 
%$\tilde{x}_0$. Hence $<v_n,\left(1,\gamma_1^i,\ldots,\gamma_n^i\right)>
%=vol\left(\tau_i\right)$, 
By the construction of the volume cocycle $\overline{c\nu}_d$ in Section 4.2.3 this implies
$$\overline{c\nu}_d\left(\sum_{i=1}^r a_i\left(1,\gamma_1^i,\ldots,\gamma_d^i\right)+\sum_{j=1}^p 
b_j\left(1,p_1^j,\ldots,p_{d-1}^j,c_{l_j}\right)\right) =$$
$$\sum_{i=1}^r a_i algvol\left(str\left(\tau_i\right)\right)+\sum_{j=1}^p b_j
algvol\left(str
\left(\kappa_j\right)\right)=algvol\left(z+Cone\left(\partial z\right)\right)=vol\left(M\right).$$

%$\Gamma_l\subset G$ consists of parabolic isometries with the same fixed point in $\partial_\infty G/K$ (see \cite[Theorem 3.1]{ebe}), thus 
%$\rho\left(\Gamma_l\right)\subset SL\left(N,{\Bbb C}\right)$ is unipotent.
By \hyperref[fb]{Lemma \ref*{fb}}b) 
%we have 
%$B\left(\rho j\right)_dEM_d^{-1}\left[M,\partial M\right]\in H_*^{simp}\left(BSL\left(N,{\Bbb C}\right)^{fb}\right)$, which by 
and \hyperref[subfield]{Definition \ref*{subfield}} 
%implies
we have
$B\left(\rho j\right)_*EM^{-1}\left[M,\partial M\right]\in H_*^{simp}\left(BSL\left(N,A\right)^{fb};{\Bbb Q}\right)$.
%\subset C_*^{simp}\left(BSL\left(N,{\Bbb C}\right)^{comp}\right).$$

%By \hyperref[preimage]{Proposition \ref*{preimage}}, 
%we have
% $$ i_*\overline{\gamma}\left(M\right)=B\left(\rho j\right)_dEM_d^{-1}\left[M,\partial M\right],$$
%where $i:BSL\left(N,{\Bbb C}\right)\rightarrow BSL\left(N,{\Bbb C}\right)^{fb}$ is the inclusion.
%\overline{\gamma}\left(M\right)\in H_d^{simp}\left(BSL\left(N,A\right); {\Bbb Q}\right),$$
%in $H_d^{simp}\left(DCone\left(\cup_{i=1}^s B\Gamma_i\rightarrow BSL\left(N,A
%\right)\right);{\Bbb Q}\right)$ equals $B\left(\rho\right)_dj_dEM_d\left[M,\partial M\right]$. 
%(Note that $G$ is perfect, thus $\rho$ has image in $SL\left(N,{\Bbb C}\right)$.)
%and is therefore represented by 
%$$\sum_{i=1}^r  \left(1,\rho\left(\gamma_1^i\right),\ldots,\rho\left(\gamma_n^i\right)\right) +\sum_{j=1}^p \left(1,\rho\left(p_1^j\right),\ldots,\rho\left(p_{n-1}^j\right),\rho\left(c^j\right)\right).$$ 

By \hyperref[fb]{Lemma \ref*{fb}}, there is 
$\overline{c\beta}_d: 
C_d^{simp}\left(BSL\left(N,{\Bbb C}\right)^{fb};{\Bbb R}\right) 
\rightarrow {\Bbb R}$
%$$\overline{b}_n: im\left(H_n\left(BSL\left(N,{\Bbb C}\right)^{fb};{\Bbb R}\right)\rightarrow H_n\left(BSL\left(N,{\Bbb C}
%\right)^{comp};{\Bbb R}\right)\right)\rightarrow {\Bbb R}$$
such that $ \overline{c\beta}_d\mid_{C_d^{simp}\left(BSL\left(N,{\Bbb C}\right);{\Bbb R}\right)}$ represents $comp\left(b_d
\right)$ and 
$\rho^*\overline{c\beta}_d$ represents $c_\rho \overline{cv}_d.$ (In particular, $\overline{c\beta}_d$ is well-defined on $\left(B\rho\right)_*H_d^{simp}\left(BG\left(A\right)^{comp};{\Bbb Q}\right)$.)
Then we have 
$$\left[\overline{c\beta}_d\right]\left(B\left(\rho j\right)_*EM^{-1}\left[M,\partial M\right]\right)=
\rho^*\left[\overline{c\beta}_d\right]\left( \left(Bj\right)_*EM^{-1}\left[M,\partial M\right]\right) =c_\rho\overline{cv}_d\left(Bj_*EM^{-1}\left[M,\partial M\right]\right)$$
$$=c_\rho\overline{c\nu}_d\left(\sum_{i=1}^r \left(1,\gamma_1^i,\ldots,\gamma_d^i\right)+\sum_{j=1}^p
\left(1,p_1^j,\ldots,p_{d-1}^j,c^j\right)\right)=c_\rho vol\left(M\right).$$

Let $i:BSL\left(N,A\right)\rightarrow BSL\left(N,A\right)^{fb}$ be the inclusion, then 
\hyperref[preimage]{Proposition \ref*{preimage}} gives
$$i_*\overline{\gamma}
\left(M\right)=B\left(\rho j\right)_* EM^{-1}\left[M,\partial M\right].$$  
%From application
Applying \hyperref[fb]{Lemma \ref*{fb}}a) to $C_d^{simp}\left(BSL\left(N,A\right);{\Bbb R}\right)\subset C_d^{simp}\left(BSL\left(N,{\Bbb C}\right);{\Bbb R}\right)$ we obtain  
%and 
$$i^*\overline{c\beta}_d=comp\left(b_d\right).$$
%follows by restriction from \hyperref[fb]{Lemma \ref*{fb}}.
Thus, 
%with $i:BSL\left(N,{\Bbb C}\right)\rightarrow BSL\left(N,{\Bbb C}\right)^{comp}$ the inclusion, and 
confusing
$\overline{\gamma}\left(M\right)$
%\in H_d\left(BSL\left(N,A\right);{\Bbb Q}\right))$ 
with its image in 
$H_d\left(BSL\left(N,{\Bbb C}\right);{\Bbb R}\right))$ 
%and $ B\left(\rho j\right)_dEM_d\left[M,\partial
%M\right]\in H_d^{simp}\left(DCone\left(\cup_{i=1}^sB\Gamma_i^\prime\rightarrow BSL\left({\Bbb F}\right)
%$\right);{\Bbb Q}\right)$ with its image in $H_d^{simp}\left(BSL\left({\Bbb C}\right)^{comp};{\Bbb Q}\right)$, 
we have
$$<b_d,\overline{\gamma}\left(M\right)>=comp\left(b_d\right)\left(\overline{\gamma}\left(M\right)\right)=$$
$$\left[i^*\overline{c\beta}_d\right]\left(\overline{\gamma}\left(M\right)\right)=
\left[\overline{c\beta}_d\right]\left(i_*\overline{\gamma}\left(M\right)\right)=$$
%<i^*\overline{b}_d,\overline{\gamma}\left(M\right)>$$
%$$=<\overline{b}_d,i_d\overline{\gamma}\left(M\right)>=
$$\left[\overline{c\beta}_d\right]\left(B\left(\rho j\right)_*EM^{-1}\left[M,\partial M\right]\right)
=c_\rho vol\left(M\right).$$
%By the van Est theorem there is an isomorphism $H_c^*\left(G\right)\cong H^*\left(\g,\kk\right)$ under which $v_n$ corresponds 
%to the volume form $dvol$, and $b_n$ corrsponds to some form which we will denote $dbol$. We have $\rho^*b_n=c\rho v_n$, hence $\rho^*dbol=c_\rho dvol$. $dvol$ abd $dbol$ pull back to differential forms on $BG^{comp}$. We have
%$$<b_n,\gamma\left(M\right)>=<i^* dbol,\gamma\left(M\right)>=<dbol,i_*\gamma\left(M\right)>=<dbol,\sum_{i=1}^ra_i\rho\left(\kappa_i\right)>$$
%$$=<\rho^* dbol, \sum_{i=1}^r a_i\kappa_i>=
%<c_\rho dvol, \sum_{i=1}^r a_i\kappa_i>
By \hyperref[proj]{Lemma \ref*{proj}} this implies $<b_{d},\gamma\left(M\right)>=c_\rho vol\left(M\right)$.

                \end{pf}

%Remark: The same argument would work for any relative fundamental cycle, not necessarily coming from a triangulation.

{\bf Examples.}
Cusped hyperbolic 3-manifolds were discussed to some extent in \cite{ny}.

If $M$ is any hyperbolic 3-manifold of finite volume, then $\pi_1M$
can be conjugated to a subgroup of $SL\left(2,{\Bbb F}\right)$, where ${\Bbb F}$ is an at most quadratic extension of the
trace field (\cite{mac}), thus
one gets an element in
$K_3\left({\Bbb F}\right)\otimes{\Bbb Q}$.
In \cite[Section 9]{ny} some examples of this construction
are given. (The discussion in \cite{ny} is about elements $\beta\left(M\right)\in
B\left({\Bbb F}\right)\otimes{\Bbb Q}$
for the Bloch group $B\left({\Bbb F}\right)$ but of course, using Suslin's isomorphism $B\left({\Bbb F}\right)\otimes{\bf Q}\cong K_3^{ind}\left({\Bbb F}
\right)\otimes {\bf Q}$ from \cite{su} and the isomorphism $K_3^{ind}\left({\Bbb F}\right)\otimes {\bf Q}
\cong K_3\left({\Bbb F}\right)\otimes{\bf Q}$ for number fields,
this construction yields elements in $K_3\left({\Bbb F}\right)\otimes {\Bbb Q}$ associated to the respective manifolds and it can 
actually be shown that $\gamma\left(M\right)$ corresponds to $\beta\left(M\right)$ under this isomorphism.)

For example (see \cite[Section 9.4]{ny}) for any number field ${\Bbb F}$
with just one complex place there exists a hyperbolic 3-manifold of finite volume, such that its
invariant trace field equals ${\Bbb F}$. The associated $\gamma\left(M\right)$ gives a nontrivial
element, and actually {\bf a generator}, in $K_3\left({\Bbb F}\right)\otimes{\Bbb Q}$.

\end{document}